\documentclass{ps}
\usepackage{amsmath}
\usepackage{amssymb}
\usepackage{accents}

\usepackage[colorinlistoftodos]{todonotes}
\usepackage{hyperref}
\numberwithin{equation}{section}

\newcommand{\EF}{{\mathrm{eff}}}
\newcommand{\FF}{{\mathrm{ff}}}
\newcommand{\PM}{{\mathrm{pm}}}
\newcommand{\TR}{{\mathrm{tr}}}
\newcommand{\WA}{{\mathrm{wall}}}
\newcommand{\DN}{{\mathrm{dn}}}
\newcommand{\JU}{{\mathrm{sj}}}
\newcommand{\BJ}{{\mathrm{BJ}}}
\newcommand{\dyn}{{\mathrm{dyn}}}
\newcommand{\OUT}{{\mathrm{outflow}}}

\renewcommand{\vec}[1]{{\ensuremath{\boldsymbol{\mathrm #1}}}}
\newcommand{\ten}[1]{\ensuremath{\boldsymbol{\mathsf{#1}}}}

\newcommand{\R}{\mathbb{R}}
\newcommand{\N}{\mathbb{N}}
\newcommand{\sspace}{\Omega}
\newcommand{\E}{\ensuremath{\mathbb{E}}}
\newcommand{\mean}{\mu}
\newcommand{\std}{\sigma}
\newcommand{\var}{\std^2}
\newcommand{\smean}[1]{\mean_{\surr{#1}}}

\newcommand{\svar}[1]{{\std^2_{\surr{#1}}}}
\newcommand{\aMR}{\text{aMR-PC}}
\newcommand{\Nr}{\ensuremath{N_r}}
\newcommand{\No}{\ensuremath{{N_o}}}
\newcommand{\Nc}{\ensuremath{N_{cf}}}
\newcommand{\Lp}[1]{\ensuremath{L^{#1}}}
\newcommand{\Prj}[1]{\ensuremath{P_{#1}}}  
\newcommand{\Pm}{\ensuremath{\mathcal{P}}}
\newcommand{\rf}{\ensuremath{Y}}
\newcommand{\rv}{\ensuremath{\theta}}
\newcommand{\vrv}{\vec{\rv}}
\newcommand{\rvi}[1]{{\rv_{#1}}}
\newcommand{\sdim}{\ensuremath{M}}
\newcommand{\dx}[1]{\ensuremath{\,\mathrm{d}\,{#1}}}
\newcommand{\cN}{\ensuremath{\mathcal{N}}}
\newcommand{\sobi}[1]{{S_{#1}}}
\newcommand{\sobit}[1]{S^T_{#1}}
\newcommand{\surr}[1]{\ensuremath{\mathcal{S}^{#1}}}
\newcommand{\pol}[1]{\varphi_{#1}}
\newcommand{\Pol}[1]{\Phi_{#1}}
\newcommand{\mrpol}[2][\Nr]{\varphi^{#1}_{#2}}
\newcommand{\mrPol}[2][\Nr]{\Phi^{#1}_{#2}}
\newcommand{\pcf}[2]{{#1}_{#2}}
\newcommand{\mrcf}[3][\Nr]{{#2}^{#1}_{#3}}
\newcommand{\sprod}[2]{\left\langle #1,#2\right\rangle}
\newcommand{\cdf}[1]{\ensuremath{F_{#1}}}
\newcommand{\quant}[1]{\ensuremath{Q_{#1}}}
\newcommand{\SD}{\text{SD}}
\newcommand{\idxSet}{\mathcal{I}^{\sdim}_{\Nr}}
\newcommand{\midx}[1]{\mathfrak{#1}}
\newcommand{\midxi}[2][l]{{\midx{#1}_{#2}}}
\newcommand{\subs}[2][\Nr]{\sspace_{#1,#2}}
\newcommand{\msubs}[2][\sdim]{\subs{#2}^{#1}}
\newcommand{\abs}[1]{{\left|#1\right|}}
\newcommand{\norm}[2]{\left\|#1\right\|_{#2}}
\newcommand{\pidxSet}[1][{}]{\Lambda^{\sdim,\No{#1}}}
\newcommand{\mrpSet}[1][{}]{\accentset{\circ}{\Lambda}^{\sdim,\No{#1}}}
\usepackage{caption}
\usepackage{subcaption}
\usepackage{lineno} 

\begin{document}
\title{Surrogate-assisted global sensitivity analysis of \\a hybrid-dimensional Stokes--Brinkman--Darcy model}
\thanks{The first and fourth authors thank the Deutsche Forschungsgemeinschaft (DFG, German Research Foundation) -- Project Number 490872182 
for the financial support.}\thanks{The second and third authors thank the Deutsche Forschungsgemeinschaft (DFG, German Research Foundation) -- Project Number 535321881 for the financial support.}
\thanks{The third and fourth authors thank the Deutsche Forschungsgemeinschaft (DFG, German Research Foundation) -- Project Number 327154368 -- SFB 1313 
for the financial support.}
%
\author{Linheng Ruan}\address{Institute of Applied Analysis and Numerical Simulation, University of Stuttgart, Stuttgart, Germany; \\\email{linheng.ruan@ians.uni-stuttgart.de \ \& \ iryna.rybak@ians.uni-stuttgart.de }}
\author{Ilja Kröker}\address{Institute for Modelling Hydraulic and Environmental
Systems, University of Stuttgart, Stuttgart, Germany; \\ \email{ilja.kroeker@iws.uni-stuttgart.de \ \& \ sergey.oladyshkin@iws.uni-stuttgart.de}}
\author{Sergey Oladyshkin}\sameaddress{2}
\author{Iryna Rybak}\sameaddress{1}
\date{\today}
\begin{abstract} 
Development of new multiscale mathematical models often entails considerable complexity and multiple undetermined parameters, typically arising from closure relations. To enable reliable simulations, one must quantify how uncertain physical parameters influence model predictions. We propose surrogate-assisted global sensitivity analysis that combines computational efficiency with a rigorous assessment of parameter influence. In this work, we analyze the recently proposed hybrid-dimensional Stokes--Brinkman--Darcy model, which describes  fluid flows in coupled free-flow and porous-medium systems with arbitrary flow directions at the fluid--porous interface. The model results from vertical averaging and contains several unknown parameters. We perform surrogate-assisted global sensitivity analysis using Sobol’ indices to investigate the sensitivity of the model to variations of physical parameters for two test cases: filtration and splitting flows.  However, constructing surrogates for higher-dimensional random fields requires either many training runs or sophisticated sampling strategies. To address this, we compare polynomial chaos surrogates, including sparse and multi-resolution representations, for their efficiency in global sensitivity analysis, using a predefined Sobol’ sequence of training samples. Across the tested cases, multi-resolution approach delivers the most accurate estimation of Sobol’ indices.
\end{abstract}
%
%

\maketitle
\section{Introduction}

Fluid flow and transport processes in coupled systems consisting of a free-flow region in contact with a porous medium arise frequently in science and engineering. Examples include surface water infiltration into soil, contaminant transport, industrial drying, blood flow in biological tissues, etc. Since resolving the detailed pore-scale geometry is computationally demanding, the free-flow region and the porous-medium domain are typically modeled as two distinct continua separated by a fluid--porous interface. This approach requires either suitable coupling conditions or a dimensionally reduced model at the interface, or alternatively a full-dimensional representation of the transition zone between the two flow domains (Fig.~\ref{fig1:intro}). Depending on the application and flow regime, different  models can be applied in the two subdomains. In this paper, we focus on the problem most extensively studied over the past 25 years, namely the Stokes equations in the free-flow region coupled with Darcy’s law in the porous medium.

Various coupling strategies for the Stokes--Darcy flow problems have been proposed in the literature~\cite{Goyeau_Lhuillier_etal_03, Angot_etal_17, Discacciati_Miglio_Quarteroni_02, Kang_Wang_2024, Lasseux_Valdes-Parada_Bottaro_2021, Lacis_Bagheri_17}. The classical set of interface conditions, which comprise mass conservation across the sharp interface, the balance of normal forces, and a variant of the Beavers--Joseph condition for the tangential velocity component, has been widely used, e.g.~\cite{Beavers_Joseph_67, Saffman, Discacciati_Quarteroni_09, Jaeger_Mikelic_09}. However, these conditions have been shown to lose accuracy for arbitrary flow directions at the fluid--porous interface~\cite{Strohbeck-Eggenweiler-Rybak-23, Eggenweiler_Rybak_20}. To enhance the classical Stokes--Darcy coupling, several generalized formulations have been  proposed in recent studies~\cite{Eggenweiler_Rybak_MMS20, Angot_etal_20, Angot_etal_17, Lacis_Bagheri_17, Sudhakar_21, Lasseux_Valdes-Parada_Bottaro_2021, zampogna2020effective,Ruan_Rybak_AMC,Ruan_Rybak_TIPM}. However, many of these theoretically derived coupling conditions have not been validated for arbitrary flow directions near the fluid–porous interface (e.g., \cite{Angot_etal_17, Angot_etal_20, Lacis_etal_20, Goyeau_Lhuillier_etal_03}), and some have not been validated at all, largely because certain model parameters remain uncertain and require quantification. 

Sensitivity analysis is a widely used tool to estimate the impact of single model parameters or combinations of parameters on the model output. The two commonly used types of sensitivity analysis are \emph{local}, which evaluates the variation of the model output in terms of the variation of the model parameters located in the immediate area of interests, and \emph{global} sensitivity analysis (GSA), such as the Sobol'~\cite{MR1052836,MR1823119,HOMMA19961} or the Shapley~\cite{shapley1953,owen2014} indices, which address the variation of the model response over the entire domain of input parameters. For a comprehensive review, we refer the reader to~\cite{daveiga_book_2021,sullivan_book}. The present work is based on global (variance-based) sensitivity, specifically Sobol’ indices \cite{MR1052836} to account for interactions across the entire parameter domain. However, the classical Monte Carlo (MC) based approach usually leads to very high computational costs in particular for complex numerical models. Therefore, the use of surrogate models can allow us to reduce the computational effort of GSA \cite{daveiga_book_2021}. Moreover, surrogates based on the orthonormal polynomial representations are attractive because Sobol’ indices can be obtained directly from the spectral coefficients. Specifically, as shown by Sudret \cite{SUDRET2008964} and extended to multi-element/multi-resolution settings in \cite{kroker2023global}, polynomial chaos expansion (PCE) and arbitrary multi-resolution polynomial chaos (\aMR) compute Sobol’ indices directly from the coefficients of orthonormal (piecewise) polynomials, thereby eliminating additional surrogate evaluations. The idea of using the orthogonal structure of the PCE to compute the Sobol' sensitivity indices out of polynomial coefficients was initially introduced by Sudret in \cite{SUDRET2008964} and later extended to various versions of the sparse PCE by several authors, e.g.~\cite{BLATMAN20101216,buerkner2022}. The specifics of using the piecewise polynomial surrogate models for GSA were addressed in~\cite{kroker2023global}, where the computation of the Sobol' sensitivity indices using the \aMR-expansion coefficients was introduced. PCE/\aMR{}-assisted GSA is accurate with appropriately chosen training designs \cite{BLATMAN20101216,SUDRET2008964,kroker2023global}. When training data can be generated on demand this is rarely restrictive. However, with given (non-optimal) datasets, PCE-based surrogates can exhibit oscillations (Runge phenomenon), leading to reduced accuracy. This problem can be addressed using  several  regularization concepts, e.g. \cite{MR2764550,buerkner2022}, at the price of significantly increased  computational costs, which become prohibitively high for larger random fields. In this study, we employ the training points generated by a Sobol' sequence \cite{sobol1967} and convert them to the assumed probability distributions based on the quasi-Monte-Carlo (QMC) approach \cite{sullivan_book}. This design provides a generic, well space-filling training sequence and enables a fair assessment of non-regularized PCE-based surrogates for GSA.

The objective of the present study is to investigate GSA for the hybrid-dimensional Stokes--Brinkman--Darcy model recently derived in~\cite{Ruan_Rybak_AMC, Ruan_Rybak_TIPM}. This model is able to accurately describe fluid flows in coupled systems with arbitrary directions to the fluid--porous interface, thereby extending the applicability of the coupled Stokes--Darcy model to a broader range of applications. The model involves five uncertain parameters and has been validated against several test cases in which the pore geometry is known and certain coefficients can be computed from pore-scale information. The validation is carried out through comparisons between pore-scale and macroscale simulation results. Here, we investigate the Stokes--Brinkman--Darcy model using surrogate-assisted global sensitivity analysis based on Sobol’ indices with respect to five uncertain physical parameters. To demonstrate applicability of the proposed framework, we consider two representative test cases: a filtration problem and a splitting flow problem.

The manuscript is organized as follows. In Section~\ref{sec:HD-SBD}, we introduce the flow system of interest, present the coupled mathematical model, specify the parameter ranges together with their probability distributions, and define two different test cases. In Section~\ref{sec:GSAonPols}, we outline the variance-based global sensitivity analysis, introduce the fundamentals of PCE/\aMR-based surrogate models, and present the PCE/\aMR-surrogate-assisted GSA. In Section~\ref{sec:numex}, we address the approximation and analysis of the physical models using PCE/\aMR{}-based surrogate models, introduced in Section~\ref{sec:GSAonPols}. The discussion and conclusions follow in Section~\ref{sec:conclusions}.

\section{Hybrid-dimensional Stokes--Brinkman--Darcy model}\label{sec:HD-SBD}

\subsection{Flow and geometry setting}
We consider single-phase, steady-state flows where the Reynolds number is sufficiently low such that inertial effects can be neglected. The fluid is assumed to be incompressible with constant viscosity and it fully saturates the entire flow domain. In the porous-medium region, the material is taken to be homogeneous and non-deformable.

\begin{figure}[h!]
    \centering
    \includegraphics[scale=1.1]{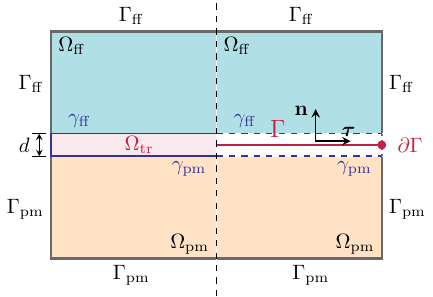}
    \caption{Schematic of coupled free-flow and porous-medium systems with a transition region $\Omega_\TR$ (left) and a complex interface $\Gamma$ (right)}
    \label{fig1:intro}
\end{figure}

In the coupled fluid--porous problem, the free-flow $\Omega_\FF$ and porous-medium $\Omega_\PM$ regions are two distinct continua with $\overline{\Omega}_\FF \cap \overline{\Omega}_\PM =  \emptyset$.  
To bridge these two subdomains,  we take into account an equi-dimensional transition zone $\Omega_\TR=\left\{ \vec{s} + t\vec{n} | \vec{s}\in \Gamma,\, t\in [-d/2, d/2] \right\}$ between them (Fig.~\ref{fig1:intro}, left), where $d>0$ is the thickness of the transition region, $\Gamma$ is the hyperplane and $\vec{n}$ is the unit normal vector pointing from hyperplane $\Gamma$ into the free-flow subdomain~$\Omega_\FF$. There are two different interfaces $\gamma_\FF=\overline{\Omega}_\TR \cap \overline{\Omega}_\FF$ and $\gamma_\PM=\overline{\Omega}_\TR \cap \overline{\Omega}_\PM$ on the top and bottom of the transition zone, respectively. 

In the hybrid-dimensional setting, the thickness of the transition zone is assumed to be sufficiently small compared to the size of the fluid domain $\overline{\Omega}_\mathrm{F}= \overline{\Omega}_\FF \cup \overline{\Omega}_\TR \cup \overline{\Omega}_\PM$. Therefore, we can treat it as a lower-dimensional inclusion and use the complex interface $\Gamma$ to approximate the transition region (Fig.~\ref{fig1:intro}, right). Note that the small gap is still existing between the free-flow region ($\overline{\Omega}_\FF \cap \overline{\Omega}_\PM =  \emptyset$). The mass and momentum are stored in and transferred along the complex interface $\Gamma$, which is modeled by averaging the corresponding conservation equations across the transition zone. 

\subsection{Model description}
In this section, we introduce the recently derived hybrid-dimensional Stokes--Brinkman--Darcy model~\cite{Ruan-Rybak-23, Ruan_Rybak_AMC}. The hybrid-dimensional model consists of the Stokes equations in the free-flow region, the averaged Brinkman equations on the complex interface $\Gamma$ and Darcy's law in the porous medium. In the free-flow domain, the laminar flow is modeled by the \textit{Stokes equations} involving the mass and momentum conservation 
\begin{align}
    \nabla \cdot \vec{v}_\FF &= 0 \hspace{+3.5ex} \textnormal{ in } \Omega_\FF,\label{eqn:Stokes1}\\ 
    -\nabla \cdot \ten{T} (\vec{v}_\FF, p_\FF) &= \vec{f}_\FF \quad \textnormal{ in } \Omega_\FF. \label{eqn:Stokes2}
\end{align}
Here, $\vec{v}_\FF$ and $p_\PM$ are the free-flow velocity and pressure, $\ten{T} (\vec{v},p) = \mu_\dyn \nabla \vec{v} - p\ten{I}$ is the pseudo stress tensor, $\mu_\dyn$ is the dynamic viscosity, $\ten{I}$ is the identity tensor, and $\vec{f}_\FF$ denotes the source term in the free flow. On the external boundary of the free-flow region $\Gamma_\FF=\Gamma_\FF^\mathrm{D}\cup\Gamma_\FF^\mathrm{N}$, either Dirichlet (D) or Neumann (N) conditions are considered
\begin{equation}
    \vec{v}_\FF = \overline{\vec{v}}_\FF \quad \textnormal{on }\Gamma_\FF^\mathrm{D}, \qquad \ten{T}(\vec{v}_\FF, p_\FF) \cdot \vec{n}_\FF = \overline{\vec{t}}_\FF \quad \textnormal{on }\Gamma_\FF^\mathrm{N},
\end{equation}
where $\overline{\vec{v}}_\FF$ and $\overline{\vec{t}}_\FF$ are the boundary values and $\vec{n}_\FF$ is the outer normal unit vector on the boundary $\Gamma_\FF^\mathrm{N}$.

In the porous-medium domain, the \textit{Darcy's law} is used to describe the slow flow
\begin{align}
    \nabla \cdot \vec{v}_\PM &= f_\PM \hspace{+11ex} \textnormal{ in } \Omega_\PM,\label{eqn:Darcy1}\\ 
    \vec{v}_\PM &= -\frac{\ten{K}_\PM}{\mu_\dyn} \nabla p_\PM \quad \textnormal{ in } \Omega_\PM, \label{eqn:Darcy2}
\end{align}
where $f_\PM$ is the source term, $\vec{v}_\PM$ is the porous-medium velocity and $p_\PM$ is the porous-medium pressure. $\ten{K}_\PM$ is the permeability tensor, which is symmetric positive definite.  On the external boundary of the porous medium $\Gamma_\PM = \Gamma_\FF^\mathrm{D} \cup \Gamma^\mathrm{N}_\FF$, the following conditions are taken into account
\begin{equation}
    \vec{v}_\PM \cdot \vec{n}_\PM=  \overline{v}_\PM \quad \textnormal{on }\Gamma_\PM^\mathrm{N}, \qquad p_\PM = \overline{p}_\PM \quad \textnormal{on } \Gamma_\PM^\mathrm{D},   
\end{equation}
where $\vec{n}_\PM$ is the outer normal unit vector on the boundary $\Gamma_\PM^\mathrm{N}$, $\overline{v}_\PM$ and $\overline{p}_\PM$ are the boundary data values.

In the full-dimensional setting, the Brinkman equations, which act as a bridge between the Stokes equations and Darcy's law, are considered in the transition region. In hybrid-dimensional model, the Brinkman equations are averaged in the vertical direction across the transition zone. The \textit{averaged Brinkman equations} on the complex interface yield
\begin{align}
    \vec{v}_\FF \cdot  \vec{n} |_{\gamma_\FF} - \vec{v}_\PM \cdot \vec{n} |_{\gamma_\PM}  &= - d \frac{\partial  V_{\vec{\tau}} }{\partial \vec{\tau}} \hspace{+34ex}\textnormal{ on } \Gamma, \label{eqn:AverageBrinkmanMass}\\    
    \left.\left(  \vec{n} \cdot  \ten{T}\left(\vec{v}_\FF ,p_\FF \right) \cdot  \vec{n} -\frac{\mu_\dyn }{\sqrt{K_{\Gamma}}} \left(\vec{\beta}\vec{v}_\FF \right)\cdot\vec{n} \right)\right|_{\gamma_\FF} + p_\PM \big|_{\gamma_\PM }&= d \left(  
    \mu_\dyn (\ten{K}_{\Gamma}^{-1} \vec{V})\cdot \vec{n}- \mu_{\EF} \frac{\partial^2 V_{\vec{n}}}{\partial \vec{\tau}^2}-F_{\vec{n}}\right) \hspace{+2.6ex}\textnormal{ on } \Gamma,\label{eqn:AveragedBrinkmanMomentumNormal} \\
    \left.\left(\vec{n} \cdot \ten{T}\left(\vec{v}_\FF, p_\FF \right) \cdot  \vec{\tau} -  \frac{\mu_\dyn }{\sqrt{K_{\Gamma}}} \left(\vec{\beta}\vec{v}_\FF \right) \cdot  \vec{\tau} \right)\right|_{\gamma_\FF}&-
        \frac{\alpha_\BJ \mu_\EF (6 V_{\vec{\tau}}- 2 \vec{v}_\FF \cdot  \vec{\tau}|_{\gamma_\FF})}{\alpha_\BJ d +4\sqrt{K_\PM}}\nonumber \\
        = d &\left( 
        \mu_\dyn (\ten{K}_{\Gamma}^{-1} \vec{V}) \cdot  \vec{\tau}
        -  \mu_\EF \frac{\partial^2 V_{\vec{\tau}}}{\partial \vec{\tau}^2} + \frac{\partial P}{\partial \vec{\tau}}-F_{\vec{\tau}} \right) \;\textnormal{ on } \Gamma,\label{eqn:AveragedBrinkmanMomentumTangential}
\end{align}
where $\vec{V}=(V_\vec{\tau}, V_\vec{n})^\top$ and $P$ denote the interfacial velocity and pressure obtained from the vertical averaging of the corresponding quantities across the transition zone, $\ten{K}_\Gamma$ is the permeability tensor on the complex interface~$\Gamma$, $\alpha_\BJ>0$ is the slip coefficient, $\vec{\beta}$ is the stress jump between $\Gamma$ and $\Omega_\FF$ and $F_\vec{n}$, $F_\vec{\tau}$ are the source terms.  Note that $\ten{K}_\Gamma$ is symmetric positive definite and has the same structure as the porous-medium permeability tensor~$\ten{K}$. The corresponding permeability values are $K_\Gamma:=\|\ten{K}_\Gamma\|_{\infty}$, $K_\PM:= \vec{\tau} \cdot \ten{K}_\PM \cdot \vec{\tau}$.

In the hybrid-dimensional Stokes--Brinkman--Darcy problem, additional transmission conditions are required to extrapolate the velocity and the pressure values on $\gamma_\FF$ and $\gamma_\PM$ and get a closed model. The following \textit{transmission conditions} based on the \textit{a-priori} assumptions on the velocity and pressure profiles across $\Gamma$ are taken into account
\begin{align}
     \left.\left( \vec{n} \cdot  \ten{T}(\vec{v}_\FF, p_\FF ) \cdot  \vec{n}-\frac{\mu_\dyn }{\sqrt{K_{\Gamma}}} \left(\vec{\beta}\vec{v}_\FF \right)\cdot\vec{n} \right)\right|_{\gamma_\FF} =& -\frac{\mu_\EF}{d}\left(6 V_{\vec{n}}  - 4 \vec{v}_\FF\cdot \vec{n} |_{\gamma_\FF} - 2 \vec{v}_\PM \cdot \vec{n} |_{\gamma_\PM}\right)-P \hspace{+4.3ex}\textnormal{ on } \Gamma,\label{eqn:TransmissionGammaFFNormal}
      \\
    \left.\left(\vec{n}\cdot \ten{T}(\vec{v}_\FF, p_\FF ) \cdot  \vec{\tau}-\frac{\mu_\dyn }{\sqrt{K_{\Gamma}}} \left(\vec{\beta}\vec{v}_\FF \right) \cdot  \vec{\tau} \right)\right|_{\gamma_\FF}=& - \frac{ 6\mu_\EF( \alpha_\BJ d + 2\sqrt{K_\PM}) }{d(\alpha_\BJ d+4\sqrt{K_\PM})} V_{\vec{\tau}} \nonumber\\
    & \hspace{+13.5ex}+ \frac{  4 \mu_\EF ( \alpha_\BJ d + 3\sqrt{K_\PM}) }{ d (\alpha_\BJ d+4\sqrt{K_\PM}) } \vec{v}_\FF \cdot  \vec{\tau}  \big|_{\gamma_\FF}  \hspace{+2.2ex}\textnormal{ on } \Gamma,\label{eqn:TransmissionGammaFFtangential}\\
    p_\PM \big|_{\gamma_\PM } =& -\frac{\mu_\EF}{d} \left(6V_{\vec{n}}-2\vec{v}_\FF \cdot \vec{n}|_{\gamma_\FF} -4\vec{v}_\PM \cdot \vec{n} |_{\gamma_\PM} \right) + P \hspace{+4.5ex}\textnormal{ on } \Gamma.\label{eqn:TransmissionGammaPM}
\end{align}
Note that we assume a quadratic velocity and a constant pressure profile across the transition regions in~\eqref{eqn:TransmissionGammaFFNormal}--\eqref{eqn:TransmissionGammaPM}. The use of a quadratic velocity profile was shown to be a good approximation with the validation against the pore-scale simulation results~\cite{Ruan_Rybak_TIPM}. 

On the external boundary nodes of the complex interface $\partial \Gamma = \partial \Gamma^\mathrm{D} \cup \partial \Gamma^\mathrm{N}$ (Fig.~\ref{fig1:intro}, right), the following boundary conditions are taken into account
\begin{equation}
    \vec{V} = \overline{\vec{V}} \hspace*{+1ex}\textnormal{ on }\partial \Gamma^\mathrm{D} , \quad    \mu_\EF \frac{\partial V_\vec{n}}{\partial \vec \tau}= \overline{T}_{\vec{n}}, \quad \mu_\EF \frac{\partial V_\vec{\tau}}{\partial \vec \tau}-P=\overline{T}_\vec{\tau} \hspace*{1ex}\textnormal{ on } \partial \Gamma^\mathrm{N}, \label{eqn:gammaBC}
\end{equation}
where $\overline{\vec{V}}$, $\overline{T}_\vec{n}$ and $\overline{T}_\vec{\tau}$ are the given boundary data.

Note that in the full-dimensional setting, there are two interfaces, $\gamma_\FF$ and $\gamma_\PM$, located at the top and bottom of the transition region (Fig.~\ref{fig1:intro}, left), respectively, where different interface conditions are imposed. On the interface $\gamma_\FF$ between the free flow and the transition region, we consider the continuity of velocity and the stress jump conditions proposed in~\cite{Angot_etal_17}. These conditions involve the stress jump tensor $\ten{\beta}$, the effective viscosity $\mu_\EF$, and the  permeability tensor $\ten{K}_\Gamma$. On the interface $\gamma_\PM$ between the transition region and the  porous medium, we impose mass conservation, balance of normal forces, and the Beavers--Joseph--Saffman condition~\cite{Saffman}. These coupling conditions contains the effective viscosity $\mu_\EF$, the Beavers–Joseph coefficient $\alpha_\BJ$, and the porous-medium permeability tensor $\ten{K}_\PM$. The coupling conditions~\eqref{eqn:AverageBrinkmanMass}–\eqref{eqn:TransmissionGammaPM} on the complex interface $\Gamma$ are derived taking the interface conditions on $\gamma_\FF$ and $\gamma_\PM$ into account. We refer the reader to~\cite{Ruan-Rybak-23,Ruan_Rybak_AMC} for detailed derivation.  

\subsection{Model parameters}
In this work, we make sensitivity analysis for the hybrid-dimensional Stokes--Brinkman--Darcy model~\eqref{eqn:Stokes1}--\eqref{eqn:gammaBC}. We restrict ourselves to isotropic porous media: $\ten{K}_\PM= k_\PM \ten{I}$, $\ten{K}_\Gamma=k_\Gamma \ten{I}$, $\ten{\beta}=\beta_\JU\ten{I}$ with $k_\Gamma \ge k_\PM >0$, $\beta_\JU\ge 0$. Therefore, we obtain $K_\Gamma= k_\Gamma$ and $K_\PM = k_\PM$ in Eqs.~\eqref{eqn:AveragedBrinkmanMomentumNormal}--\eqref{eqn:TransmissionGammaFFtangential}. The following parameters in the hybrid-dimensional model are uncertain: stress jump parameter $\beta_\JU$, permeability on the complex interface $k_\Gamma$, effective viscosity $\mu_\EF$,  slip parameter $\alpha_\BJ$, and permeability in the porous medium $k_\PM$. 
\begin{table}[h!]
\caption{Typical values and intervals for the uncertain parameters for the hybrid-dimensional model.\label{table1:UQtable}}
\begin{tabular}{p{0.3cm} p{4cm}  p{3cm}  p{3cm}} 
 \hline
 $I$ &Parameter & Interval & Distribution \\ [0.5ex] 
 \hline\hline
 1&Stress jump $\beta_\JU$ & $\left[ 0, 10\right]$ & uniform  \\
2&Permeability $k_{\Gamma}$ & $ ( 10^{-5}, 10^{-2} )$ & log-normal  \\
 3&Effective viscosity $\mu_\EF$ & $ \left[ 0.1, 10\right]$ & beta  \\ 
4&Beavers--Joseph $\alpha_\BJ$ & $( 0, 10)$ & beta  \\
5&Permeability $k_\PM$ & $( 10^{-8}, 10^{-5})$ & log-normal \\
 \hline
\end{tabular}
\end{table}
\begin{figure}[b]
    \centering
    \begin{subfigure}{0.32\linewidth}
    \includegraphics[height=3.7cm]{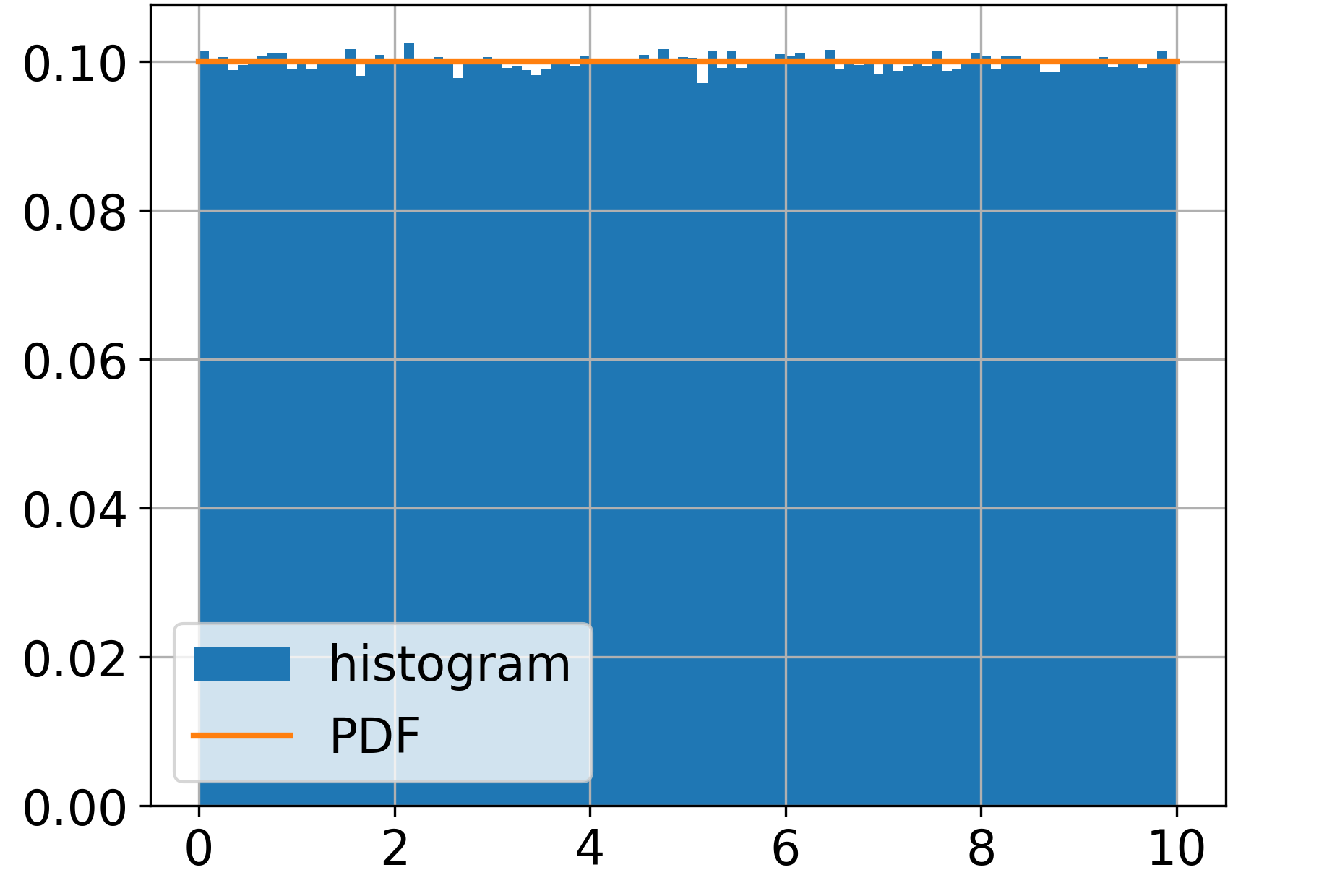}
    \caption{Stress jump $\beta_\JU$}
    \end{subfigure}
    \begin{subfigure}{0.32\linewidth}
    \includegraphics[height=3.7cm]{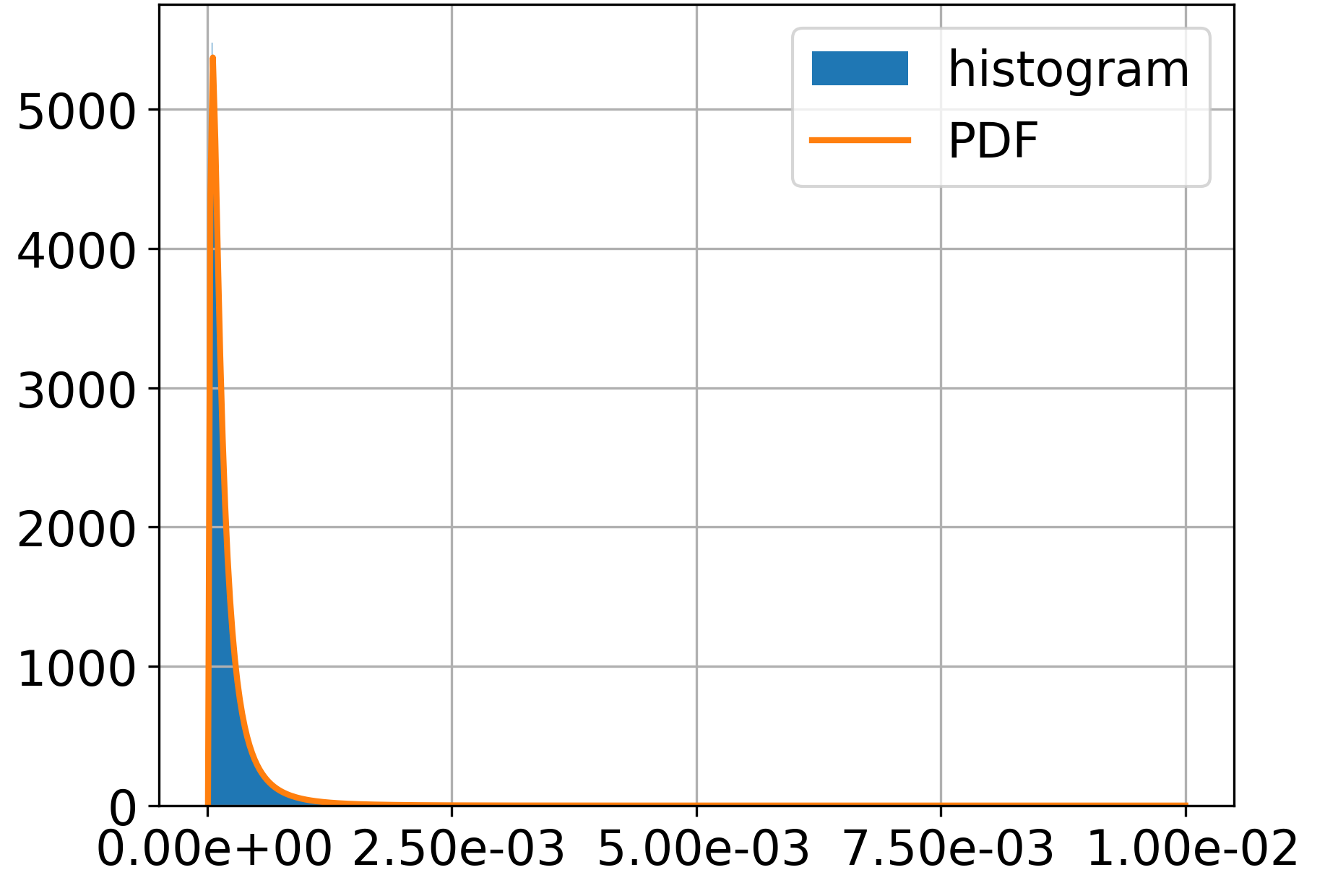}
    \caption{Permeability $k_\Gamma$}
    \end{subfigure} 
    \begin{subfigure}{0.32\linewidth}
    \includegraphics[height=3.7cm]{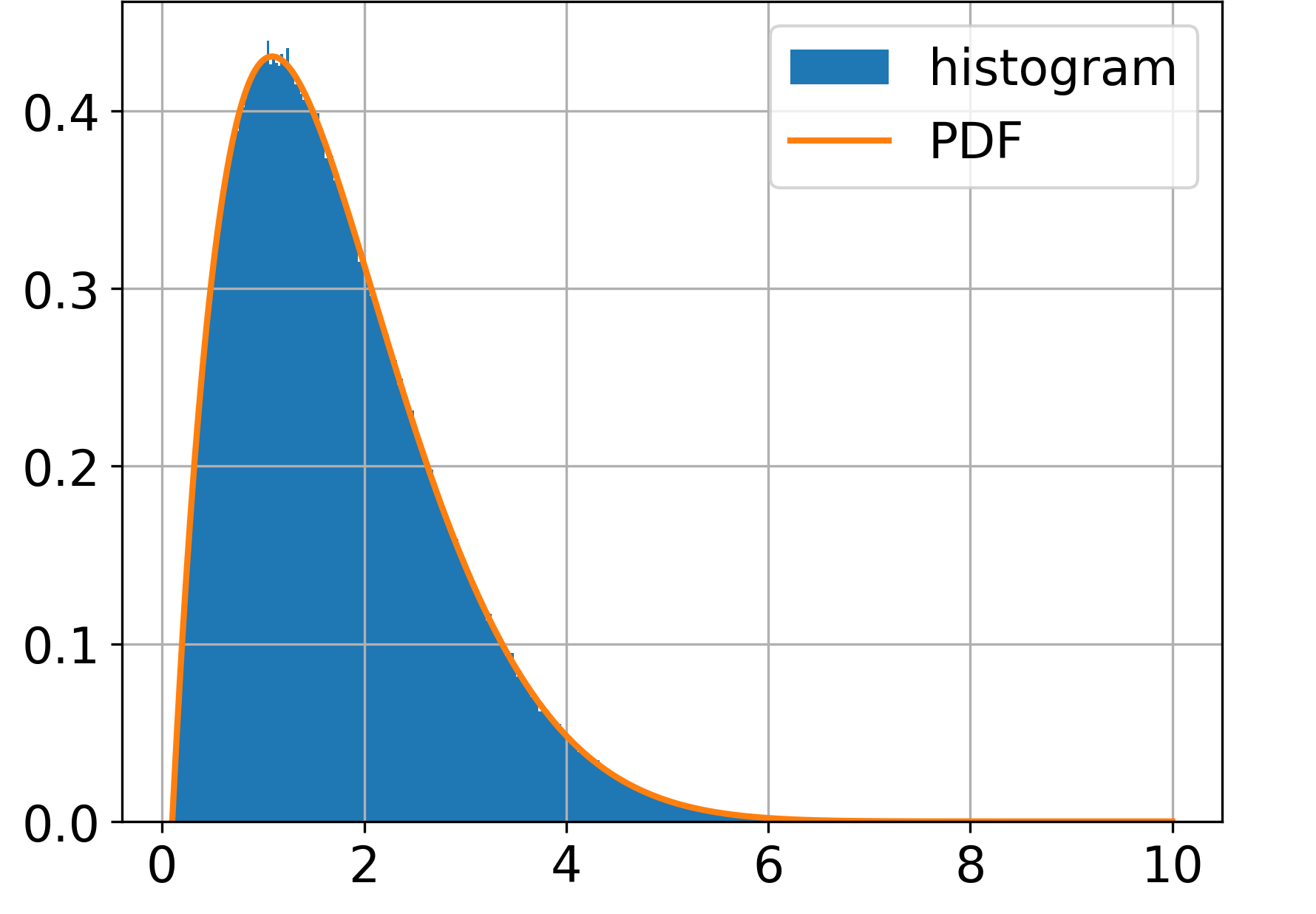}
    \caption{Effective viscosity $\mu_\EF$}
    \end{subfigure}\\[4mm]
    \begin{subfigure}{0.32\linewidth}
    \includegraphics[height=3.7cm]{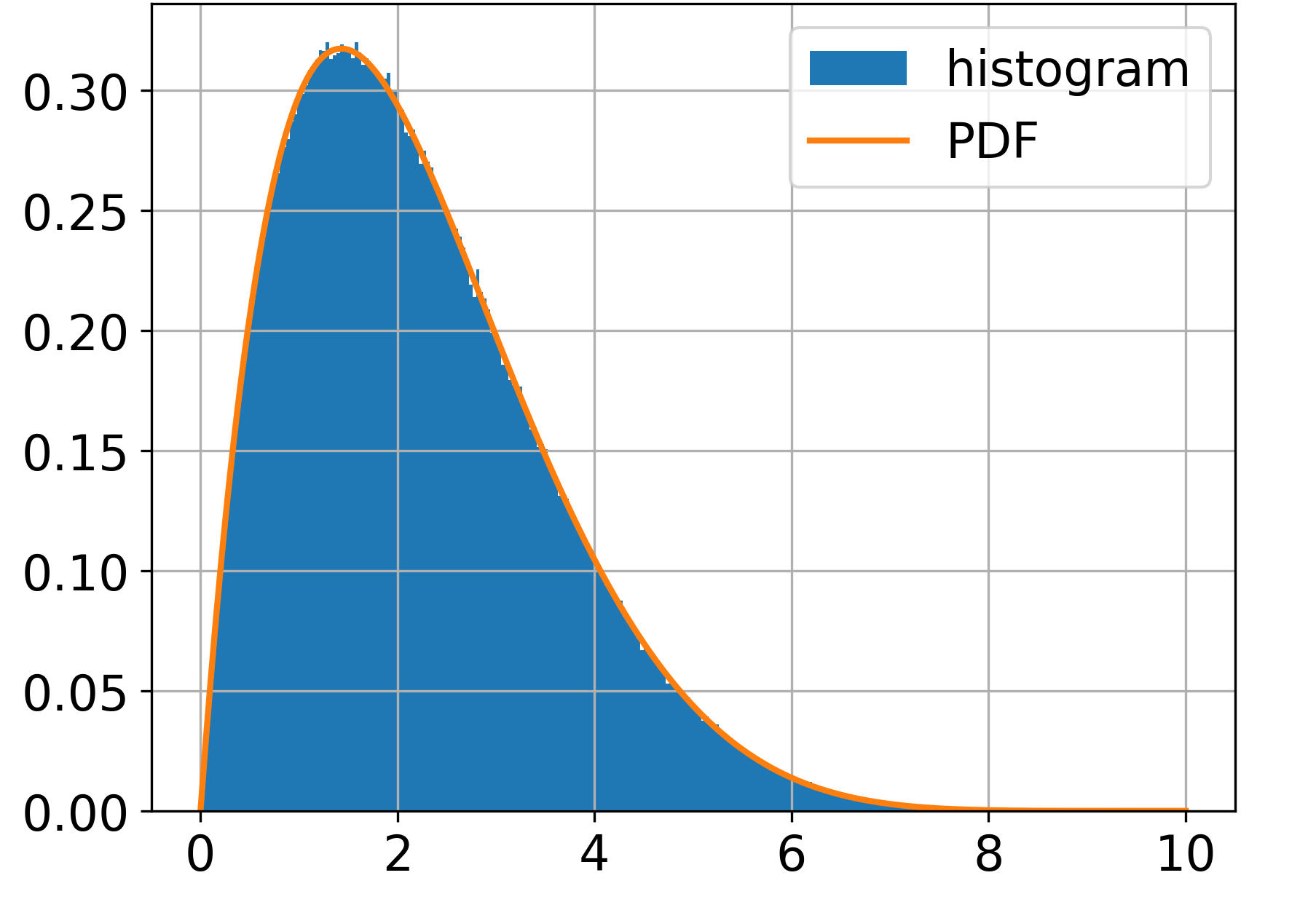}
    \caption{Beavers--Joseph $\alpha_\BJ$}
    \end{subfigure}
    \hfil
    \begin{subfigure}{0.32\linewidth}
    \includegraphics[height=3.7cm]{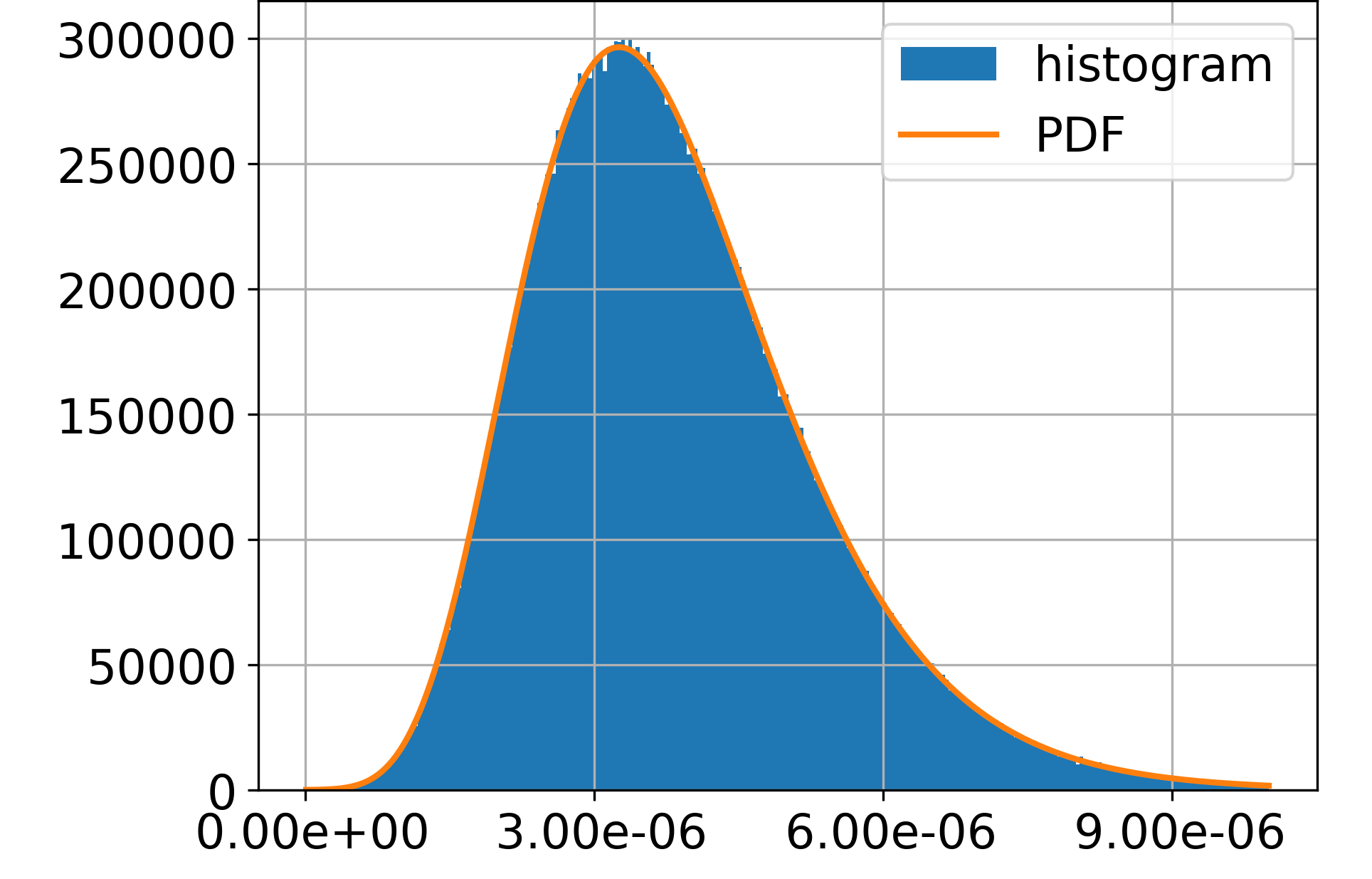}
    \caption{Permeability $k_\PM$}
    \end{subfigure}
    \caption{Probability density functions (PDF) and histograms of the uncertain parameters.}
    \label{fig:PDFparameters}
\end{figure}

The stress jump parameter $\beta_\JU$ is theoretically estimated to be of the order $\mathcal{O}(\sqrt{K_\Gamma}/d, d/\sqrt{K_\Gamma})$~\cite{OCHOATAPIA19952635,OCHOATAPIA19952647}. This represents a very large range and is not suitable for sensitivity analysis. Through the validation against the pore-scale simulation results~\cite{Ruan_Rybak_TIPM}, the stress jump parameter $\beta_\JU$ is treated as fitting parameter and its range can be reduced to $[0,10]$. We apply a uniform distribution for $\beta_\JU$, since no prior knowledge about its distribution is available.

The values of porous-medium permeability are taken into account $k_\PM \in (10^{-8}, 10^{-5})$ to reflect the commonly used porous materials. The log-normal distribution is applied to generate the initial data and reduce it to the interval of interest~\cite{RICCIARDI2005248,kroker2023global}. The permeability on the complex interface is expected to be larger than or equal to that in the porous-medium region, while retaining the same structural characteristics as in the porous regions.
 In this work, we consider the permeability on the complex interface in the range $k_\Gamma \in (10^{-5}, 10^{-2})$ and apply the log-normal distribution as for the porous-medium permeability in this interval. 

The value of the permeability viscosity $\mu_\EF$ in the Brinkman equations remains a debated issue. Depending on the type of porous medium, numerical simulations have shown that the effective viscosity may be either smaller or larger than the dynamic viscosity $\mu_\dyn$~\cite[Chap.~3.4.3]{vafai2005handbook}. Based on the volume-averaging approach of the Navier–Stokes equations~\cite{OCHOATAPIA19952635}, the effective viscosity is expressed as $\mu_\EF = \mu_\dyn/\epsilon$, where $\epsilon$ denotes the porosity. For applications of the Brinkman model, the porosity is typically in the range $\epsilon \in [0.8,0.95]$. In highly porous cases, it is often accepted that the effective viscosity should coincide with the dynamic viscosity, i.e., $\mu_\EF = \mu_\dyn$. In this work, we choose $\mu_\dyn=1$ and select a wider admissible range for the effective viscosity, $\mu_\EF \in [0.1,10]$. To capture plausible physical values, we employ a suitably shaped beta distribution that concentrates most of the probability mass around the interval $\left[\mu_\dyn, \mu_\dyn/\epsilon\right]$, with $\epsilon \in [0.8,0.95]$.

The slip coefficient $\alpha_\BJ$ is chosen originally from the Beavers--Joseph condition and the most commonly used value is $\alpha_\BJ = 1$ in the literature.  We consider a larger range of slip coefficient $\alpha_\BJ \in(0,10)$. We set the mode of the distribution function at $\alpha_\BJ=1$ and use a scaled beta distribution to show the suitable shape of probability density function. The intervals and distributions for these uncertain parameters are summarized in Tab.~\ref{table1:UQtable}. The probability density functions and the histograms of the uncertain parameters are plotted in Fig.~\ref{fig:PDFparameters}. In sections~\ref{sec:testcase1FP} and \ref{sec:testcase2SF}, we will introduce two test case used for sensitivity analysis: filtration problem and splitting flow. For both test cases, we study the impact of uncertain parameters from hybrid-dimensional model on the flow velocity $\vec{v}= (u,v)^\top$.

\subsection{Test case 1: Filtration problem} \label{sec:testcase1FP}
As the first test case, we introduce a filtration problem. In this case, a parabolic inflow is considered on the bottom of the porous-medium region, the flow passes through the complex interface and exits through the free-flow domain. Note that we have different sizes of the outflow on the left and right external free-flow boundary that causes the arbitrary flow near the fluid--porous interface. 

\begin{figure}[b]
    \centering
    \includegraphics{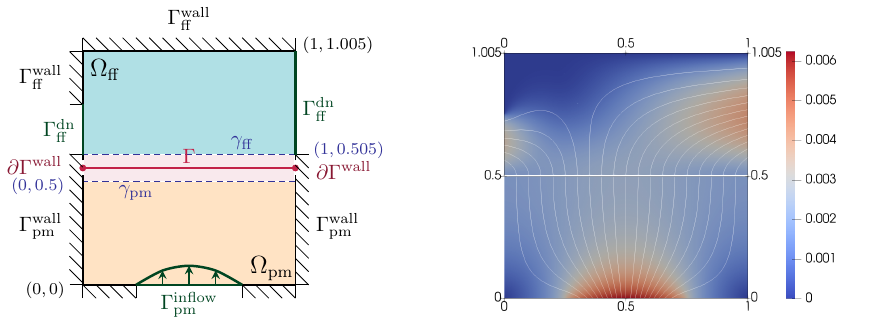}
    \caption{Schematic setting (left) and velocity (right) for filtration problem.}
    \label{fig:testcaseFP}
\end{figure}

We consider the flow domain $\Omega_\mathrm{F} = [0, 1]\times [0, 1.005]$ with the free-flow domain $\Omega_\FF= [0,1]\times [0.505, 1.005]$, the complex interface $\Gamma=[0,1] \times \{0.5025\}$ and the porous-medium region $\Omega_\PM= [0,1]\times [0,0.5]$ (Fig.~\ref{fig:testcaseFP}, left). The position of the transition region is $[0,1]\times[0.5,0.505]$ with thickness $d= 0.005$. On the bottom boundary, the following ``inflow" boundary condition
\begin{equation}
        \vec{v}_\PM \cdot \vec{n}_\PM= 0.1(x_1 - 0.25)(x_1-0.75) \quad \textnormal{on } \Gamma_{\PM}^\mathrm{inflow} = [0.25, 0.75] \times \{0\}
\end{equation}
is imposed. On the remaining bottom boundary and the side boundaries of the porous-medium domain $\Gamma_\PM^\WA= (\{0\}\times[0,0.5])\cup([0,0.25]\times\{0\}) \cup ([0.75, 1] \times \{0\})\cup (\{1\}\times [0,0.5])$, we take into account the no-slip boundary conditions
\begin{equation}
    \vec{v}_\PM \cdot \vec{n}_\PM = 0 \quad \textnormal{on } \Gamma_{\PM}^\WA. \label{eqn:porousmediumwallBC}
\end{equation}
On the boundary of the complex interface $\Gamma$, we apply the following no-slip boundary condition
\begin{equation}
     \vec{V}\cdot\vec{\tau} = 0, \quad \mu_\EF\frac{\partial V_\vec{n}}{\partial \vec{\tau}} =0  \quad \textnormal{on } \partial \Gamma^\WA=\{(0,0.5025),(1,0.5025)\}.
\end{equation}
On the right external boundary of the free-flow region, the ``do-nothing" boundary condition is considered
\begin{equation}
    \ten{T}(\vec{v}_\FF, p_\FF) \cdot \vec{n}_\FF = \vec{0}\quad \textnormal{on }\Gamma_{\FF}^\DN =\left( \{0\}\times [0.505, 0.755]\right) \cup\left(\{1\} \times [0.505, 1.005]\right),
\end{equation}
which allows the flow to exist the domain.  On the top and the remaining left boundaries of the free-flow region $\Gamma_{\FF}^\WA = \left(\{0\}\times [0.755, 1.005]\right)\cup( [0,1] \times \{1.005\})$, we have
\begin{equation}
    \vec{v}_\FF = \vec{0} \quad \textnormal{on } \Gamma_{\FF}^\WA. \label{eqn:freeflowwallBC}
\end{equation}
The source terms are  $\vec{f}_\FF=\vec{0}$, $F_{\vec{\tau}}=F_\vec{n}=0$ and $f_\PM=0$. The fluid viscosity is set as $\mu_\dyn=1$. The typical values and the intervals of the stress jump parameter $\beta_\JU$ and permeability of the complex interface $k_\Gamma$ in Eqs.~\eqref{eqn:AveragedBrinkmanMomentumNormal}--\eqref{eqn:TransmissionGammaFFtangential}, effective viscosity $\mu_\EF$ in Eqs.~\eqref{eqn:AveragedBrinkmanMomentumNormal}--\eqref{eqn:TransmissionGammaPM}, slip coefficient $\alpha_\BJ$ in Eqs.~\eqref{eqn:AveragedBrinkmanMomentumTangential}, \eqref{eqn:TransmissionGammaFFtangential} and porous-medium permeability $k_\PM$ in Eqs.~\eqref{eqn:Darcy2}, \eqref{eqn:AveragedBrinkmanMomentumTangential}, \eqref{eqn:TransmissionGammaFFtangential} are provided in Tab.~\ref{table1:UQtable}.

The hybrid-dimensional Stokes--Brinkman--Darcy problems are discretised by the second order finite volume method on staggered rectangular grid of size $h_x\times h_y$ in free-flow and porous-medium regions and second-order finite difference method along the complex interface $\Gamma$ with the grid size $h_x$~\cite{Ruan_Rybak_TIPM}. The grid size is fixed $h=h_x = h_y =1/100$.  The velocity magnitude and the streamlines for the filtration problem are visualized in Fig.~\ref{fig:testcaseFP} (right).

\subsection{Test case 2: Splitting flow problem}\label{sec:testcase2SF}
Generally, the velocity in the free-flow region is much larger than in the porous-medium region. Therefore, as an additional test case, we consider a splitting flow problem, where flow enters the domain from the top boundary and exits through the side boundaries of the free-flow domain. In this setting, the flow in the free-flow region is significantly faster than that in the porous medium. Similar to the case presented in~\ref{sec:testcase1FP}, different sizes of the outflow boundaries on the external free-flow domain lead to a non-parallel flow to the fluid–porous interface.
\begin{figure}[b]
    \centering
    \includegraphics{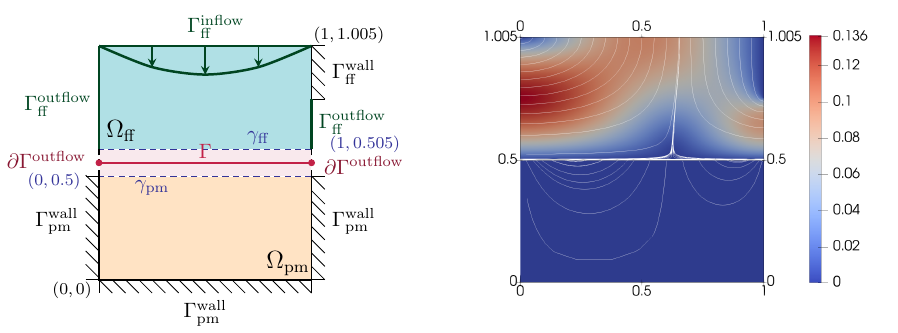}
    \caption{Schematic setting (left) and velocity (right) for splitting flow problem.}
    \label{fig:testcase2SFgeo}
\end{figure}

In this case, we keep the same geometric setting as the filtration problem in section~\ref{sec:testcase1FP}, but different settings of boundary conditions are taken into account (Fig.~\ref{fig:testcase2SFgeo}, left). On the top of the free-flow region, we prescribe the following ``inflow'' boundary condition
\begin{equation}
\vec{v}_\FF = \left(-0.1 \sin{(\pi x_1)},0\right)^\top  \quad \textnormal{on } \Gamma^\mathrm{inflow}_\FF = [0,1] \times \{1.005\}. 
\end{equation}
On the side boundaries of the free-flow domain $\Gamma_\FF^\OUT = \left(\{0\}\times [0.505, 1.005]\right) \cup \left(\{1\} \times [0.505, 0.755]\right)$, the following ``outflow'' boundary condition 
\begin{equation}
    \vec{\tau} \cdot \ten{T} (\vec{v}_\FF , p_\FF)  \cdot \vec{\tau} = 0, \quad \vec{v}_\FF \cdot \vec{n} = 0 \quad \textnormal{on } \Gamma_\FF^\OUT
\end{equation}
is applied. Similarly, on the boundary nodes of the complex interface $\Gamma$, we set
\begin{align*}
            \mu_\EF \frac{\partial V_\vec{\tau}}{\partial \tau} - P=0, \quad \vec{V}_\vec{n} =0 \quad \textnormal{on }\partial \Gamma^\OUT=\{(0,0.5025),(1,0.5025)\}.
\end{align*}
We take into account ``wall'' boundary condition~\eqref{eqn:freeflowwallBC} on the remaining right boundary of the free-flow region $\Gamma_\FF^\WA= \{1\}\times[0.755, 1.005]$ and the no-slip boundary condition~\eqref{eqn:porousmediumwallBC} on all external boundaries of the porous medium~$\Gamma_\PM^\WA = \left(\{0\}\times [0,0.5]\right)\cup\left([0,1]\times \{0\}\right)\cup\left(\{1\}\times[0,0.5]\right)$.

The source terms are set to $\vec{f}_\FF=\vec{0}$, $F_{\vec{\tau}}=F_\vec{n}=0$, and $f_\PM=0$. The ranges of the stress jump $\beta_\JU$, the interface permeability $k_\Gamma$, the effective viscosity $\mu_\EF$, the slip coefficient $\alpha$, and the porous-medium permeability $k_\PM$ can be found in Tab.~\ref{table1:UQtable}. The coupled problem is discretised using the same numerical method as in section~\ref{sec:testcase1FP}. The velocity for the splitting flow is visualised in Fig.~\ref{fig:testcase2SFgeo} (right).

\section{\aMR/PCE-assisted global sensitivity analysis.}
\label{sec:GSAonPols}

In this section, we briefly recall variance-based GSA, introduce PCE/\aMR{} surrogate fundamentals, and then show how to compute GSA directly from PCE/\aMR{} coefficients.
\subsection{GSA using Sobol' indices}
\label{sec:gsa}
The origins of variance-based sensitivity analysis go back
to the concepts of ANOVA (analysis of variance) originally proposed in the work of Fisher~\cite{Fisher_1921,fisher1923manurial} published in the 1920s in the context of
agriculture.  We briefly recall the ANOVA projection framework to create the link with Sobol’ indices. The core idea is as follows.
Let $\rf(\vrv):=\rf(\rvi1,\ldots,\rvi\sdim)$
be a square-integrable random variable or a random field. 
Let $\Prj{I}$, for $I\subseteq\cN:=\left\{1,\ldots,\sdim\right\}$
be an orthogonal projection operator satisfying
\begin{align*}
    \Prj\emptyset\rf&:=\E[\rf],\\
    \Prj{\{j\}}\rf &:=\int \rf\dx{\Pm(\vrv)}_{-\{j\}} - \Prj\emptyset \rf,
    \quad j\in\cN,\\
    \Prj{I}\rf &:=\int \rf \dx{\Pm(\vrv)}_{-\{I\}}
     -\sum_{J\subsetneq I}\Prj{J}\rf,\quad I\subseteq\cN.
\end{align*}
Here, we denote by $\int \cdot \dx{\Pm(\vrv)}_{-\{j\}}$ the integration across all variables $\rvi{i}$ for each $i \in \cN$, where $i \neq j$. This leads to the following decomposition of $\rf$:
\begin{align*}
    \rf(\rvi1,\ldots,\rvi\sdim) & = \Prj\emptyset\rf 
    + \sum_{i=1}^\sdim \Prj{\{i\}}\rf(\rvi{i})
    +\sum_{1\leq i\leq j\leq\sdim} \Prj{\{i,j\}}\rf(\rvi{i},\rvi{j})
    +\ldots=\sum_{I\subseteq\cN}\Prj{I}\rf(\rvi{I}).
\end{align*}
This also implies the following decomposition of the variance:
\begin{equation}\label{eq:var-dec}
    \var = \sum_{I\subseteq\cN} \var_I,
\quad\text{where}\;\; \var_{\emptyset}:=0,
\quad \var_I:=\int \left(\Prj{I}\rf(\vrv) \right)^2\dx{\Pm(\vrv)}.
\end{equation}
The normalization by the variance leads to the definition of Sobol' sensitivity indices~\cite{MR1052836, MR1823119}:
\begin{equation*}
    \sobi{I}:= \frac{\var_{I}}{\var},\quad\text{for } I\subseteq\cN.
\end{equation*}
Using Sobol' sensitivity indices $\sobi{I}$, we can define the total sensitivity
indices $\sobit{i}$ describing the total contribution of each input parameter $i\in\cN$
to the variance. The total sensitivity indices are given by
\begin{equation*}
    \sobit{i}:=\sum_{I\in J_i} \sobi{I},\quad
    \text{for}\; 
    J_i=\left\{(i_1,\ldots,i_s)\subseteq\cN \mid \exists\,k, 1\leq k \leq s,\, i_k=i \right\}.
\end{equation*}

\subsection{Introduction to PCE and \aMR}
The concept of PCE, also the use of orthogonal polynomials
to build a surrogate model, 
was initially introduced by Wiener~\cite{wiener1938} and has since been widely used and extended to date.
In the context of uncertainty quantification (UQ), PCE was introduced by Ghanem and Spanos in \cite{Ghanem91StochFE}, and since then applied in several areas, e.g.~\cite{MR2121216,MR2501693,Xiu2003137}.
However, the core idea of the PCE of the second-order
random field $\rf(\rv)\equiv\rf(\vec{x}, \rv(\cdot))\in\Lp2(\sspace)$
has the following form
\begin{align*}
    \rf(\vec{x},\rv)\approx\sum_{p=0}^\No \pcf{c}{p}(\vec{x})\pol{p}(\rv),
    \quad
    \pcf{c}{p}(\vec{x})&=\sprod{\rf(\vec{x},\cdot)}{\pol{p}}.
\end{align*}
Here,  $\No$ is the maximal polynomial degree, and at least in the classical formulation, $\pcf{c}{p}(\vec{x})$ or $\pcf{c}{p}$ denotes the deterministic function or the deterministic coefficient in the case of $\rf$ being a random variable,
$\sprod{\cdot}{\cdot}$ is equivalent to expectation $\E[\cdot]$
with respect to the probability measure $\Pm$, such that
orthonormal polynomials $\pol{p},\, p\leq\No$ satisfy
\begin{align*}
    \sprod{\pol{p}}{\pol{q}}&=\int_\sspace \pol{p}(\rv)\pol{q}(\rv)\dx{\Pm(\rv)}
    =\delta_{p,q},\quad p,q=0,\ldots,\No.
\end{align*}
The selection or construction of the polynomial basis $\{\pol{p}\}_{p\leq\No}$ depends directly on the probability distribution of the input parameter $\rv$.

Apart from typical probability distributions like the uniform distribution (associated with Legendre polynomials) or the normal distribution (linked with Hermite polynomials), there exist various construction methods of the corresponding orthogonal polynomials. 
There is a three-term recurrence relation based on the work of Favard~\cite{favard1935polynomes}.
Another concept is the generalized PCE \cite{Xiu2003137} or the
Wiener--Askey polynomial chaos \cite{MR1951058}, which relies on the Askey polynomials \cite{MR783216}.
The third method, which we would like to mention here, is  
arbitrary polynomial chaos (aPC), which builds orthogonal polynomials
directly from the Hankel matrix \cite{OLADYSHKIN2012179}, 
and therefore relies only on the raw moments of the underlying random variable $\rv$.
However, note that the three-term recurrence relation
can also be constructed using the Hankel matrix only. 
Here we refer to \cite{MR2061539}.
Existence of orthogonal polynomial bases depends on the properties of
the raw moments of $\rv$ and in particular on the positive definiteness
of the Hankel matrix. Here, we refer the readers to  \cite{MR2855645,MR2061539,sullivan_book} for details.

However, independent of the basis construction strategy,
one of the drawbacks of the classical PCE approach is their susceptibility
to oscillations due to irregularities of the model response (Gibbs' phenomena)
or suboptimal location of the training points (the Runge phenomena).
There exist different approaches to reduce or avoid oscillations of the PCE-surrogate
model response. First of all, there are the different regularization and sparsity strategies such as LARS-LASSO regression \cite{MR2764550,Luethen2021,KUSCH2020109073},
Bayesian sparse PCE \cite{buerkner2022} or PCE-Kriging \cite{Schobi_2015}.
Another approach is the use of multi-resolution \cite{MR2063905}
and multi-element \cite{MR2240796} frameworks,
or a combination of the multi-resolution / multi-element framework and regularization~\cite{kohlhaas2023gaussian,kroeker2025bayesian3}.
The typical trade-off between both strategies is, at least in our experience, the higher computational effort for training the PCE-surrogate model due to
additional optimization steps in
case of regularization and higher amount of training samples necessary for building the multi-resolution / multi-element
surrogates.
Since in our application cases, on the one hand, computational costs of physical models are moderated, and on the other hand, the requirement to deal with a $100\times100$
random fields yielding $10000$ PCE expansions, we omit regularization
and proceed with multi-resolution / multi-element PCE Ansatz.

The initial a multi-resolution concept was proposed by Le~Ma{\^i}tre et al.~\cite{MR2063905} 
and later extended to the  multi-element by Wan et al.~\cite{MR2240796}.
Since then it has been used in intrusive approaches \cite{MR3023712,Pettersson2016367,Buerger2012,Koeppel2017}
and later also in non-intrusive approaches \cite{Marelli_2021,WAGNER2022102179,DREAU2023117920}.
Here, we use arbitrary multi-resolution polynomial chaos (\aMR)
proposed in \cite{kroker2022arbitrary}. 
The related Python implementation is available in~\cite{ikroeker_aMR-PC_2023}.

Following the idea introduced by Le~Ma{\^i}tre et al.~\cite{MR2063905}
we split the sample space / dataset $\sspace$ into
stochastic subdomains (\SD) using dyadic decomposition
\begin{equation*}
    \subs{l}:=\cdf{\rv}^{-1}\left([2^{-\Nr}l,\,2^{-\Nr}(l+1)] \right),
    \quad
    l=0,\ldots, 2^{\Nr}-1.
\end{equation*}
Here, $\cdf\rv(\cdot)$ denotes the cumulative distribution function, and
$\Nr$ is the number of dyadic refinements.
However, by assuming $\sspace\subseteq\R$ 
(this assumption is valid when working with datasets), we can use the quantile function defined as follows
\[
\quant{\rv}(r) :=\inf\left\{s\in\R\mid r\leq\cdf\rv(s) \right\},
\]
to obtain \SD's 
\begin{equation*}
    \subs{l}=\quant{\rv}\left([2^{-\Nr}l,\,2^{-\Nr}(l+1)] \right),
    \quad
    l=0,\ldots, 2^{\Nr}-1.
\end{equation*}
In the next step, we use one of the moment-based methods (aPC~\cite{OLADYSHKIN2012179} or the three-term method described by Gautschi \cite{MR2061539}) to construct orthogonal polynomial bases
in each \SD. 
By normalizing the polynomials, we obtain the piecewise polynomial basis
$
    \left\{\mrpol{l,p}\mid l=0,\ldots,2^{\Nr}-1,\, p=0,\ldots,\No \right\}
$,
satisfying the orthonormal relation on the whole stochastic domain $\sspace$
as follows
\begin{equation*}
    \sprod{\mrpol{l,p}}{\mrpol{k,q}}=\delta_{l,k}\delta_{p,q},
    \quad
    k,l=0,\ldots,2^{\Nr}-1, \; \; p,q=0,\ldots,\No.
\end{equation*}

Now we describe the multivariate formulation
which provides \aMR{}-surrogate models with multiple uncertain parameters.
Let $\sdim$ be the number of uncertain parameters of the physical model.
Since in the present work we use the same refinement level $\Nr\in\N_0$
for each uncertain parameter, we assume it for notation.
Let us define the set of multi--indices $\idxSet$ addressing all \SD's as follows
\[
\idxSet:=\left\{\midx{l}=\left(\midxi{1},\ldots,\midxi{\sdim}\right)^T\in\N_0^\sdim
\mid \midxi{i}=0,\ldots,2^{\Nr}-1,\;  i=1,\ldots,\sdim \right\}.
\]
In other words, we obtain
$
\sspace^\sdim = \bigcup_{\midx{l}\in\idxSet} \msubs{\midx{l}}
$,
where
$
\msubs{\midx{l}} := \subs{\midxi{1}}\times\cdots\times\subs{\midxi{\sdim}}
$.
To complete the description, we introduce polynomial degree-related multi-indices.
Here, we follow the standard notation and use Greek letters
$\alpha=\left(\alpha_1, \ldots,\alpha_\sdim\right)^T\in\N_0^\sdim$.
This means that for the multi-dimensional variable $\vrv=(\rvi{1},\ldots,\rvi\sdim)$
the monomial $\vrv^\alpha$ has the form
$\vrv^\alpha=\rvi1^{\alpha_1}\cdots\rvi{\sdim}^{\alpha_\sdim}$.
The term $\abs{\alpha}=\alpha_1+\ldots+\alpha_\sdim$ represents the total degree of the monomial $\vrv^\alpha$.
Using the notation above, we obtain the polynomial basis
$\left\{\mrPol{\midx{l},\alpha}\mid \midx{l}\in\idxSet,\; \alpha\in\N_0^\sdim \right\}$
that is built out of the multivariate piecewise polynomials constructed as follows
\begin{equation*}
    \mrPol{\midx{l},\alpha}(\vrv):=
    \begin{cases}
    \prod_{i=1}^\sdim \mrpol{\midxi{i},\alpha_i}(\rvi{i}), & \vrv\in\msubs{\midx{l}},\\
    0, &\text{else}.
    \end{cases}
\end{equation*}
Furthermore, multivariate piecewise polynomial functions $\mrPol{\midx{l}, \alpha}$
are orthonormal in the sense of
\begin{equation*}
    \sprod{\mrPol{\midx{l},\alpha}}{\mrPol{\midx{m},\beta}}
    =\delta_{\midx{l},\midx{m}}\delta_{\alpha,\beta},
    \quad\text{for}\; \midx{l},\midx{m}\in\idxSet,\, \alpha,\beta\in\N_0^\sdim.
\end{equation*}
Now, we can define the surrogate model $\surr{\Nr,\No}\equiv\surr{\Nr,\No}(\vec{x};\vrv)$ of the
second-order random field $\rf\equiv\rf(\vec{x};\vrv)$ using the following \aMR{} expansion
\begin{align}\label{eq:amr-surr}
    \surr{\Nr,\No}(\vec{x};\vrv)&:= \sum_{\midx{l}\in\idxSet}\sum_{\alpha\in\pidxSet}
    \mrcf{\rf}{\midx{l},\alpha}(\vec{x})\mrPol{\midx{l},\alpha}(\vrv),\\
    \label{eq:amr-cf}
    \mrcf{\rf}{\midx{l},\alpha}(\vec{x})&:=\sprod{\rf(\vec{x};\cdot)}{\mrPol{\midx{l},\alpha}}.
\end{align}
Here, $\pidxSet=\left\{\alpha\in\N_0^\sdim\mid \abs{\alpha}\leq\No \right\}$
denotes the set of polynomial degrees $\alpha$ truncated by the highest polynomial
degree \No.
It should also be mentioned that the \aMR{} expansion defined above
is identical to the multivariate PCE in case $\Nr=0$.
Despite the fact that the polynomial coefficients $\mrcf{\rf}{\midx{l},\alpha}$
can be computed according to ~\eqref{eq:amr-cf} using, for example,
a numerical quadrature, the established approach for the non-intrusive
method is to use an appropriate least-squares approach, e.g.~\cite{SUDRET2008964,kroker2022arbitrary}.
In the present work, we apply the Penrose pseudoinverse \cite{penrose1956best}
to solve the linear systems and to obtain the polynomial coefficients.
One of the key properties of the PCE / \aMR{} based surrogates,
which is also relevant for the present work, is the calculation
of the mean $\smean{\Nr,\No}$ and variance $\svar{\Nr,\No}$ of the surrogate $\surr{\Nr,\No}(\cdot)$
using only the polynomial coefficients
\begin{align}\label{eq:amr-mean}
    \smean{\Nr,\No}&:=\E[\surr{\Nr,\No}]=\sum_{\midx{l}\in\idxSet}
    \mrcf[\Nr]{\rf}{\midx{l}} \sprod{\mrPol{\midx{l},0}}{\mrPol[0]{0,0}}
    =\sum_{\midx{l}\in\idxSet} \mrcf[0]{\rf}{\midx{l}} \cdot 2^{-\sdim\Nr/2},\\
    \label{eq:amr-var}
    \svar{\Nr,\No}&:=\var\left[\surr{\Nr,\No} \right]
    =\sum_{\midx{l}\in\idxSet}\sum_{\alpha\in\pidxSet}
    \left(\mrcf{\rf}{\midx{l},\alpha} \right)^2 - \smean{\Nr,\No}^2.
\end{align}

The size of the linear system that must be solved for the training
of the surrogate model $\surr{\Nr,\No}(\cdot)$ is given by the amount of training samples
and the number $\Nc$ of polynomial coefficients $\mrcf{\rf}{\midx{l},\alpha}$. 
There is a clear dependence of the necessary number of training samples on $\Nc$.
However, the number of polynomial coefficients of the full-tensor approach
is given by
\begin{equation*}
    \Nc = 2^{\sdim\Nr}\frac{(\No+\sdim)!}{\No!\,\sdim!}.
\end{equation*}
Here, the symbol $(!)$ denotes the factorial.
Since the number of polynomial coefficients $\Nc$, and therefore also the
problem size, rapidly increases in almost each non-trivial setting, the
reduction of the system size, called sparsity, was addressed by many
researchers.
One of the few methods that can be performed on the given training data
without additional computational effort is \emph{hyperbolic truncation} \cite{MR2764550,Luethen2021}.
The idea here is to reduce the number of \textit{interaction order terms} 
by employing the $q$-norm 
$\norm{\alpha}{q}=\left(\sum_{i=1}^\sdim \abs{\alpha_i}^q \right)^{1/q}$
for truncation limited by maximal polynomial degree.
More precisely, we replace $\pidxSet$ in \aMR{} expansion \eqref{eq:amr-surr}
by
\[
\pidxSet[,q]:=\left\{\alpha\in\N_0^\sdim\mid \norm{\alpha}{q}\leq\No\right\},
\quad\text{for}\;
q\in(0,\,1].
\]
We mention here that in the case $q=1$ the hyperbolic truncation is
equivalent to the classical one, such that
$\pidxSet[,1]\equiv\pidxSet$. For $q<1$, higher-order and high-interaction terms are preferentially pruned.

\subsection{GSA using PCE and \aMR{} expansion coefficients}

In this section, we summarize how to compute Sobol' sensitivity indices from PCE/\aMR{} expansion coefficients used in the present work.
For a more comprehensive description and validation of the methods,
we refer, e.g., to \cite{SUDRET2008964,kroker2023global}. This subsection connects Section~\ref{sec:gsa} (ANOVA) with the surrogate definitions above.
\subsubsection{Sobol' sensitivity indices using PCE}
Let us start with the calculation of the Sobol' indices using the (global) PCE.
Since the basis polynomials of the PCE have global support, we can
omit the first sum in Eq.~\eqref{eq:amr-surr} and obtain
$$
\surr{\No}\equiv\surr{0,\No}
=\sum_{\alpha\in\pidxSet}\pcf{\rf}{\alpha} \Pol{\alpha}(\vrv), 
$$
or
$$
\surr{\No,q}\equiv\surr{0,\No,q}
=\sum_{\alpha\in\pidxSet[,q]}\pcf{\rf}{\alpha} \Pol{\alpha}(\vrv),
$$
if the hyperbolic truncation was used. 
Since $\surr{\No}\equiv\surr{\No,1}$ and $\pidxSet\equiv\pidxSet[,1]$, we
will proceed with the notation related to the hyperbolic truncation.
Now we can exploit the structure of the polynomial degrees of the multivariate
PCE to obtain the orthogonal projection operator $\Prj{I}$ introduced
in Section~\ref{sec:gsa}.
Taking advantage of the fact that the degree of polynomial $\alpha_i=0,\, i=1,\ldots,\sdim$
corresponds to the mean in the related direction, 
define the set $\pidxSet[,q]_I$ of tuples
$\alpha=\left(\alpha_1,\ldots,\alpha_\sdim\right)\in\N_0^\sdim$ with non-zero entries only in
the indices $I=(i_i,\ldots,i_s)\subseteq\cN$:
\[
\pidxSet[,q]_I :=\left\{\alpha\in\N_0^\sdim \mid 
\norm{\alpha}{q}\leq\No,
\begin{array}{ccc}
    \alpha_i > 0, & \forall\; i\in\cN, & i\in I \\
    \alpha_i = 0, & \forall\; i\in\cN & i\notin I
\end{array}
 \right\},
 \quad I\subseteq \cN,\, q\in(0,1].
\]
For example, it is easy to see that for $\alpha \in\pidxSet[,q]_{\{i\}}$ the multivariate
polynomial $\Pol{\alpha}(\vrv)$ depends only on the uncertain parameter $\rvi{i}$.
This leads to the Sobol' decomposition of the PCE surrogate model $\surr{\No,q}$:
\[
\surr{\No,q}(\vrv)=\sum_{I\subseteq\cN} \Prj{I}\surr{\No,q}(\rvi{I})
=\pcf{\rf}{0} + \sum_{I\subseteq\cN}\sum_{\alpha\in\pidxSet[,q]_I} 
\pcf{\rf}{\alpha}\Pol{\alpha}(\vrv).
\]
Using~\eqref{eq:amr-var} we obtain the Sobol' sensitivity indices 
$\sobi{I}[\surr{\No,q}]$
 of the PCE surrogate model $\surr{\No,q}$, 
which approximate the Sobol' sensitivity indices $\sobi{I}$ of the random field \rf:
\[
\sobi{I}\approx\sobi{I}[\surr{\No,q}]=\sum_{\alpha\in\pidxSet[,q]_I}
\pcf{\rf}{\alpha}^2\underbrace{\E[\Pol{\alpha}^2]}_{=1}/\svar{\No,q}
=\sum_{\alpha\in\pidxSet[,q]_I}
\pcf{\rf}{\alpha}^2/\svar{\No,q},
\quad\text{for}\; I\in\cN.
\]
Therefore, the total sensitivity index $\sobit{i},\, i\in\cN$
can be approximated by
\[
\sobit{i}\approx\sobit{i}[\surr{\No,q}]
=\sum_{J\in\cN, i\in J} \sobi{J}[\surr{\No,q}]
=\sum_{J\in\cN, i\in J} \sum_{\alpha\in\pidxSet[,q]_J}
\pcf{\rf}{\alpha}^2/\svar{\No,q}.
\]
\subsubsection{Sobol' sensitivity indices using \aMR{}}
Now, we can extend the idea of the calculation of the Sobol' sensitivity
indices directly from the expansion coefficients discussed above
to the piecewise polynomial setting of \aMR{}.
Here, we skip the step-by-step calculation and refer to \cite{kroker2023global}
where the approach was published.
Nevertheless, although the underlying principle of this approach is similar to the PCE scenario, the primary distinction lies in the fact that zero-degree polynomials contribute to the variance when $\Nr>0$.

Again, we define the set of polynomial degree tuples $\alpha\in\N_0^\sdim$
corresponding to the index set $I\subseteq\cN$ by
\begin{equation*}
    \mrpSet[,q]_I:=\left\{ \alpha=(\alpha_1,\ldots,\alpha_\sdim)\in\N_0^\sdim\mid 
    \norm{\alpha}{q}\leq\No, \alpha_k=0,\, \forall\, k\in\cN,\, k\notin I\right\},
    \quad\text{for}\, I\subseteq\cN,\, q\in(0,1].
\end{equation*}
Using this notation, we can write the variance index 
(given in~\eqref{eq:var-dec})
of the \aMR{} surrogate model as follows
\begin{multline*}
    \svar{\Nr,\No,q}_{,I}:=
    \sum_{\midx{l}\in\idxSet}\sum_{\midx{m}\in\idxSet}\sum_{\alpha\in\mrpSet[,q]_I}
    \mrcf{\rf}{\midx{l},\alpha}\mrcf{\rf}{\midx{m},\alpha}
    2^{-\Nr\#(\cN\setminus I)} \delta_{\midx{l},\midx{m},I}\\
    +\sum_{J\subsetneq I} (-1)^{\#(I\setminus J)}
    \sum_{\midx{l}\in\idxSet}\sum_{\midx{m}\in\idxSet}\sum_{\alpha\in\mrpSet[,q]_J}
    2^{-\Nr\#(\cN\setminus J)} \delta_{\midx{l},\midx{m},J},\hspace{+5ex}
\end{multline*}
where $\#(A)$ denotes the count of elements in $A$.  
Further, we use the notation
\[
\delta_{\midx{l},\midx{m},I}:=
\begin{cases}
    1, & \midxi{i} = \midxi[m]{i},\, \forall\, i\in I,\\
    0, & \text{else}.
\end{cases}
\]
This allows us to approximate the Sobol' sensitivity indices $\sobi{I},\, I\subseteq\cN$
and the total sensitivity indices $\sobit{i},\, i\in\cN$ of the random field $\rf$ by
\[
\sobi{I}\approx \sobi{I}[\surr{\Nr,\No,q}]:=
\frac{\svar{\Nr,\No,q}_{,I}}{\svar{\Nr,\No,q}},
\quad
\sobit{i}\approx \sobit{i}[\surr{\Nr,\No,q}]:=
\sum_{J\subseteq\cN, i\in J} \frac{\svar{\Nr,\No,q}_{,J}}{\svar{\Nr,\No,q}}.
\]

\section{Numerical experiments}
\label{sec:numex}
In this section, we address the approximation and analysis of the physical models described in Section~\ref{sec:testcase1FP} and Section~\ref{sec:testcase2SF}
using PCE/\aMR{}-based surrogate models, introduced in Section~\ref{sec:GSAonPols}.

It is known that non-regularized (piecewise) polynomial-based surrogate
models provide an accurate surrogate model if the training samples
are selected appropriately.
Otherwise, at least for higher polynomial degrees, the surrogate model
response exhibits oscillatory behavior (Runge phenomenon).
In the worst case, this can lead to each surrogate requiring its own tailored training set.

In this study, we evaluate the utility of different polynomial-based surrogate models,
presented in Section~\ref{sec:GSAonPols},
whether piecewise or classical, that are developed using identical space-filling
sequences of training samples.  More precisely, we employ QMC Sobol' sequences~\cite{sobol1967} with training sample sizes between 512 and 16,384, corresponding to $2^9$ through $2^{14}$ samples. 
This design also allows us to assess the utility of the considered polynomial-based
surrogate models when training data are given.
In both cases, we assume that the uncertain input parameters of the physical model comply with the probability distributions listed in Tab.~\ref{table1:UQtable}.

To assess the reliability of the expansion coefficients,
we evaluated the $\Lp2$-error of the mean and standard deviation over the domain
based on the expansion coefficient of the two quantities of interest $(u,v)$
compared to the quantities based on the MC calculated from 50,000 model evaluations.
In order to evaluate the overall accuracy of the surrogate models, we calculate the root mean squared error (RMSE), which indicates the error on the scale of the quantities of interest. Additionally, we compute the relative mean square error (MSE), which is scale-independent due to normalization.
Based on this evaluation, we proceed to the GSA based on the expansion coefficients.
Here, we consider the spatial distribution of the total sensitivity indices, as well as space-averaged values of the Sobol' and total sensitivity indices.
To reduce the influence of numerical noise on the GSA, we exclude areas of vanishing
variance ($\var<10^{-10}$) while computing the space-averaged indices.

\subsection{Case~1: Postprocessing and evaluation}
We begin with the postprocessing and evaluation of the filtration problem introduced in Section~\ref{sec:testcase1FP}. In this configuration, the prescribed inflow at the bottom boundary drives the fluid in the direction of the complex interface, resulting in an extensive exchange between the porous-medium and free-flow regions (Fig.~\ref{fig:testcaseFP}). Consequently, the velocity field $\vec{v} = (u,v)^\top$ has comparable magnitude in both horizontal and vertical directions.

\subsubsection{Mean and variance}
\begin{figure}
    \centering
    \begin{subfigure}{0.4\linewidth}
    \includegraphics[trim=560pt 185pt 560pt 175pt, clip, width=1.0\linewidth]{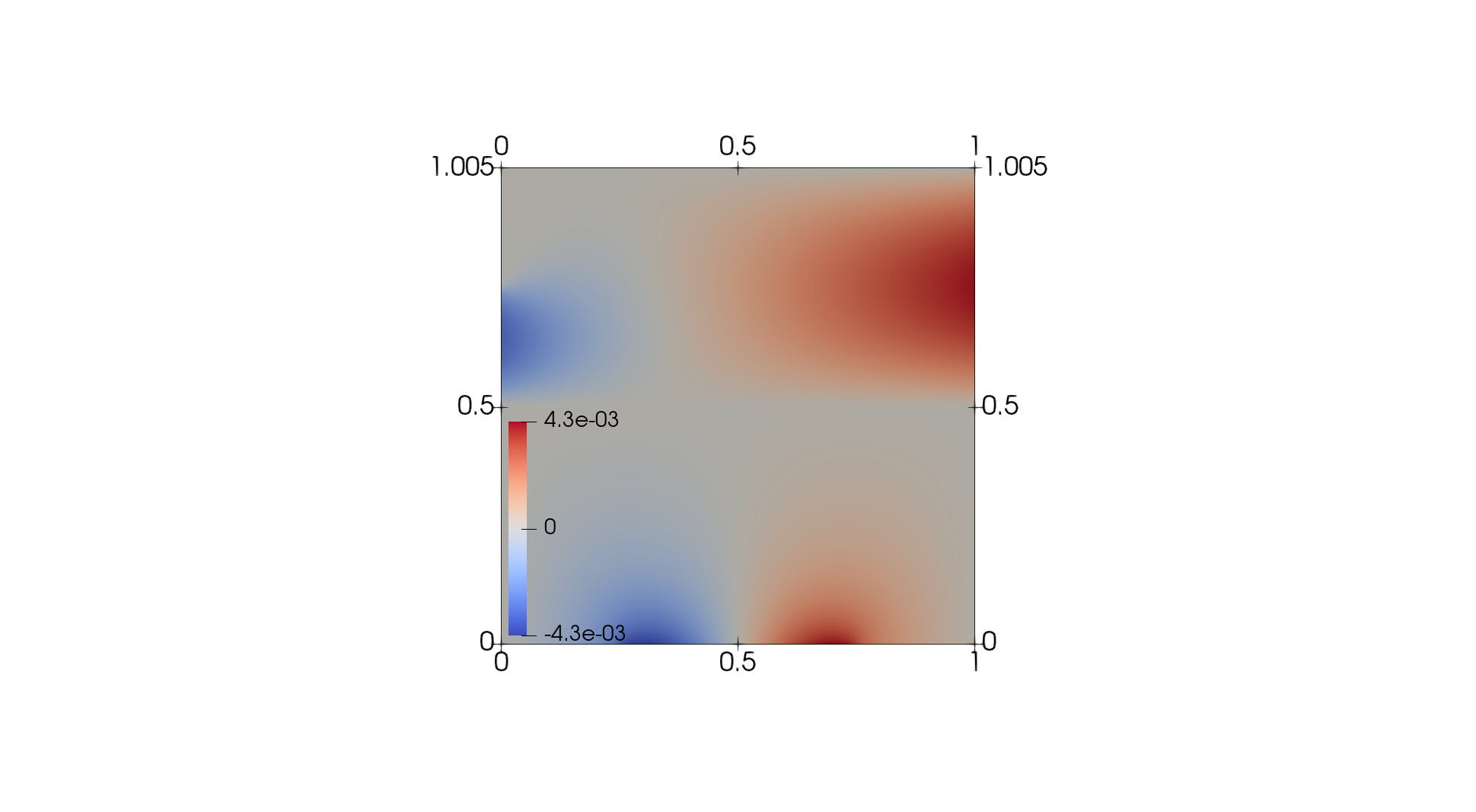}
    \caption{$\E[u]$}
    \end{subfigure}
    \hfil
    \begin{subfigure}{0.4\linewidth}
    \includegraphics[trim=560pt 185pt 560pt 175pt, clip, width=1.0\linewidth]{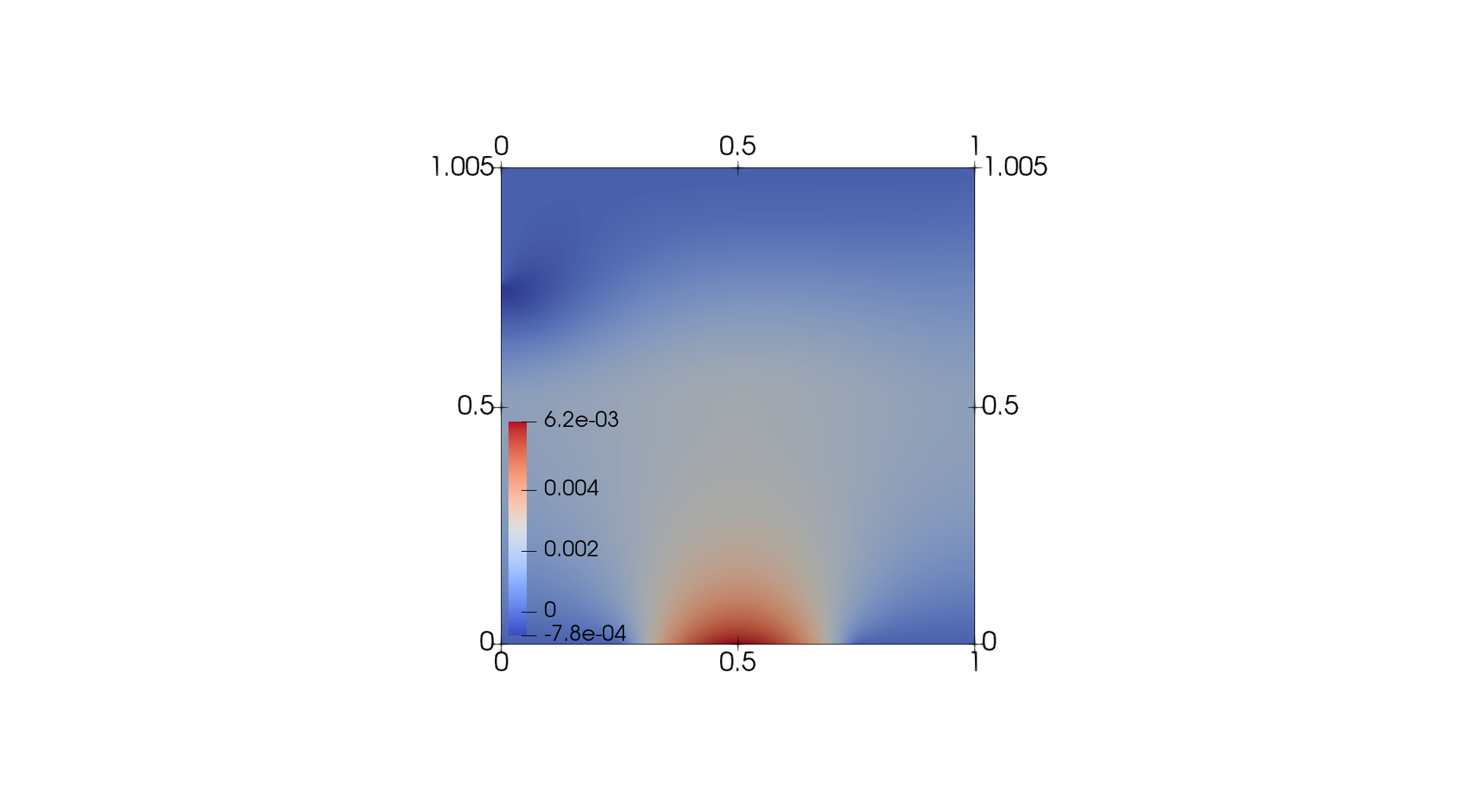}
    \caption{$\E[v]$}
    \end{subfigure}
    
    \begin{subfigure}{0.4\linewidth}
    \includegraphics[trim=560pt 185pt 560pt 175pt, clip, width=1.0\linewidth]{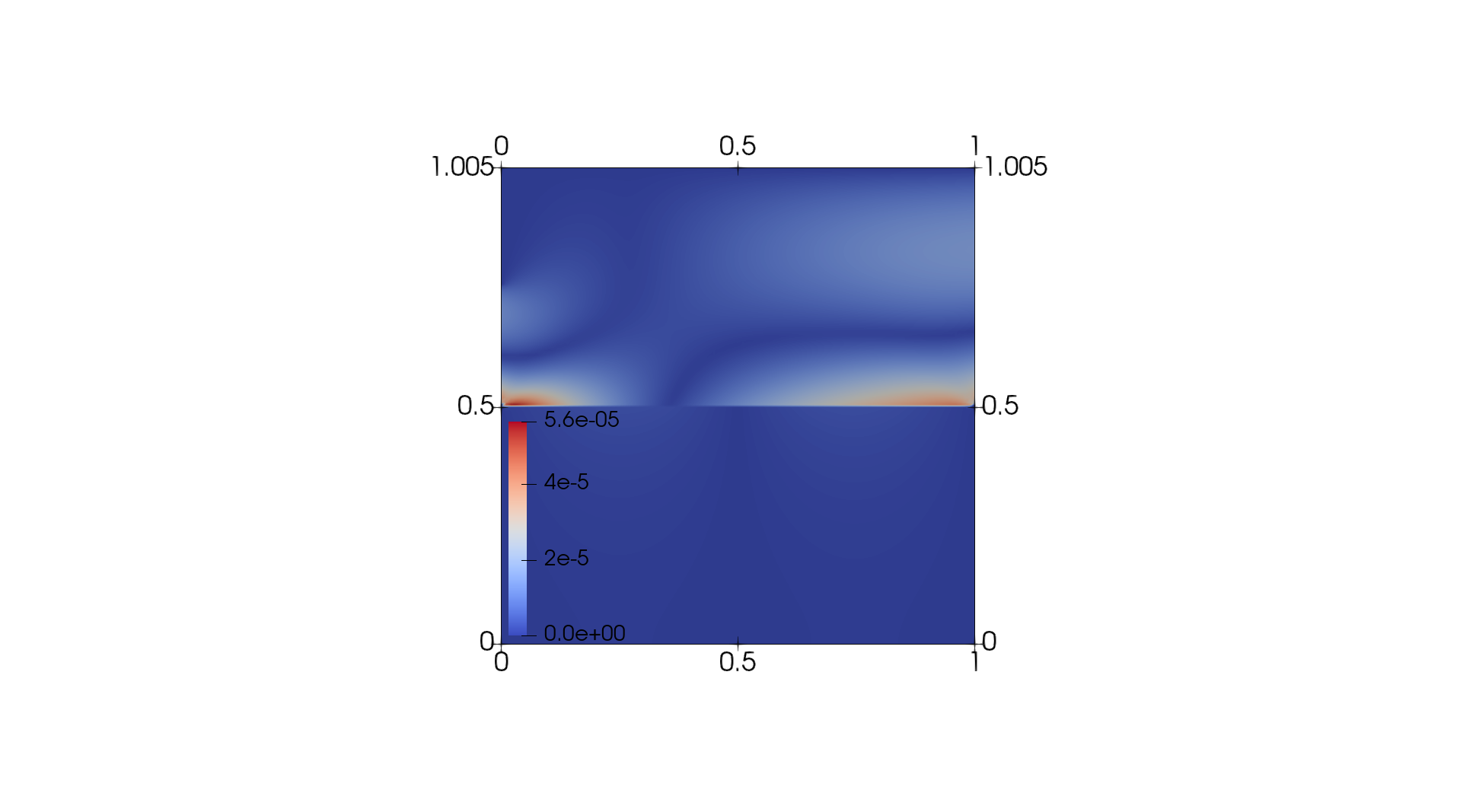}
    \caption{$\std[u]$}
    \end{subfigure}
    \hfil
    \begin{subfigure}{0.4\linewidth}
        \includegraphics[trim=560pt 185pt 560pt 175pt, clip, width=1.0\linewidth]{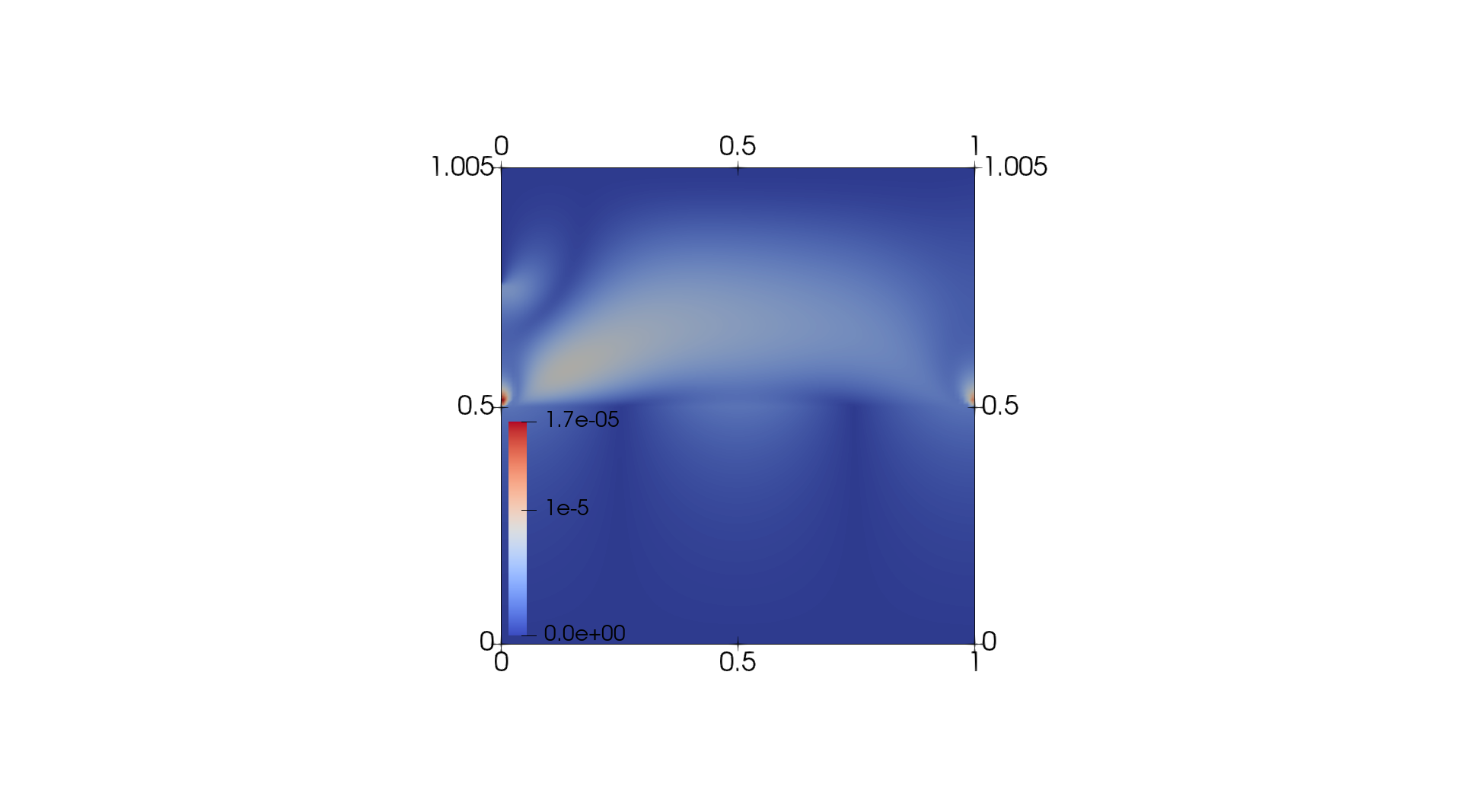}
    \caption{$\std[v]$}
    \end{subfigure}

    \begin{subfigure}{0.4\linewidth}
    \includegraphics[trim=560pt 185pt 560pt 175pt, clip, width=1.0\linewidth]{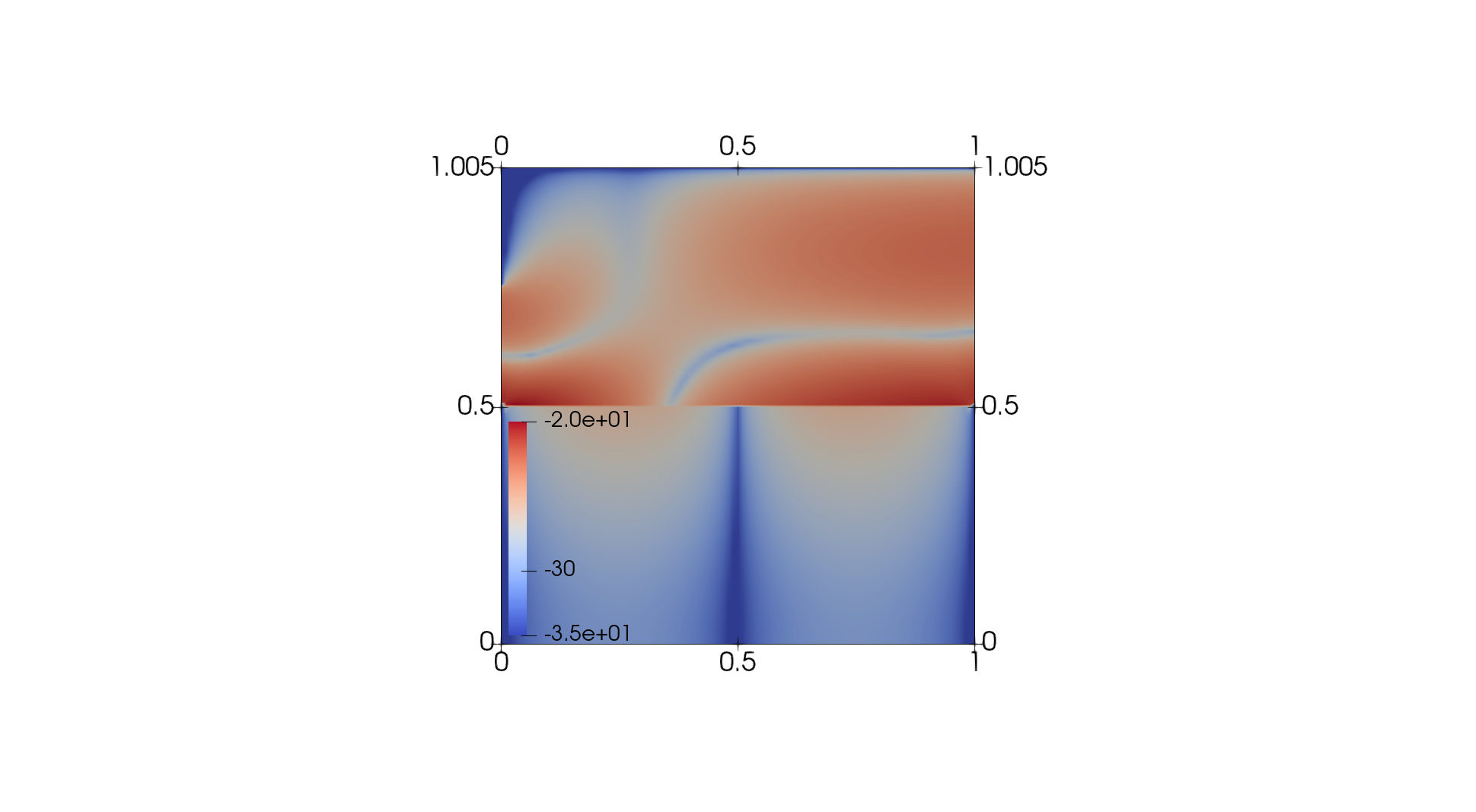}
    \caption{$\log\var[u]$}
    \end{subfigure}
    \hfil
    \begin{subfigure}{0.4\linewidth}
    \includegraphics[trim=560pt 185pt 560pt 175pt, clip, width=1.0\linewidth]{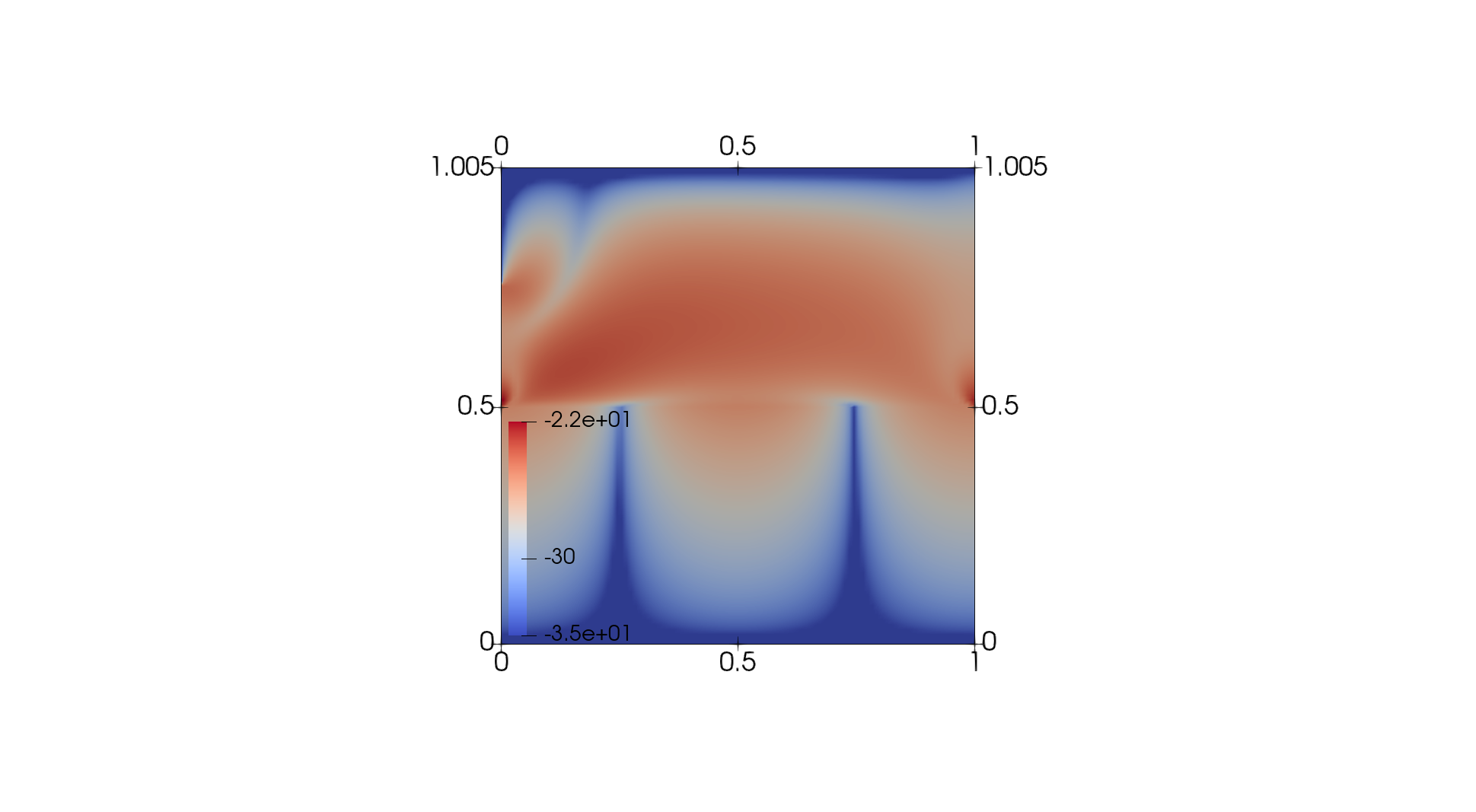}
    \caption{$\log\var[v]$}
    \end{subfigure}
    
    \caption{Mean, standard deviation and log-variance of the horizontal
    and vertical  velocities $(u,v)$
    computed from 50,000 MC samples.}
    \label{fig:mc_mean_var}
\end{figure}

\begin{figure}
    \centering
    \begin{subfigure}{0.4\linewidth}
        \includegraphics[trim=560pt 185pt 560pt 175pt, clip, width=1.0\linewidth]{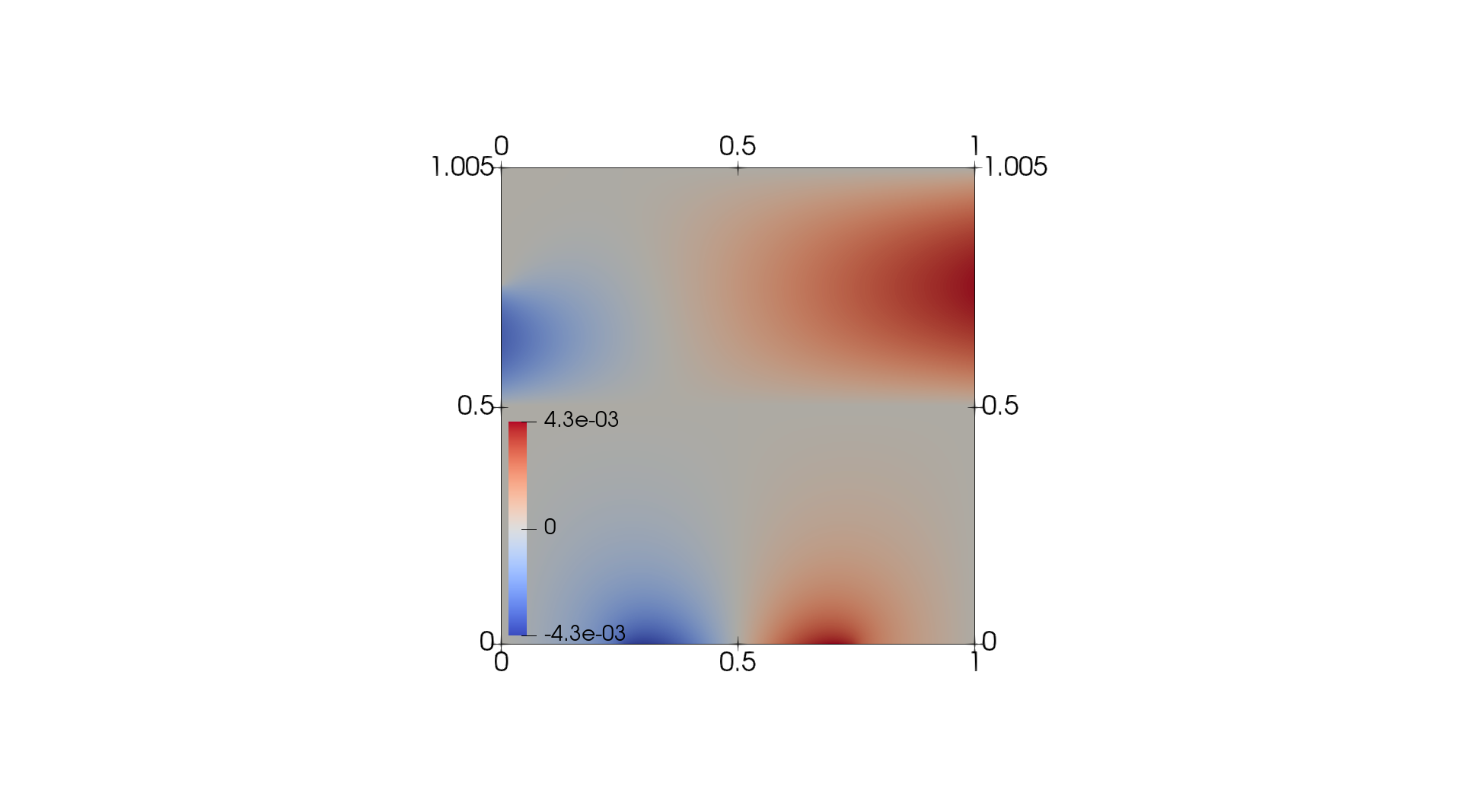}
    \caption{$\E[u]$}
    \end{subfigure}
    \hfil
    \begin{subfigure}{0.4\linewidth}
        \includegraphics[trim=560pt 185pt 560pt 175pt, clip, width=1.0\linewidth]{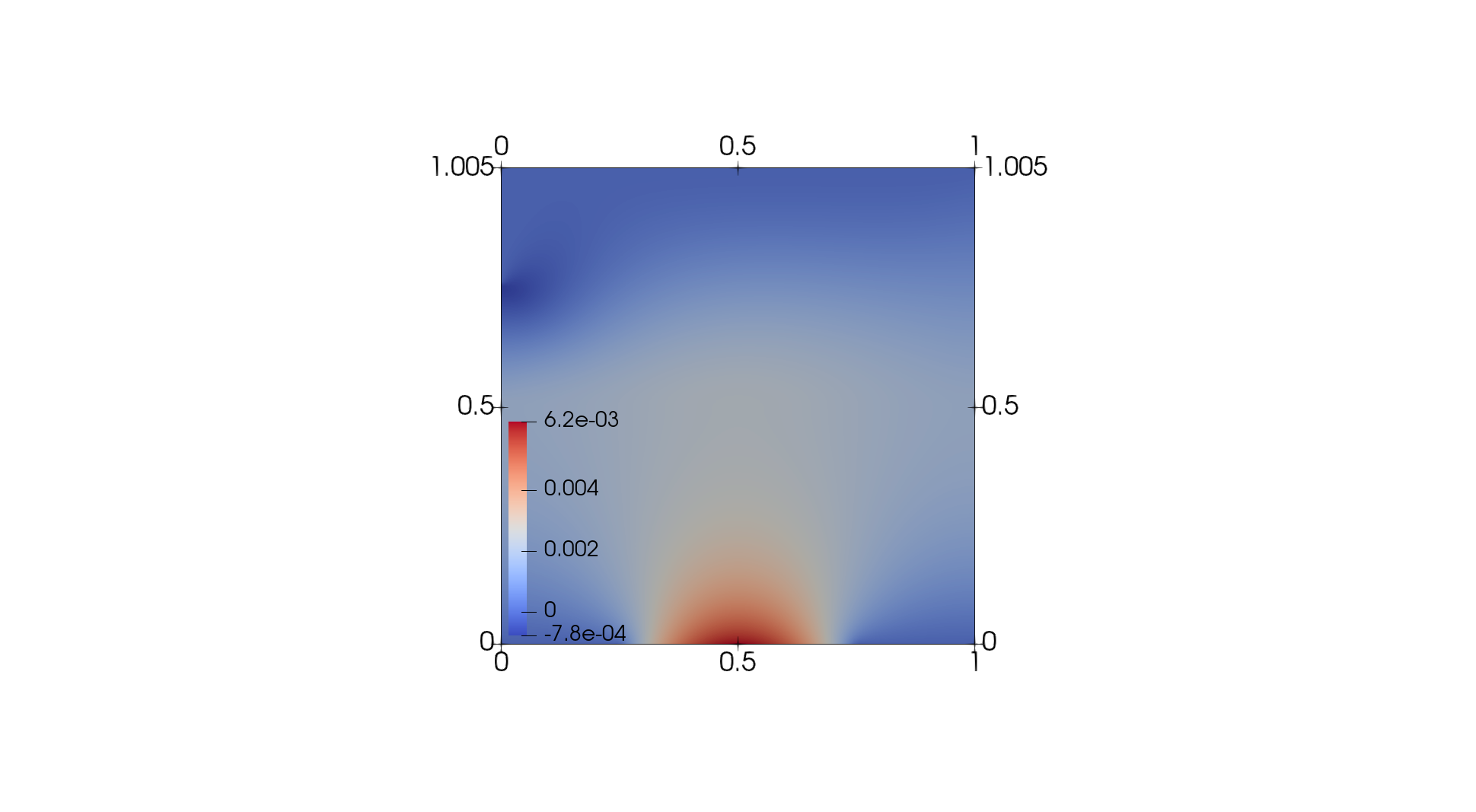}
    \caption{$\E[v]$}
    \end{subfigure}
    
    \begin{subfigure}{0.4\linewidth}
        \includegraphics[trim=560pt 185pt 560pt 175pt, clip, width=1.0\linewidth]{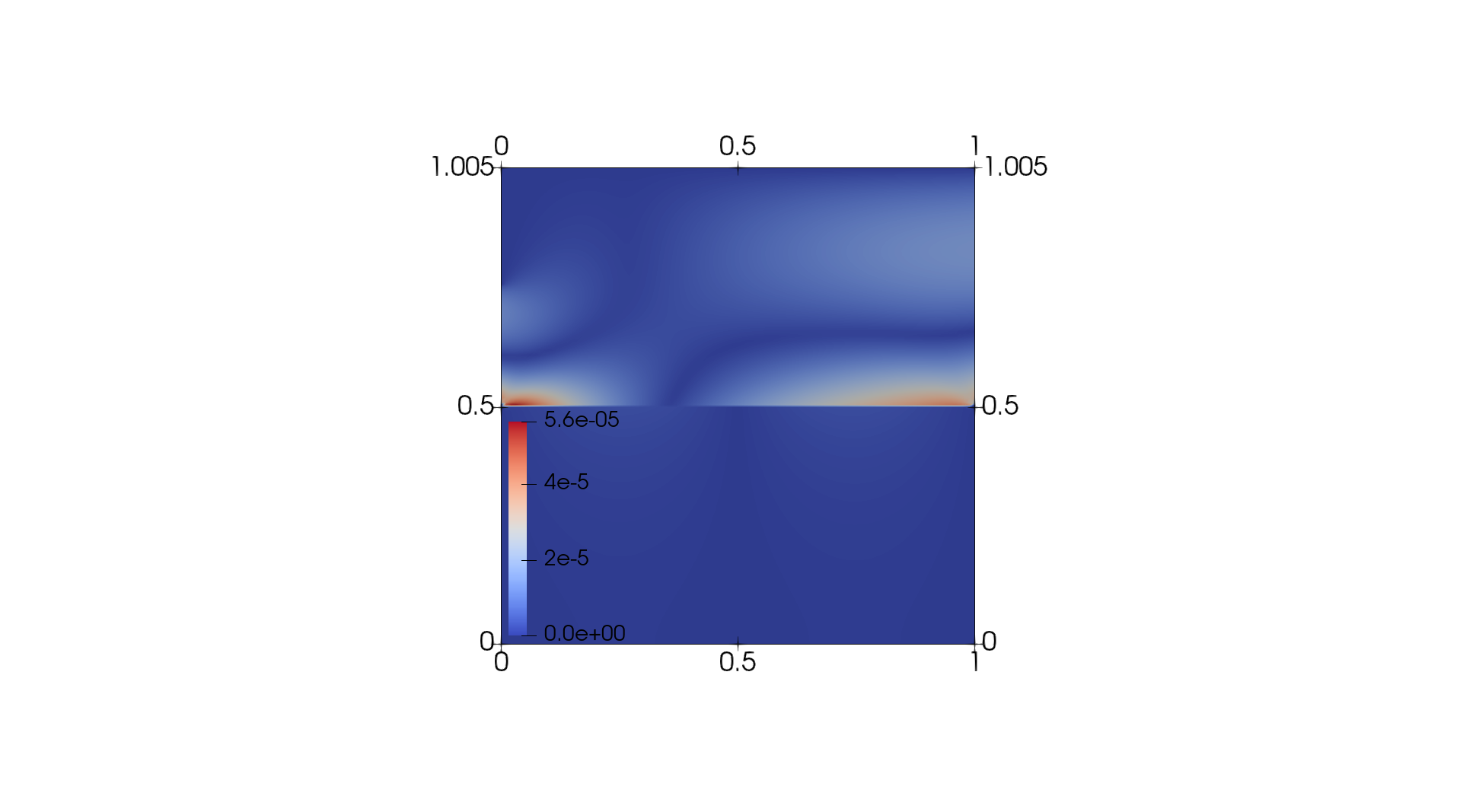}
    \caption{$\std[u]$}
    \end{subfigure}
    \hfil
    \begin{subfigure}{0.4\linewidth}
        \includegraphics[trim=560pt 185pt 560pt 175pt, clip, width=1.0\linewidth]{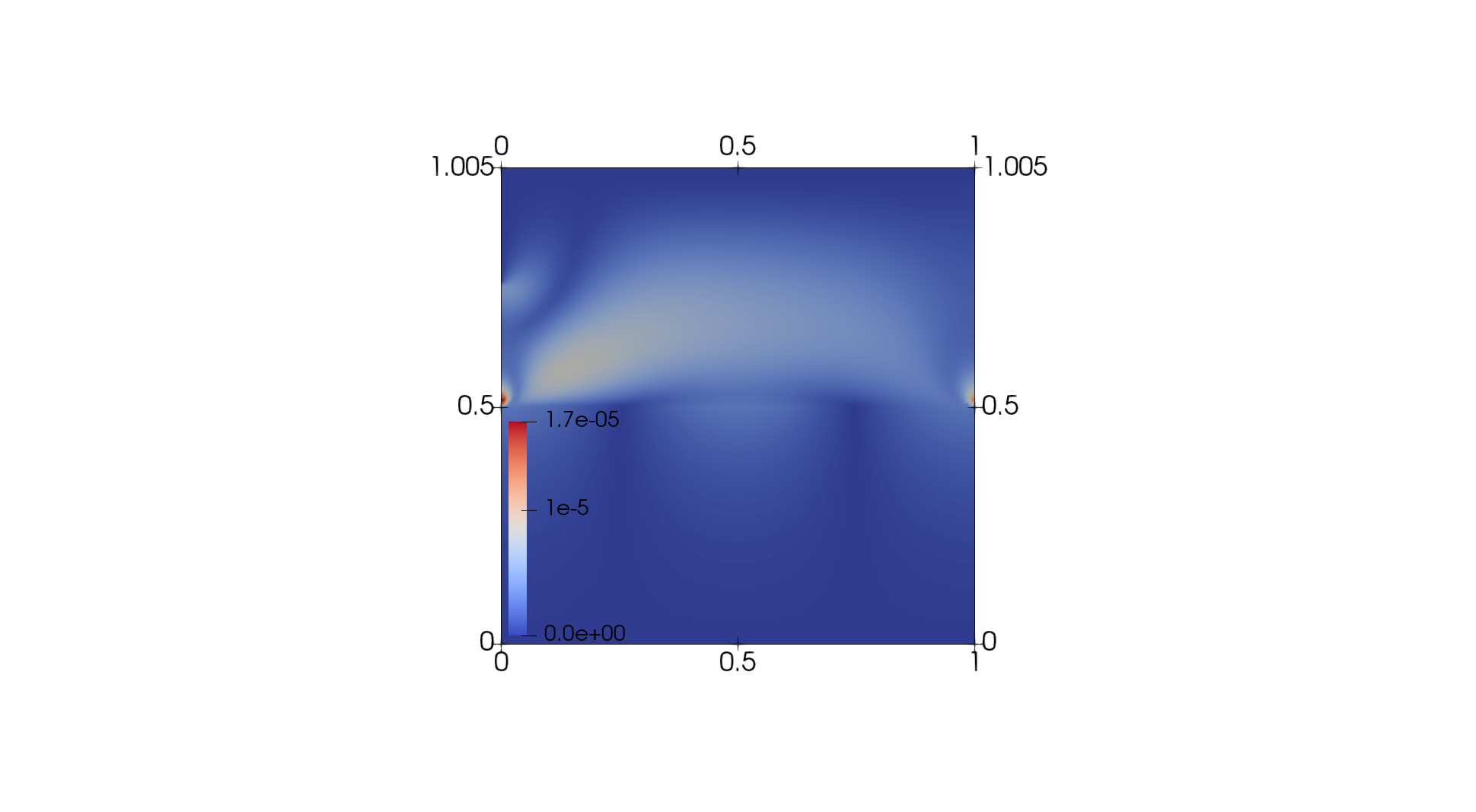}
    \caption{$\std[v]$}
    \end{subfigure}

    \begin{subfigure}{0.4\linewidth}
        \includegraphics[trim=560pt 185pt 560pt 175pt, clip, width=1.0\linewidth]{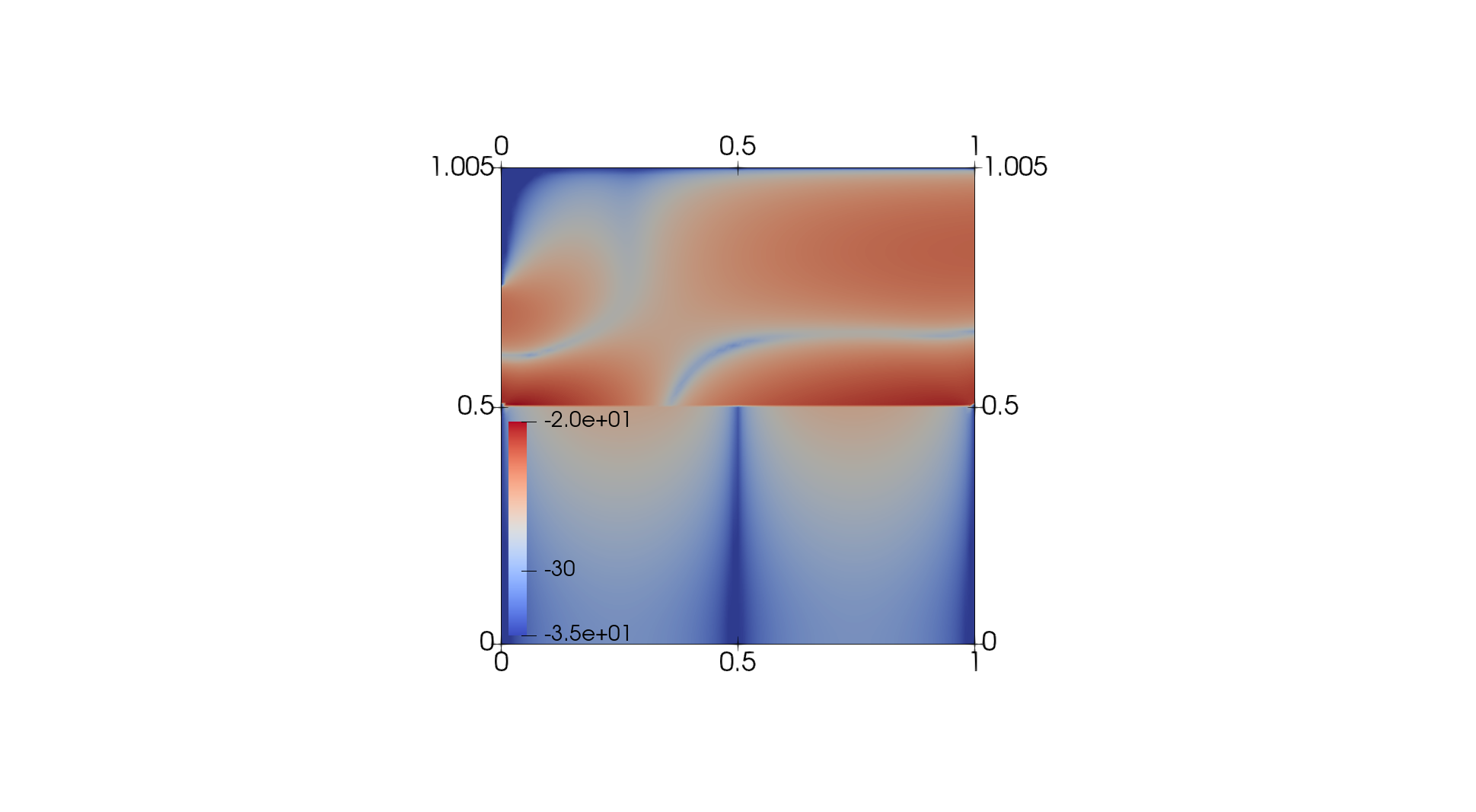}
    \caption{$\log\var[u]$}
    \end{subfigure}
    \hfil
    \begin{subfigure}{0.4\linewidth}
        \includegraphics[trim=560pt 185pt 560pt 175pt, clip, width=1.0\linewidth]{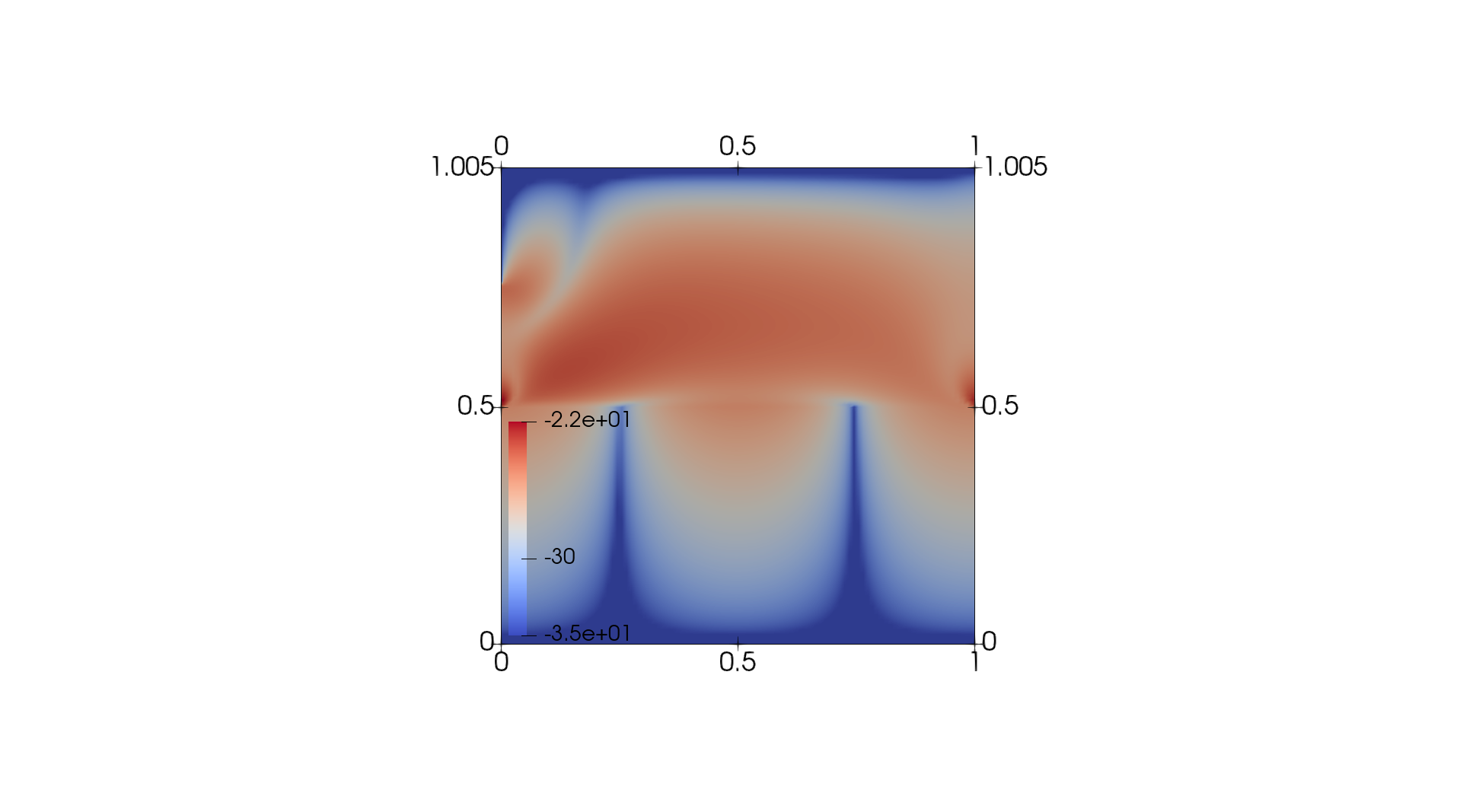}
    \caption{$\log\var[v]$}
    \end{subfigure}
    
    \caption{Mean, standard deviation and log-variance of the horizontal
    and vertical velocities $(u,v)$
    obtained from \aMR-surrogate with $\Nr=1, \No=2$,
    trained on 8192 QMC samples.}
    \label{fig:amr_mean_var}
\end{figure}

We use the mean $\E$, standard deviation $\std$, and log-variance $\log \var$ to characterize the random fields. Figure~\ref{fig:mc_mean_var} presents the corresponding statistical moments of the MC solution for horizontal and vertical velocity components, computed using $50{,}000$ realizations. Figure~\ref{fig:amr_mean_var} shows the mean, standard deviation, and log-variance obtained with the \aMR{} expansion of the polynomial degree $\No=2$ and refinement level $\Nr=1$. The logarithm of the variance enhances the visibility of regions with very small fluctuations that would otherwise be hidden in a linear scale. A comparison of the two figures demonstrates that the \aMR{} results reproduce the MC statistics with high fidelity. This agreement is particularly clear in the log-variance of the velocity components $u$ and $v$, which highlights regions of very small variance (Fig.~\ref{fig:mc_mean_var}(E,F), Fig.~\ref{fig:amr_mean_var}(E,F)). 

\subsubsection{Accuracy of the approximation}
To assess the accuracy of the approximation of the response of the physical models for our quantities of interest ($u,\,v$), we use RMSE and relative MSE. We evaluate both metrics across several (piecewise) polynomial-based surrogate models.
Despite their definitions being very similar, both metrics are relevant.
RMSE uses the same scale as the underlying quantities and is
therefore more relevant from a modeling point of view.
However, to compare different quantities of interest, the
relative MSE is preferable because of its normalization.

\begin{figure}
    \centering
    
    \begin{subfigure}{0.45\linewidth}
    \includegraphics[width=0.95\linewidth]{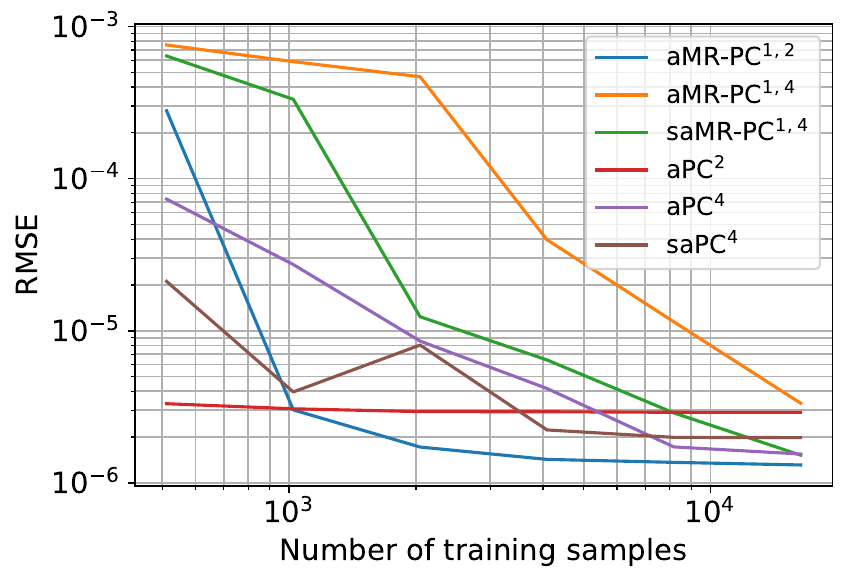}  
    \caption{RMSE of $u$}
    \end{subfigure}
    \hfil 
    \begin{subfigure}{0.45\linewidth}
    \includegraphics[width=0.95\linewidth]{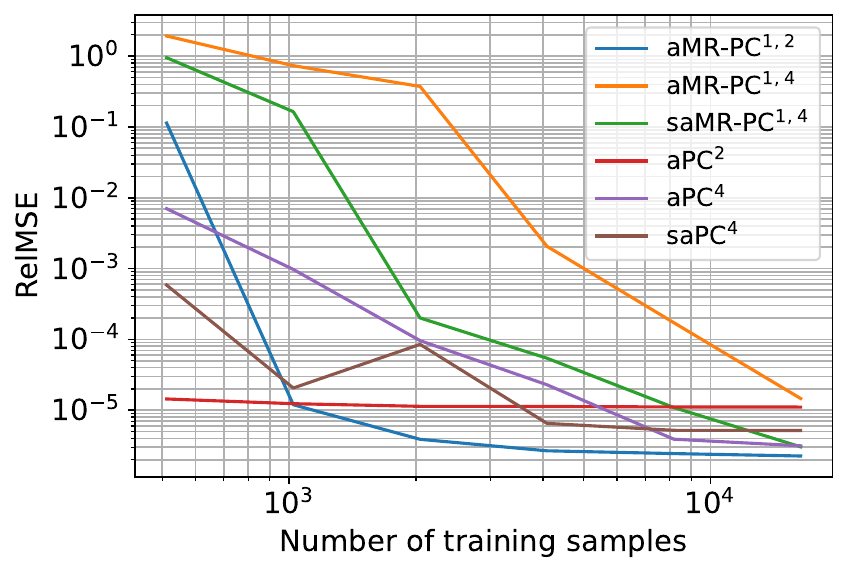}  
    \caption{Relative MSE of $u$}
    \end{subfigure}

    \begin{subfigure}{0.45\linewidth}
    \includegraphics[width=0.95\linewidth]{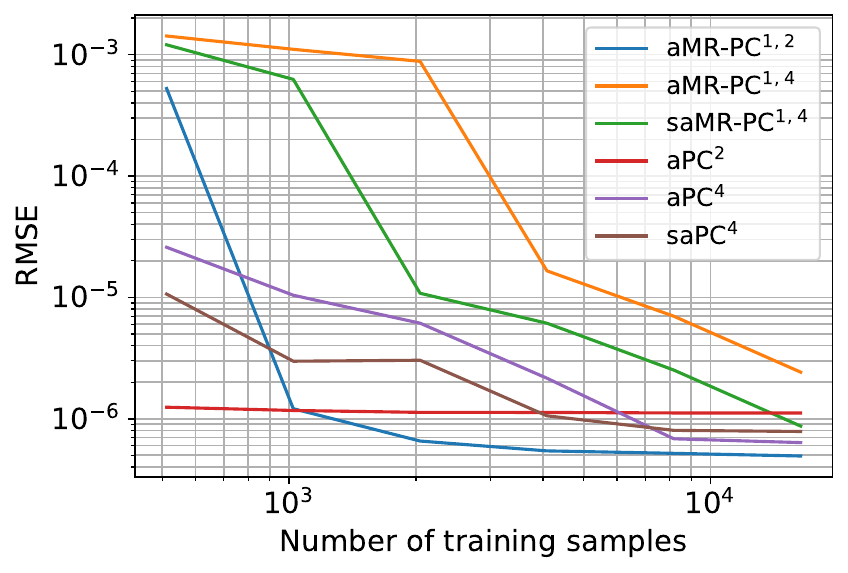}  
    \caption{RMSE of $v$}
    \end{subfigure}
    \hfil 
    \begin{subfigure}{0.45\linewidth}
    \includegraphics[width=0.95\linewidth]{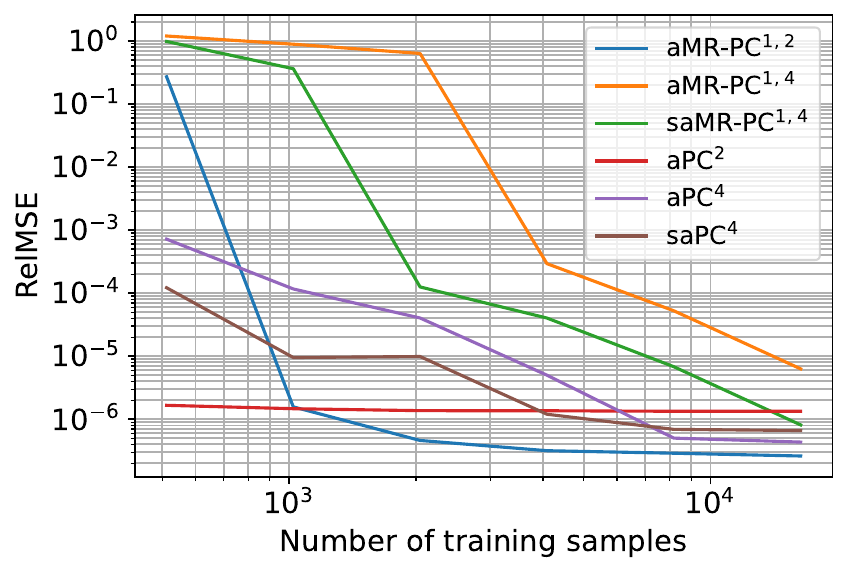}  
    \caption{Relative MSE of $v$}
    \end{subfigure}
    
    \caption{Case 1: RMSE and relative MSE error plots for $512,\ldots,16384$
    training samples.}
    \label{fig:error_i}
\end{figure}

Figure~\ref{fig:error_i} shows the above-mentioned errors evaluated over the whole space domain on 10,000 randomly selected samples.
The $x$-axis is labeled in terms of $2^9,\ldots,2^{14}$ samples used for the training of
the surrogate models. 
Here, we use the following notation: aPC$^{\No}$
is aPC with the highest polynomial degree $\No$,
$\aMR^{\Nr,\No}$ is
\aMR{} with the number of refinements $\Nr$ and the highest polynomial degree $\No$.
Further, saPC and s\aMR{} denote the sparse versions of the surrogate models
using hyperbolic truncation with $q=0.75$.

In almost all cases, we observe that the errors decrease with the increasing number of training samples. In Fig.~\ref{fig:error_i}, we can see the expected relation between the number of degrees of freedom, the number of required training samples, and the achieved accuracy. 
More precisely, aPC with ${\No=2}$ obtains its almost highest accuracy with 512 training samples. 
The \aMR{} with $\Nr=1$ and low polynomial degree $\No=2$ provides the most
accurate approximation for $u$ and $v$ starting with $2048$ training samples.
It can also be seen that the convergence is faster if hyperbolic truncation-based sparsity is used.

Since we intend to use the expansion coefficients rather than direct surrogate responses for calculation of the Sobol' coefficients,
we also evaluate the accuracy of the approximation of mean and standard
deviation of the quantities of interest by using the expansion coefficients
given in Eqs.~\eqref{eq:amr-mean} and \eqref{eq:amr-var}.
The $\Lp2$-error of the approximation of the mean calculated in
the whole spatial domain is provided in Tab.~\ref{tab:error-means_i}.
Evaluation of the $\Lp2$-error of the approximation of the standard
deviation is provided in Tab.~\ref{tab:error-std_i}.
In both cases, we use the mean and standard deviation provided by MC
computed on 50,000 samples as a reference solution.
For both velocity components $u$ and $v$, aPC with $\No=2$ and $\aMR$ for $\Nr=1,\, \No=2$ (as the number of training samples increases) outperform other surrogate models in accuracy. 
However, the accuracy of other surrogate models also increases with the increasing number
of training samples. 

\begin{table}
\caption{Case 1: $\Lp2$-error of the approximation of the mean $(\E)$.}
\label{tab:error-means_i}
\begin{flushleft}

\begin{subtable}{0.95\linewidth}
\caption{Approximation of the mean of the horizontal velocity
component $\E[u]$.}
\begin{tabular}{|r|c|c|c|c|c|c|}
\hline
Samples & $N_o=2$ & $N_o=4$ & $N_o=4 / q=0.75$ & $N_r=1 / N_o=2$ & $N_r=1 / N_o=4$ & $N_r=1 / N_o=4 / q=0.75$ \\
\hline
512 & 1.31e-07 & 4.08e-06 & 1.42e-06 & 1.30e-04 & 9.16e-04 & 6.46e-04 \\
1024 & 5.52e-08 & 1.44e-06 & 2.70e-07 & 9.50e-08 & 5.07e-04 & 1.78e-04 \\
2048 & 1.28e-08 & 5.73e-07 & 5.59e-07 & 4.78e-08 & 3.58e-04 & 9.28e-07 \\
4096 & 7.88e-09 & 2.57e-07 & 4.62e-08 & 2.57e-08 & 3.94e-06 & 3.76e-07 \\
8192 & 1.22e-08 & 5.96e-08 & 2.51e-08 & 9.67e-09 & 6.38e-07 & 1.54e-07 \\
16384 & 1.18e-08 & 2.42e-08 & 7.04e-09 & 5.47e-09 & 2.40e-07 & 6.33e-08 \\
\hline
\end{tabular}
\end{subtable}\\[4mm]

\begin{subtable}{0.95\linewidth}
\caption{Approximation of the mean of the vertical velocity component $\E[v]$.}
\begin{tabular}{|r|c|c|c|c|c|c|}
\hline
Samples & $N_o=2$ & $N_o=4$ & $N_o=4 / q=0.75$ & $N_r=1 / N_o=2$ & $N_r=1 / N_o=4$ & $N_r=1 / N_o=4 / q=0.75$ \\
\hline
512 & 4.32e-08 & 1.45e-06 & 6.80e-07 & 1.76e-04 & 1.24e-03 & 8.72e-04 \\
1024 & 1.87e-08 & 6.27e-07 & 1.97e-07 & 3.25e-08 & 6.84e-04 & 2.40e-04 \\
2048 & 6.70e-09 & 4.06e-07 & 2.07e-07 & 1.70e-08 & 4.83e-04 & 4.53e-07 \\
4096 & 6.67e-09 & 1.44e-07 & 5.15e-08 & 8.92e-09 & 1.32e-06 & 1.61e-07 \\
8192 & 6.63e-09 & 3.61e-08 & 2.08e-08 & 5.27e-09 & 3.00e-07 & 8.81e-08 \\
16384 & 5.96e-09 & 1.69e-08 & 9.47e-09 & 4.52e-09 & 9.72e-08 & 3.02e-08 \\
\hline
\end{tabular}

\end{subtable}
\end{flushleft}
\end{table}

\begin{table}
\caption{Case 1: $\Lp2$-error of the approximation of the
standard deviation $(\std)$.}
\label{tab:error-std_i}
\begin{flushleft}

\begin{subtable}{0.95\linewidth}
\caption{Approximation of the standard deviation of the
horizontal velocity component $\std[u]$.}
\begin{tabular}{|r|c|c|c|c|c|c|}
\hline
Samples & $N_o=2$ & $N_o=4$ & $N_o=4 / q=0.75$ & $N_r=1 / N_o=2$ & $N_r=1 / N_o=4$ & $N_r=1 / N_o=4 / q=0.75$ \\
\hline
512 & 2.77e-07 & 4.01e-04 & 1.52e-04 & 3.95e-04 & 6.27e-04 & 6.71e-04 \\
1024 & 3.88e-07 & 1.52e-04 & 3.10e-05 & 1.62e-06 & 6.50e-04 & 4.53e-04 \\
2048 & 6.71e-07 & 6.88e-05 & 6.37e-05 & 1.07e-07 & 5.92e-04 & 3.29e-04 \\
4096 & 6.42e-07 & 3.21e-05 & 3.91e-06 & 1.10e-07 & 6.61e-04 & 7.48e-05 \\
8192 & 6.54e-07 & 8.84e-06 & 2.90e-06 & 7.49e-08 & 8.52e-05 & 3.05e-05 \\
16384 & 6.56e-07 & 3.39e-06 & 8.34e-07 & 8.74e-08 & 3.67e-05 & 1.03e-05 \\
\hline
\end{tabular}

\end{subtable}\\[4mm]

\begin{subtable}{0.95\linewidth}
\caption{Approximation of the standard deviation of the vertical
velocity component $\std[v]$.}
\begin{tabular}{|r|c|c|c|c|c|c|c|}
\hline
Samples & $N_o=2$ & $N_o=4$ & $N_o=4 / q=0.75$ & $N_r=1 / N_o=2$ & $N_r=1 / N_o=4$ & $N_r=1 / N_o=4 / q=0.75$ \\
\hline
512 & 2.13e-07 & 1.50e-04 & 7.32e-05 & 5.37e-04 & 8.50e-04 & 9.09e-04 \\
1024 & 1.98e-07 & 7.02e-05 & 2.32e-05 & 6.26e-07 & 8.81e-04 & 6.15e-04 \\
2048 & 2.66e-07 & 4.72e-05 & 2.43e-05 & 5.98e-08 & 8.03e-04 & 1.33e-04 \\
4096 & 2.36e-07 & 1.78e-05 & 5.72e-06 & 5.43e-08 & 2.53e-04 & 3.08e-05 \\
8192 & 2.52e-07 & 5.07e-06 & 2.26e-06 & 3.96e-08 & 3.99e-05 & 1.41e-05 \\
16384 & 2.61e-07 & 2.05e-06 & 6.99e-07 & 3.81e-08 & 1.44e-05 & 4.49e-06 \\
\hline
\end{tabular}
\end{subtable}
\end{flushleft}
\end{table}

\subsubsection{GSA using Sobol' indices}\label{sec:totalSobolindicies_i}
As we have seen in the previous section, \aMR{} with $\Nr=1,\,\No=2$
demonstrates the best performance in terms of accuracy.
Therefore, we focus on this surrogate model for the GSA
of the uncertain parameters. 
To assess the reliability of the estimated sensitivity indices, we examine how the space-averaged total sensitivity indices, covering the entire spatial domain, change as the number of training samples increases, as shown in Fig~\ref{fig:tot_evo_u},\,\ref{fig:tot_evo_v}.
To avoid the influence of physically not meaningful results,
we excluded the low variance $\var<10^{-10}$ from the evaluation.
In all considered test series, we can observe that a reasonable approach of the final values was achieved with 2048 training samples.
Changes in total sensitivity indices calculated from the expansion
coefficients of the \aMR{} surrogate models obtained from more
than 4096 training samples are only minimal.
Moreover, when we compare the total sensitivity indices shown in Tab.~\ref{tab:Sobol-total} with the Sobol' sensitivity indices in Tab.~\ref{tab:Sobol-u}, \ref{tab:Sobol-v}, we notice that the trends in results derived from aPC with $\No=2$ are similar. Consequently, these trends indirectly validate the findings derived from the \aMR{} surrogate model with $\Nr=1,\,\No=2$.
Hyperbolic truncation affects not only intermediate Sobol' terms but also the total indices computed from truncated expansions, leading to unreliable behavior.

\begin{figure}
    \centering
    \begin{subfigure}{0.4\linewidth}
        \includegraphics[trim=560pt 185pt 560pt 175pt, clip, width=1.0\linewidth]{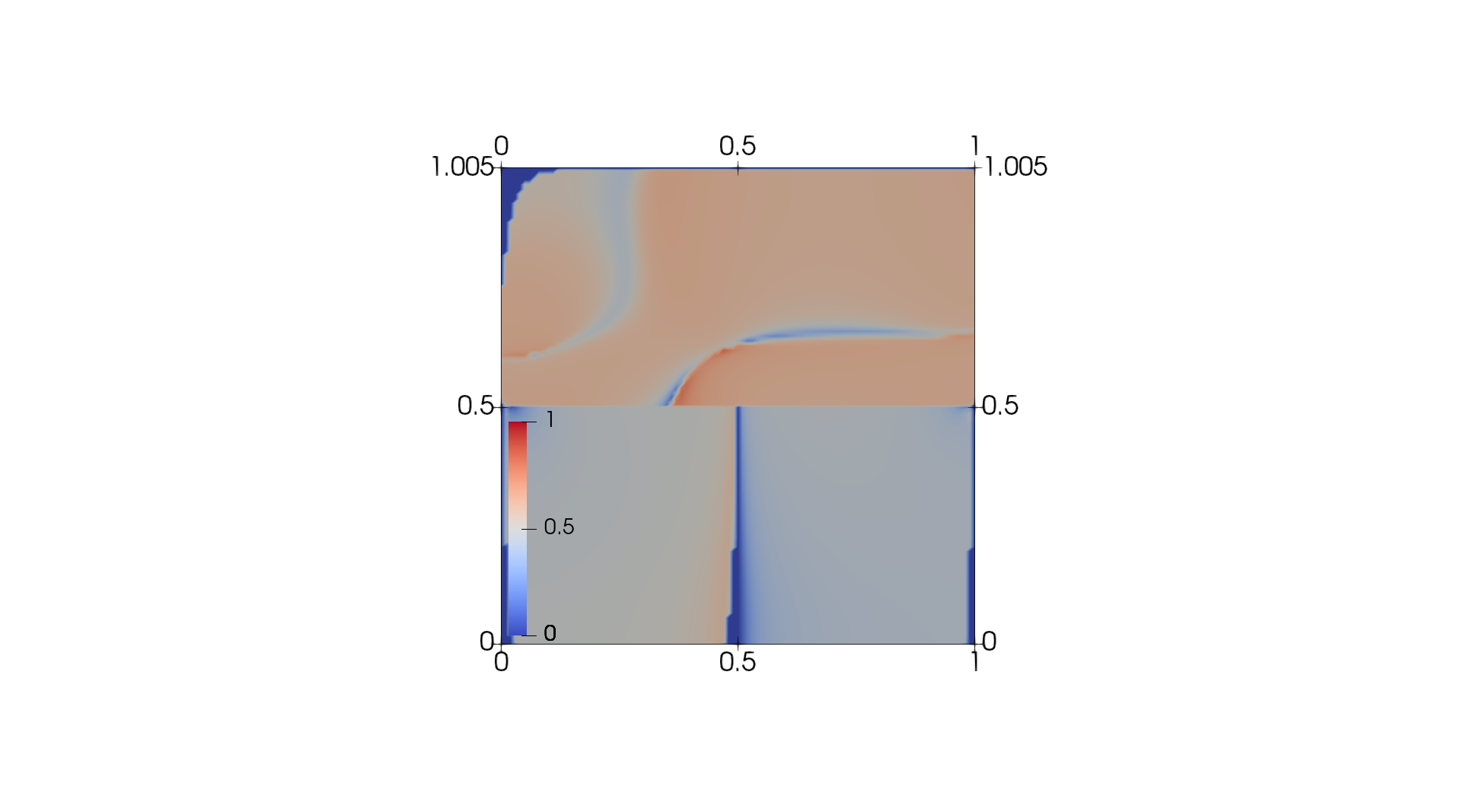}
    \caption{Velocity $u$, input parameter $I=1$}
    \end{subfigure}
    \hfil
    \begin{subfigure}{0.4\linewidth}
        \includegraphics[trim=560pt 185pt 560pt 175pt, clip, width=1.0\linewidth]{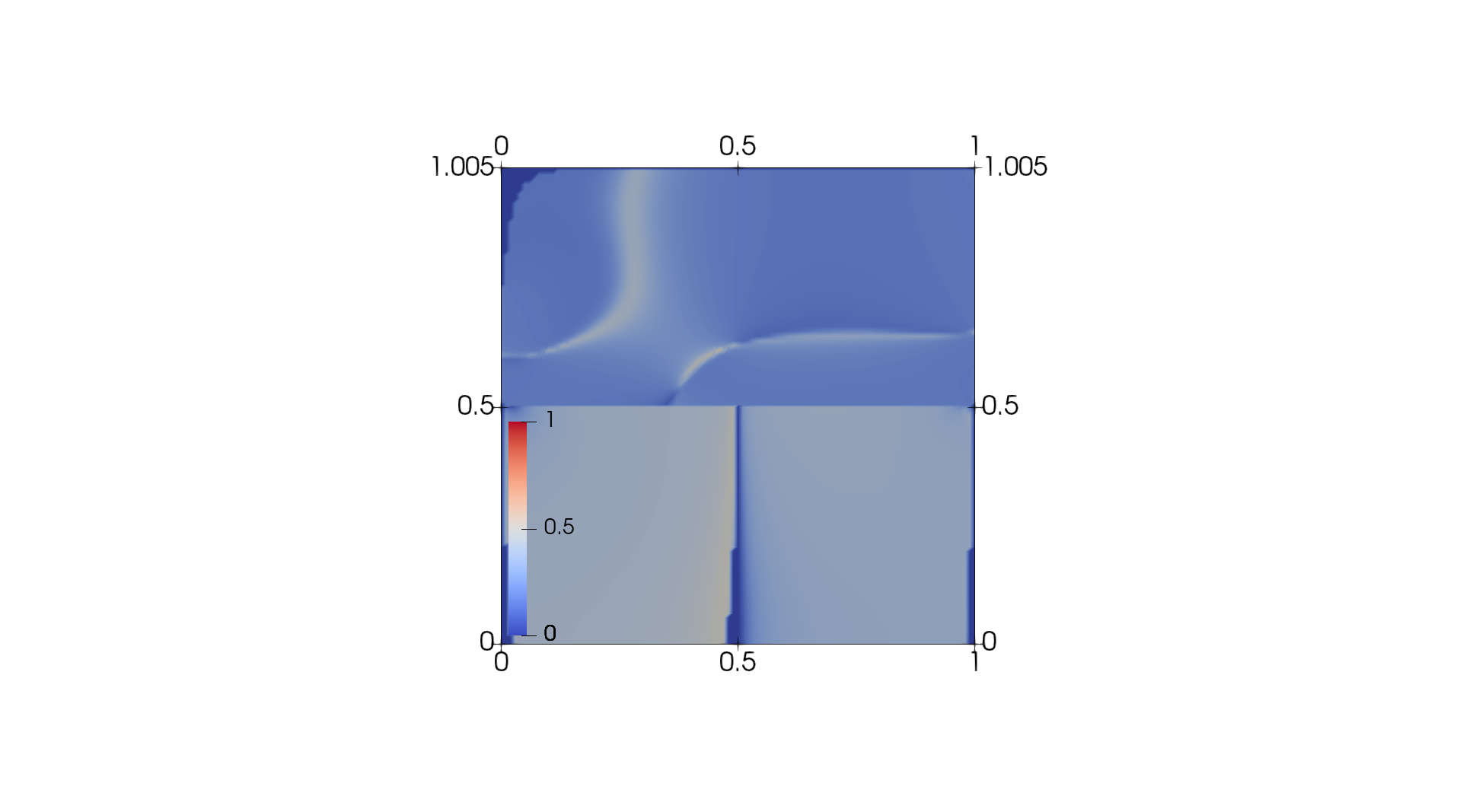}
    \caption{Velocity $u$, input parameter $I=2$}
    \end{subfigure}\\[2mm]
    
    \begin{subfigure}{0.4\linewidth}
        \includegraphics[trim=560pt 185pt 560pt 175pt, clip, width=1.0\linewidth]{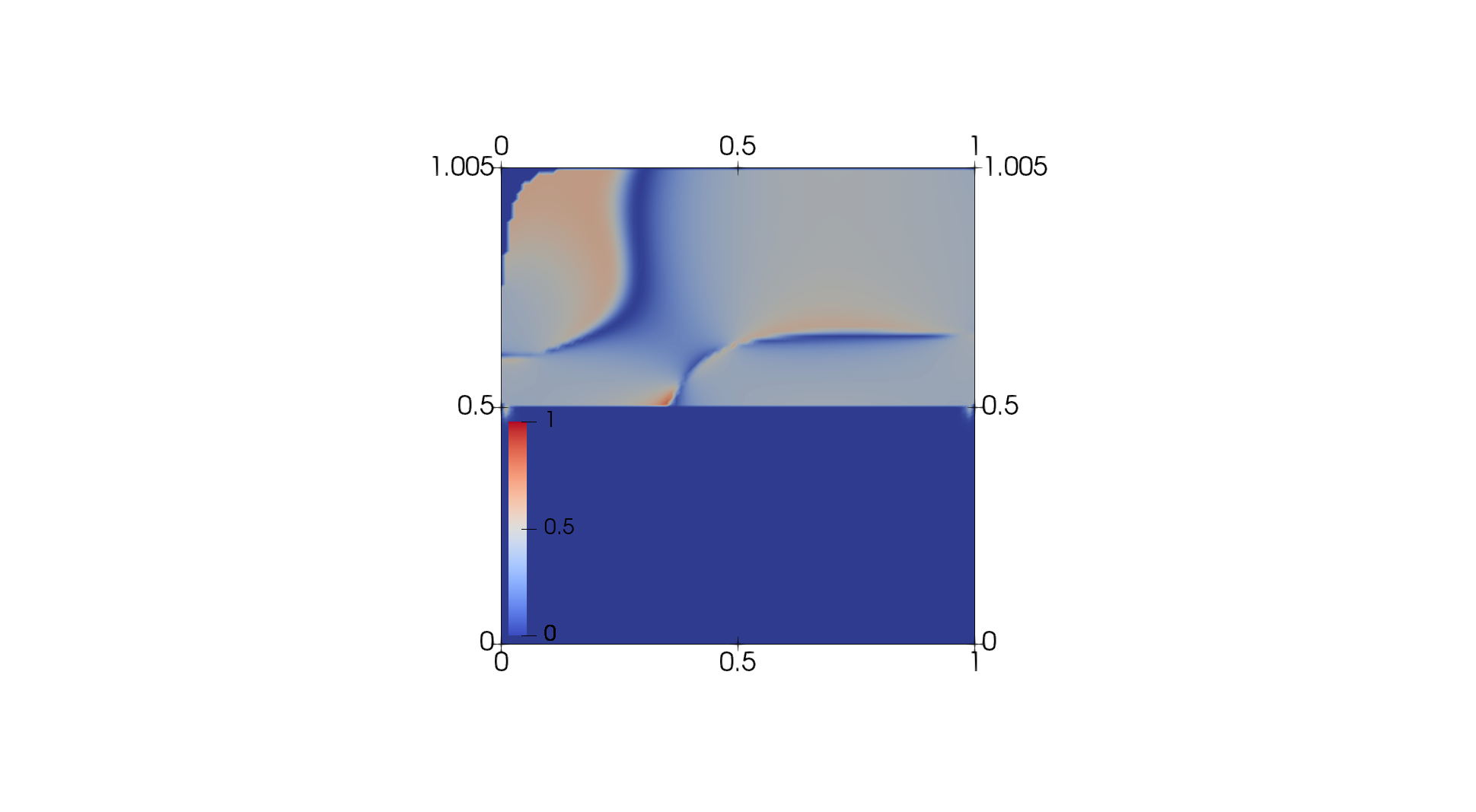}
    \caption{Velocity $u$, input parameter $I=3$}
    \end{subfigure}
    \hfil
    \begin{subfigure}{0.4\linewidth}
        \includegraphics[trim=560pt 185pt 560pt 175pt, clip, width=1.0\linewidth]{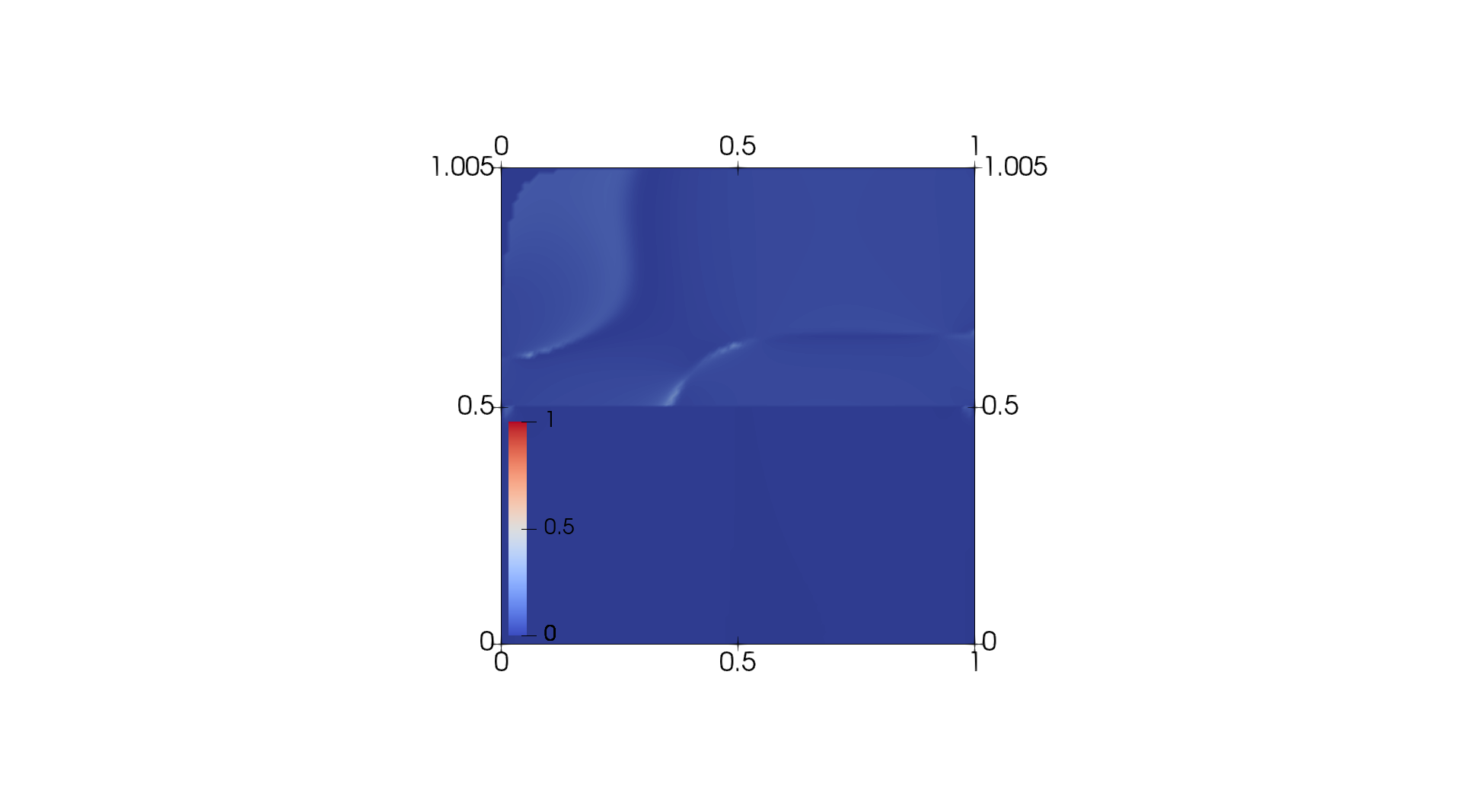}
    \caption{Velocity $u$, input parameter $I=4$}
    \end{subfigure}\\[2mm]

    \begin{subfigure}{0.4\linewidth}
        \includegraphics[trim=560pt 185pt 560pt 175pt, clip, width=1.0\linewidth]{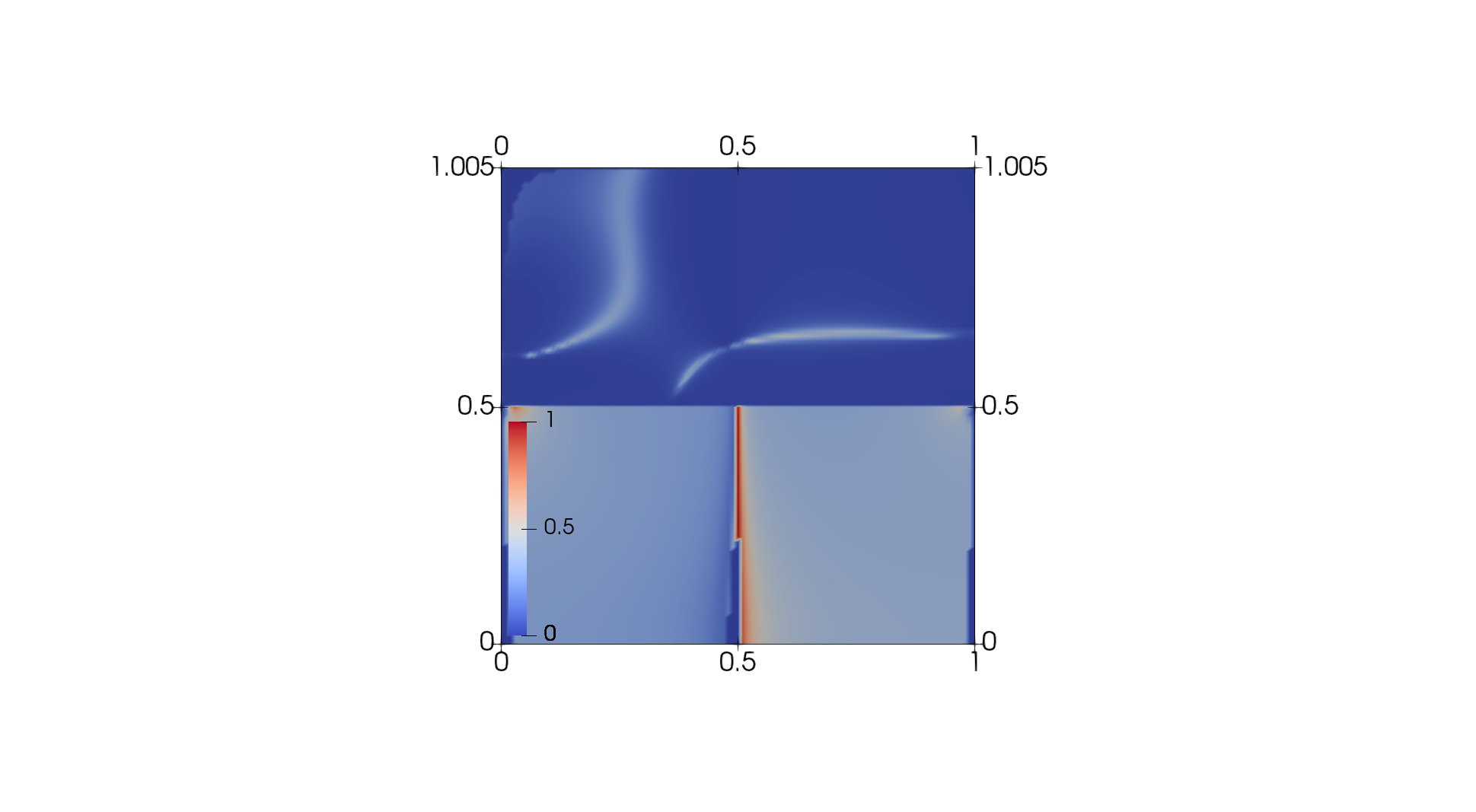}
    \caption{Velocity $u$, input parameter $I=5$ }
    \end{subfigure}
   
    \caption{Case 1: Total Sobol' indices 
    of the horizontal velocity $u$
    obtained from \aMR-surrogate with $\Nr=1, \No=2$,
    trained on 8192 QMC samples.}
    \label{fig:amr_sobol-u}
\end{figure}

\begin{figure}
    \centering
    \begin{subfigure}{0.4\linewidth}
        \includegraphics[trim=560pt 185pt 560pt 175pt, clip, width=1.0\linewidth]{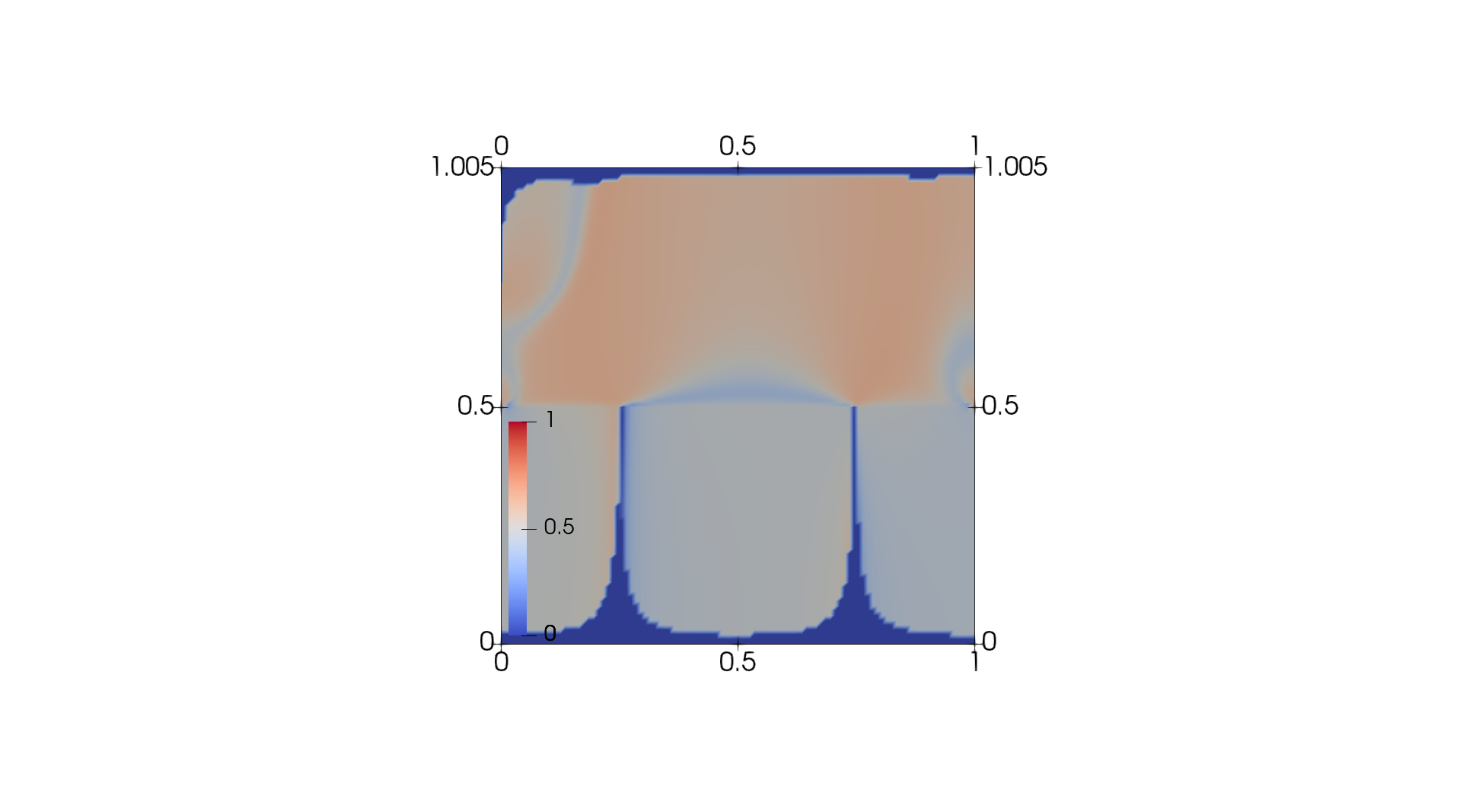}
    \caption{Velocity $v$, input parameter $I=1$}
    \end{subfigure}
    \hfil
    \begin{subfigure}{0.4\linewidth}
        \includegraphics[trim=560pt 185pt 560pt 175pt, clip, width=1.0\linewidth]{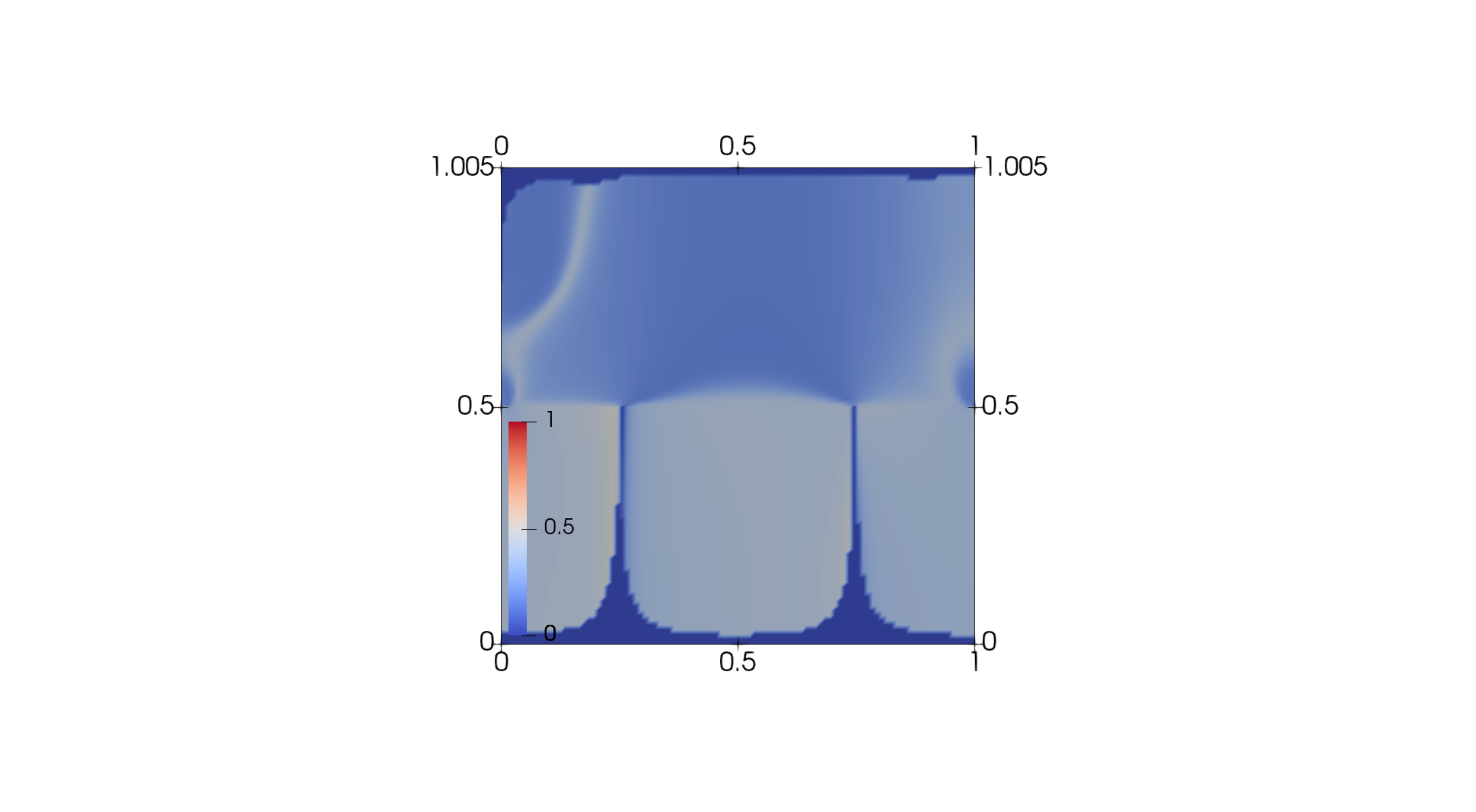}
    \caption{Velocity $v$, input parameter $I=2$}
    \end{subfigure}\\[2mm]
    
    \begin{subfigure}{0.4\linewidth}
        \includegraphics[trim=560pt 185pt 560pt 175pt, clip, width=1.0\linewidth]{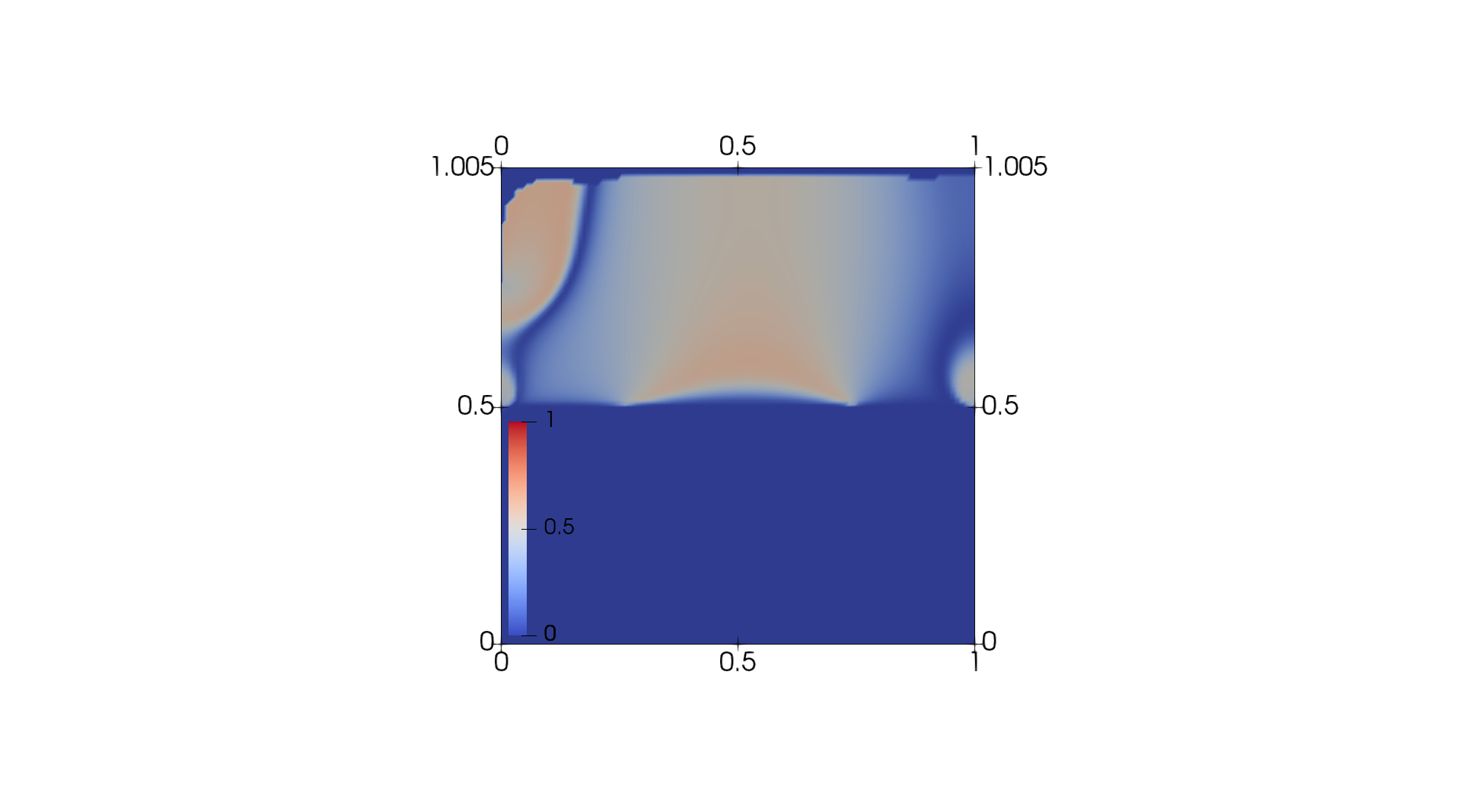}
    \caption{Velocity $v$, input parameter $I=3$}
    \end{subfigure}
    \hfil
    \begin{subfigure}{0.4\linewidth}
        \includegraphics[trim=560pt 185pt 560pt 175pt, clip, width=1.0\linewidth]{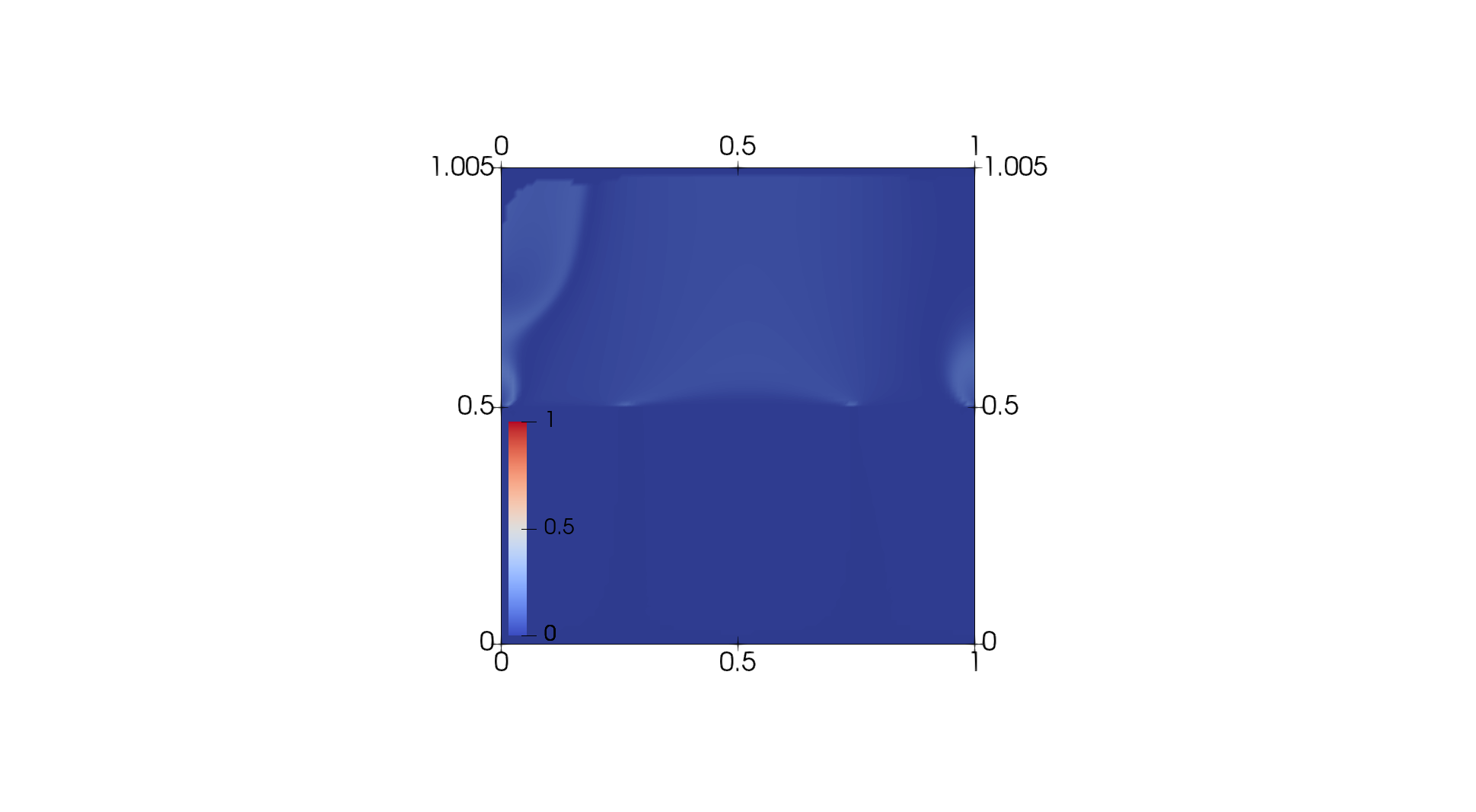}
    \caption{Velocity $v$, input parameter $I=4$}
    \end{subfigure}\\[2mm]

    \begin{subfigure}{0.4\linewidth}
        \includegraphics[trim=560pt 185pt 560pt 175pt, clip, width=1.0\linewidth]{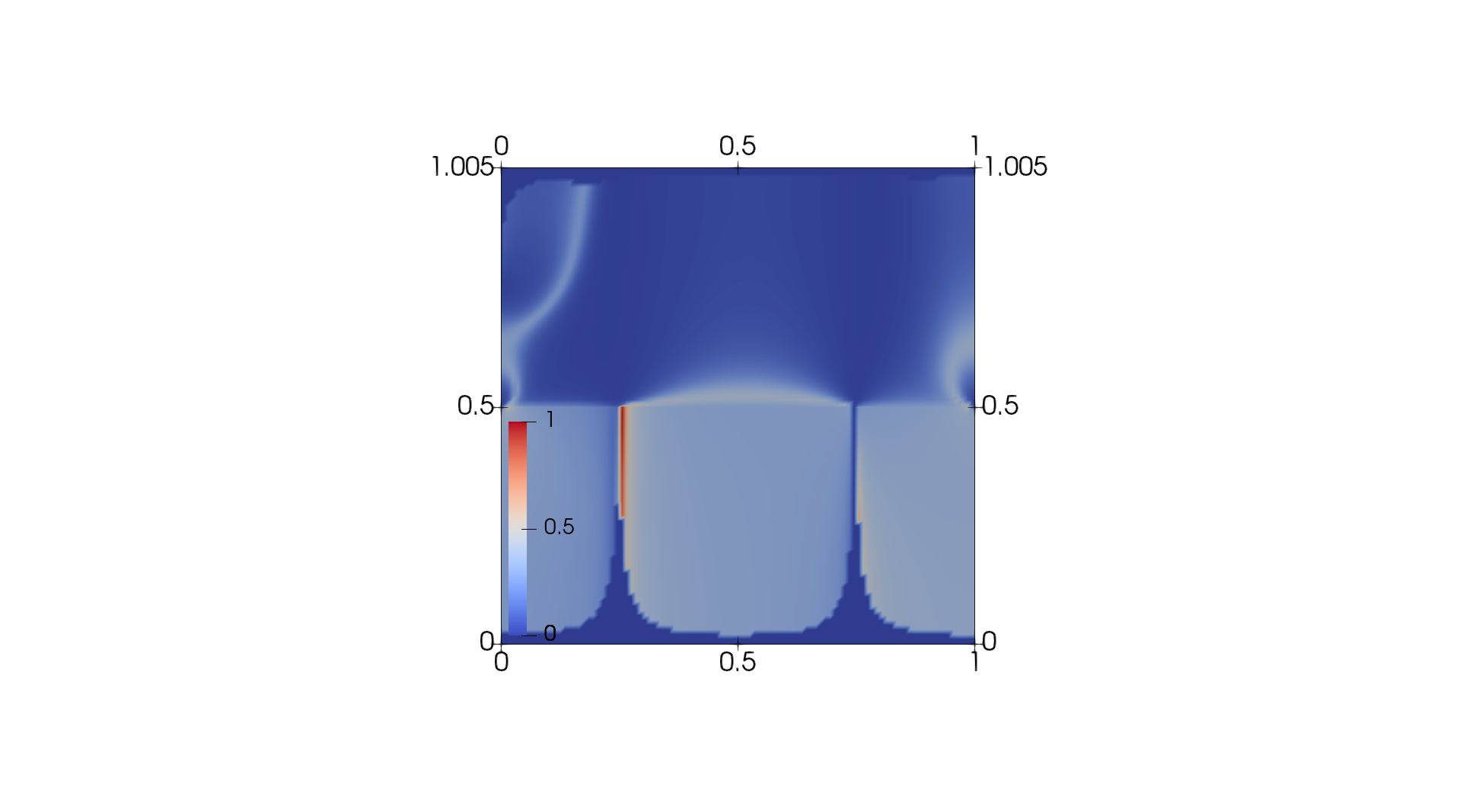}
    \caption{Velocity $v$, input parameter $I=5$}
    \end{subfigure}
   
    \caption{Case 1: Total Sobol' indices 
    of the vertical velocity $v$
    obtained from \aMR-surrogate with $\Nr=1, \No=2$,
    trained on 8192 QMC samples.}
    \label{fig:amr_sobol-v}
\end{figure}

\begin{figure}
    \centering
    \begin{subfigure}{.45\linewidth}
    \centering
    \includegraphics[width=0.65\linewidth]{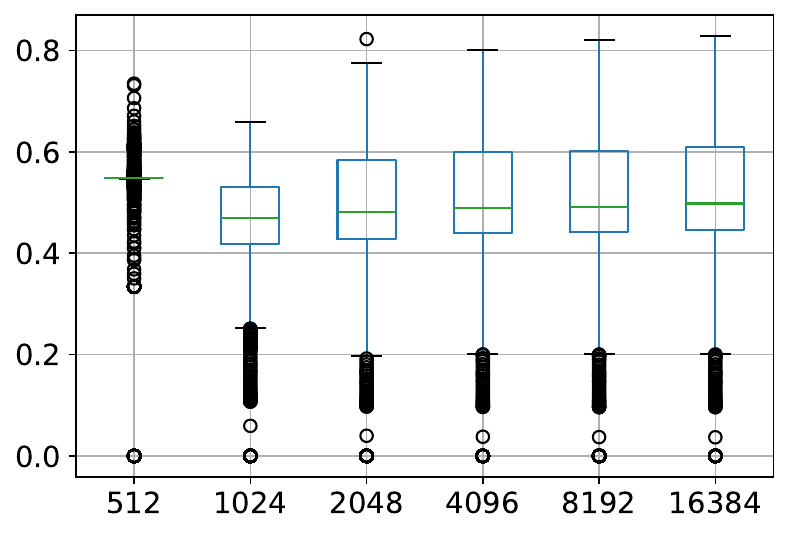}
    \caption{Input parameter $I=1$}
    \end{subfigure}
    \hfil
    \begin{subfigure}{.45\linewidth}
    \centering
    \includegraphics[width=0.65\linewidth]{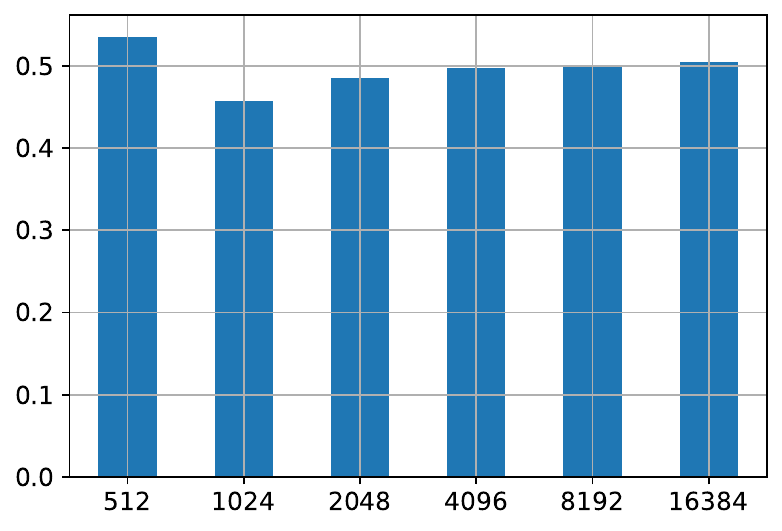}
    \caption{Input parameter $I=1$}
    \end{subfigure}

    \begin{subfigure}{.45\linewidth}
    \centering
    \includegraphics[width=0.65\linewidth]{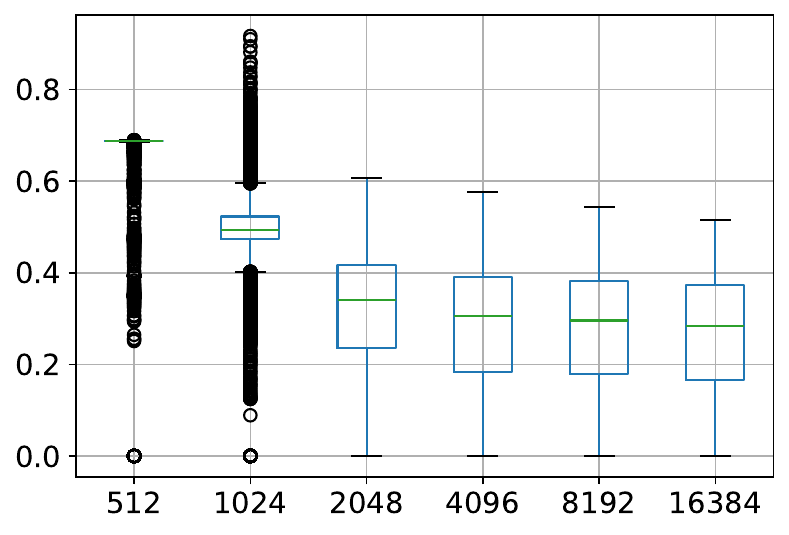}
    \caption{Input parameter $I=2$}
    \end{subfigure}
    \hfil
    \begin{subfigure}{.45\linewidth}
    \centering
    \includegraphics[width=0.65\linewidth]{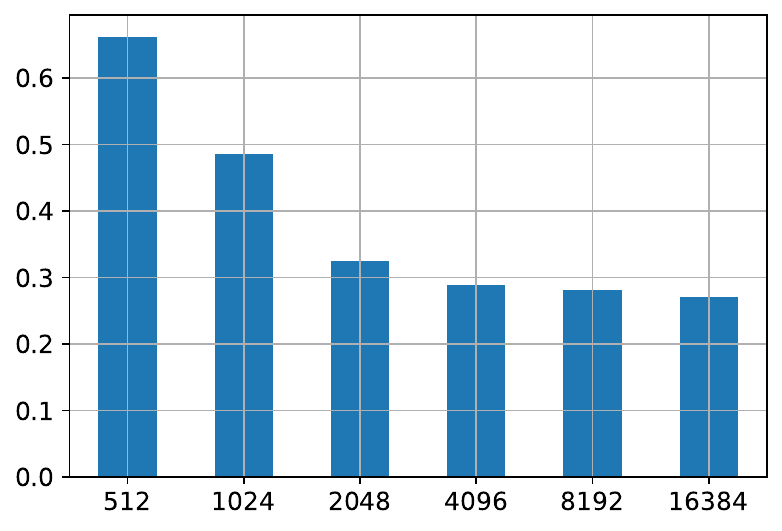}
    \caption{Input parameter $I=2$}
    \end{subfigure}

    \begin{subfigure}{.45\linewidth}
    \centering
    \includegraphics[width=0.65\linewidth]{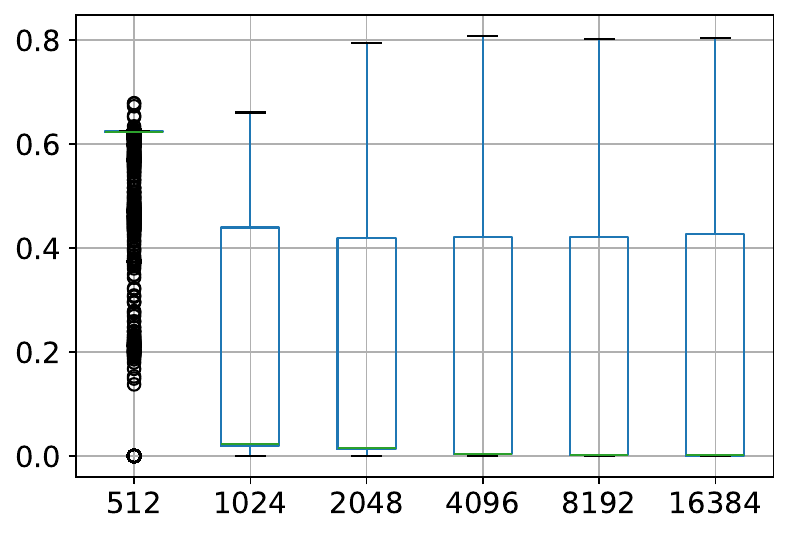}
    \caption{Input parameter $I=3$}
    \end{subfigure}
    \hfil
    \begin{subfigure}{.45\linewidth}
    \centering
    \includegraphics[width=0.65\linewidth]{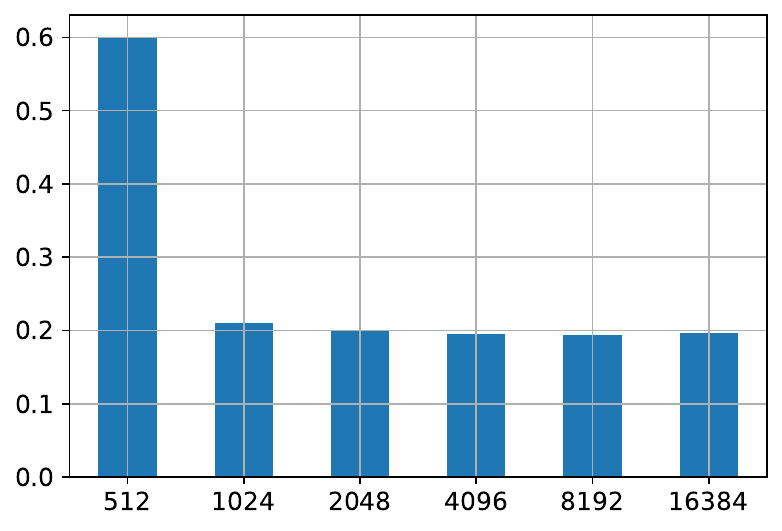}
    \caption{Input parameter $I=3$}
    \end{subfigure}

    \begin{subfigure}{.45\linewidth}
    \centering
    \includegraphics[width=0.65\linewidth]{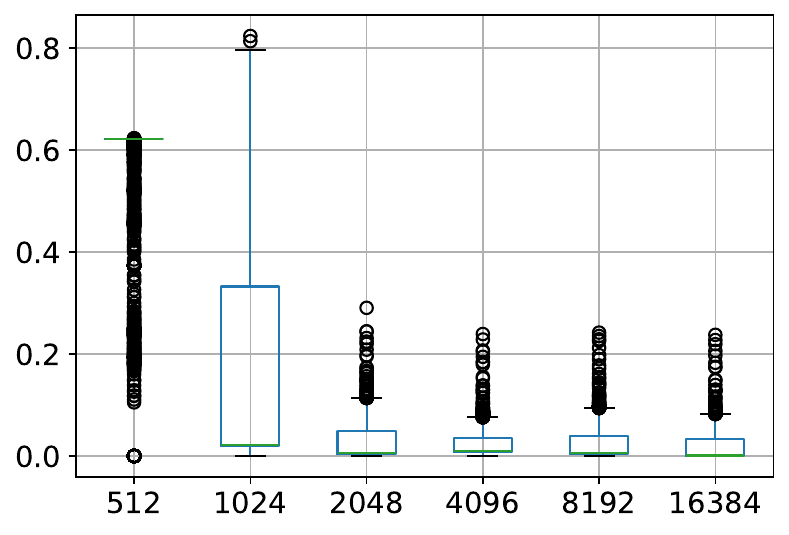}
    \caption{Input parameter $I=4$}
    \end{subfigure}
    \hfil
    \begin{subfigure}{.45\linewidth}
    \centering
    \includegraphics[width=0.65\linewidth]{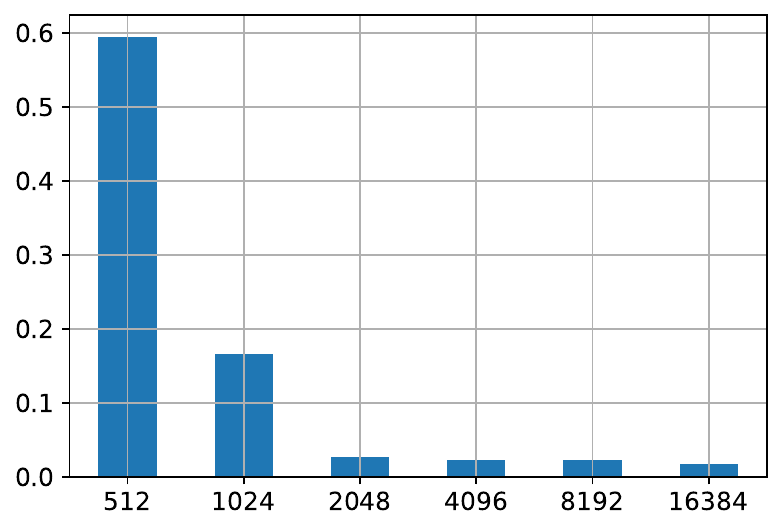}
    \caption{Input parameter $I=4$}
    \end{subfigure}

    \begin{subfigure}{.45\linewidth}
    \centering
    \includegraphics[width=0.65\linewidth]{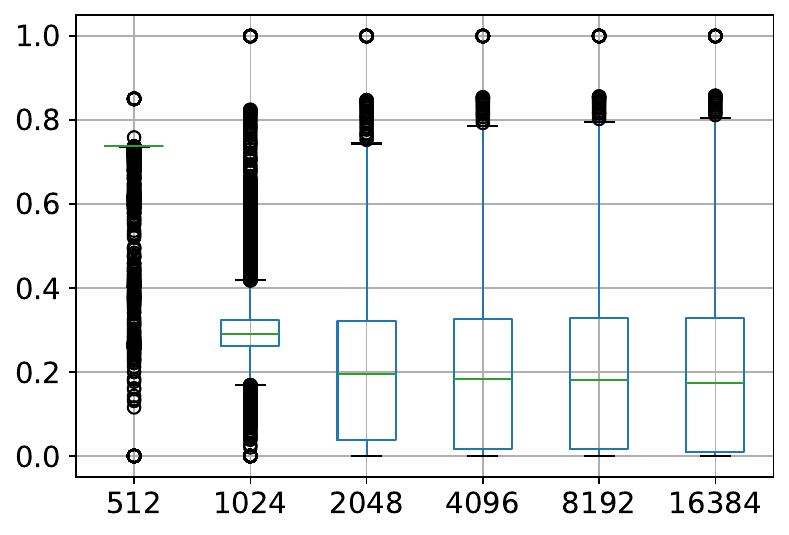}
    \caption{Input parameter $I=5$}
    \end{subfigure}
    \hfil
    \begin{subfigure}{.45\linewidth}
    \centering
    \includegraphics[width=0.65\linewidth]{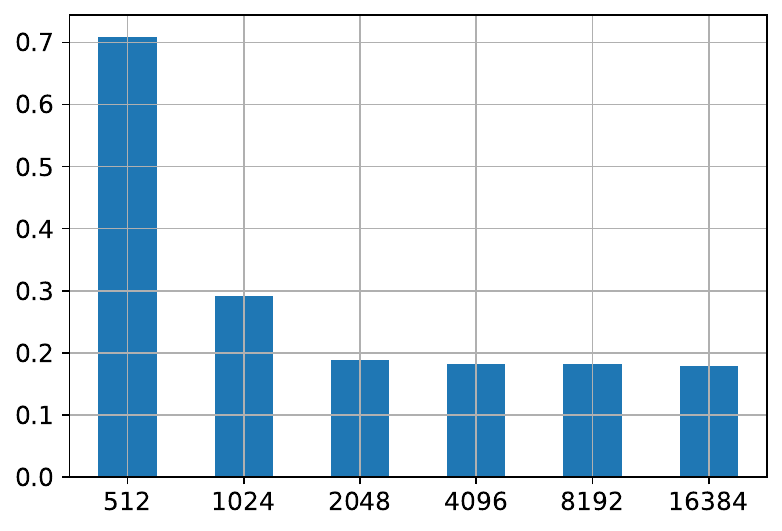}
    \caption{Input parameter $I=5$}
    \end{subfigure}
    
    \caption{Case 1: Evolution of the total sensitivity indices of the horizontal velocity $u$ obtained
    from \aMR-surrogate with $\Nr=1,\, \No=2$.}
    \label{fig:tot_evo_u}
\end{figure}

\begin{figure}
    \centering
    \begin{subfigure}{.45\linewidth}
    \centering
    \includegraphics[width=0.65\linewidth]{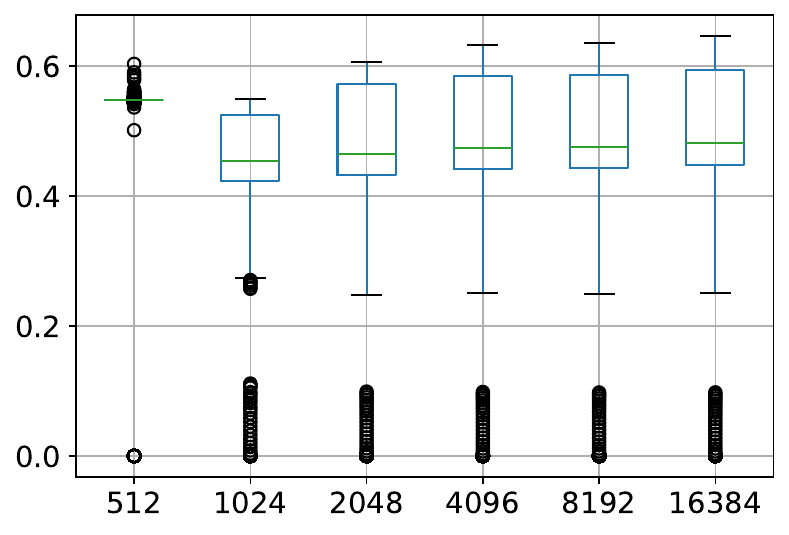}
    \caption{Input parameter $I=1$}
    \end{subfigure}
    \hfil
    \begin{subfigure}{.45\linewidth}
    \centering
    \includegraphics[width=0.65\linewidth]{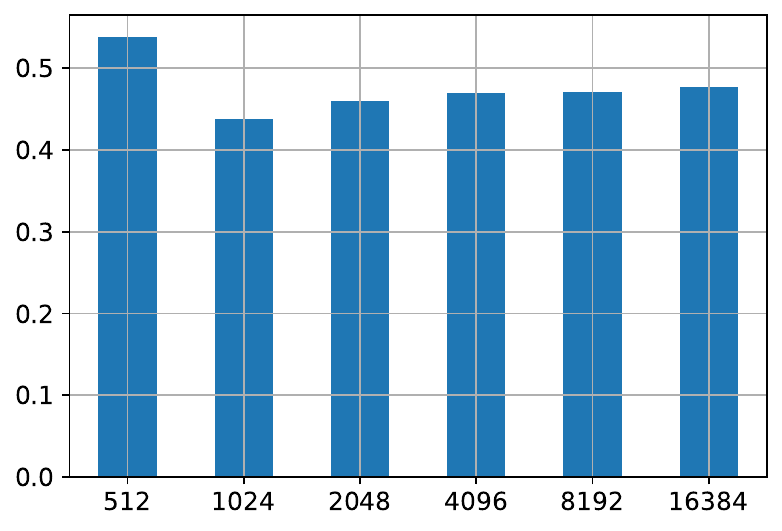}
    \caption{Input parameter $I=1$}
    \end{subfigure}

    \begin{subfigure}{.45\linewidth}
    \centering
    \includegraphics[width=0.65\linewidth]{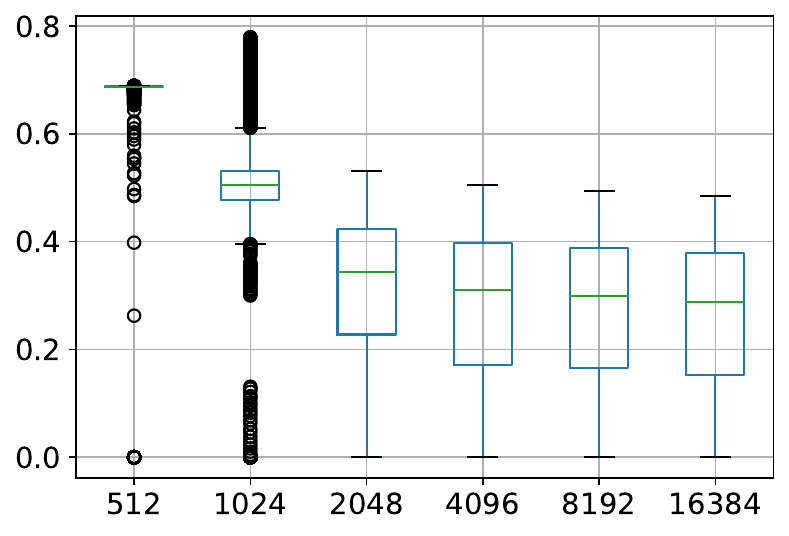}
    \caption{Input parameter $I=2$}
    \end{subfigure}
    \hfil
    \begin{subfigure}{.45\linewidth}
    \centering
    \includegraphics[width=0.65\linewidth]{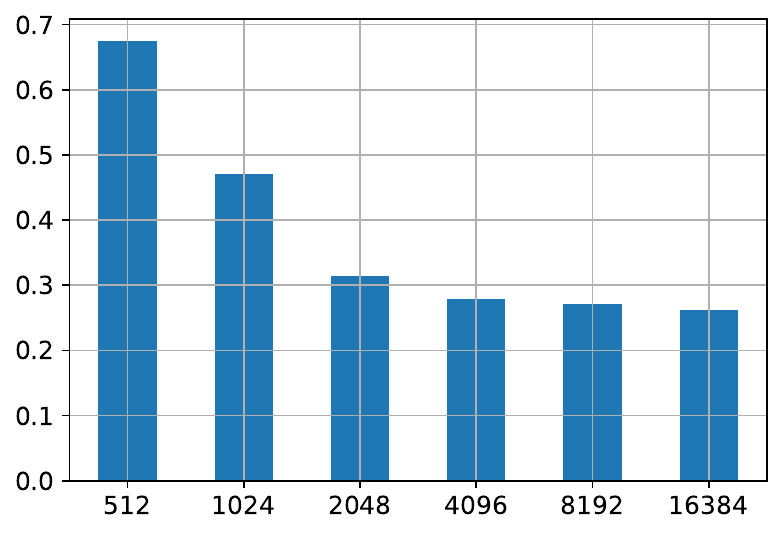}
    \caption{Input parameter $I=2$}
    \end{subfigure}

    \begin{subfigure}{.45\linewidth}
    \centering
    \includegraphics[width=0.65\linewidth]{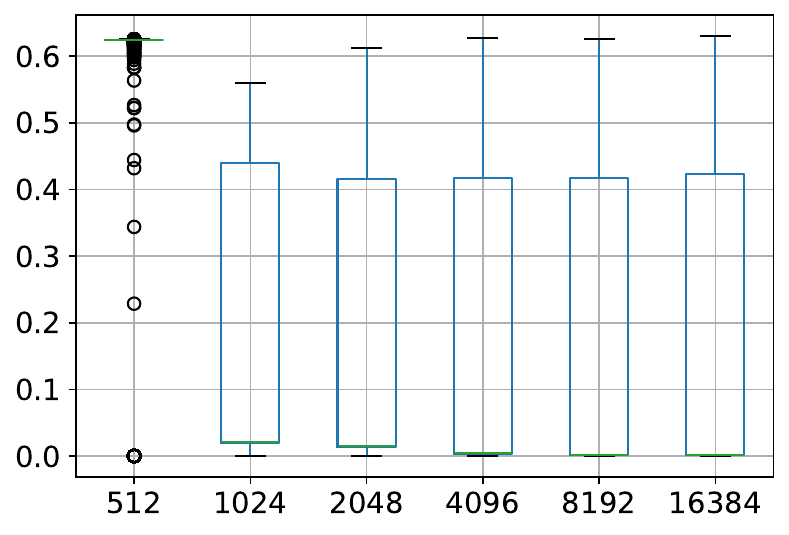}
    \caption{Input parameter $I=3$}
    \end{subfigure}
    \hfil
    \begin{subfigure}{.45\linewidth}
    \centering
    \includegraphics[width=0.65\linewidth]{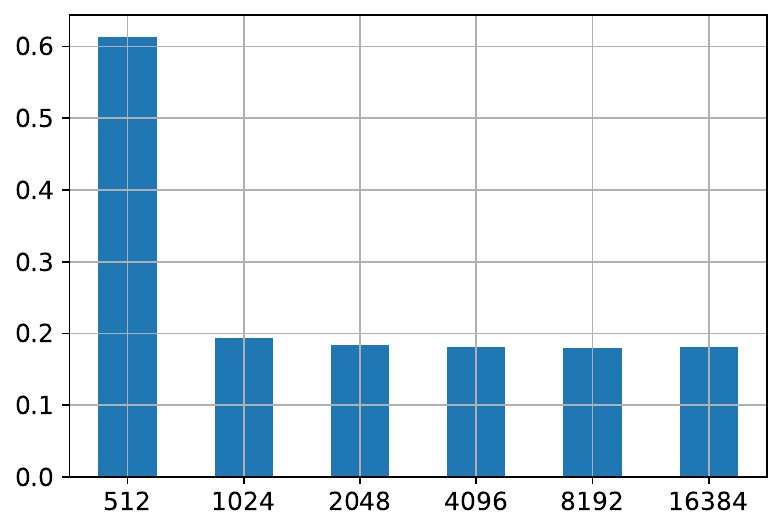}
    \caption{Input parameter $I=3$}
    \end{subfigure}

    \begin{subfigure}{.45\linewidth}
    \centering
    \includegraphics[width=0.65\linewidth]{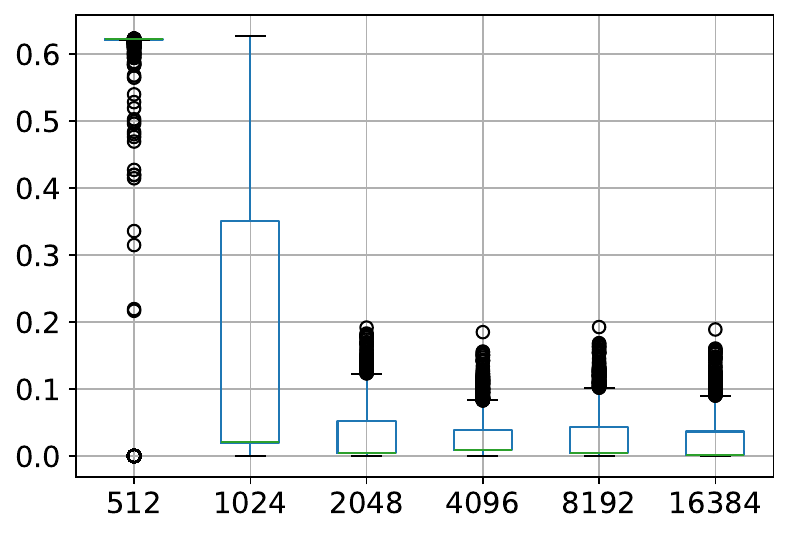}
    \caption{Input parameter $I=4$}
    \end{subfigure}
    \hfil
    \begin{subfigure}{.45\linewidth}
    \centering
    \includegraphics[width=0.65\linewidth]{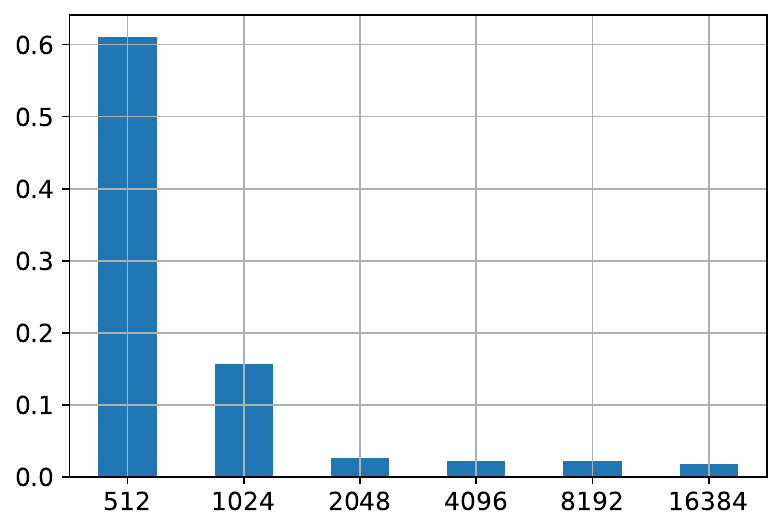}
    \caption{Input parameter $I=4$}
    \end{subfigure}

    \begin{subfigure}{.45\linewidth}
    \centering
    \includegraphics[width=0.65\linewidth]{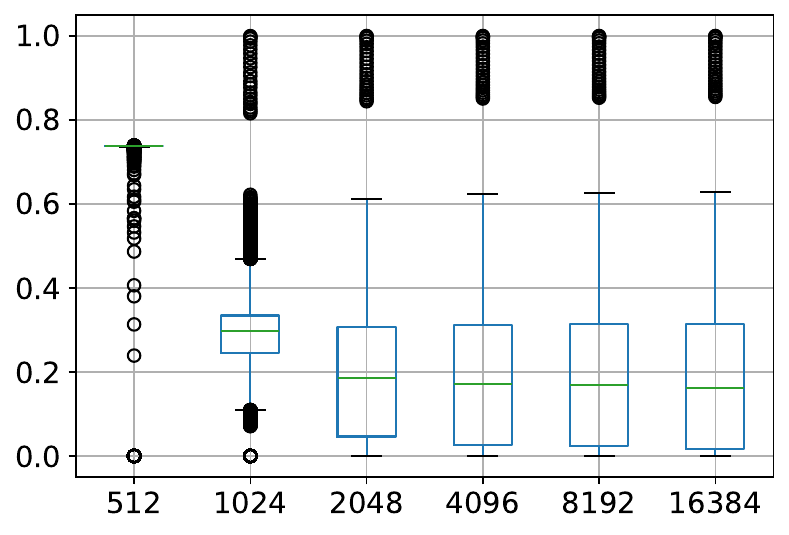}
    \caption{Input parameter $I=5$}
    \end{subfigure}
    \hfil
    \begin{subfigure}{.45\linewidth}
    \centering
    \includegraphics[width=0.65\linewidth]{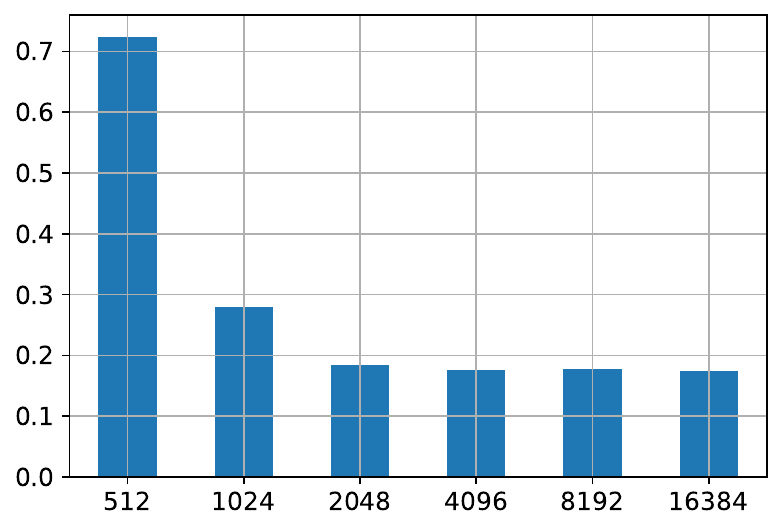}
    \caption{Input parameter $I=5$}
    \end{subfigure}
    
    \caption{Case 1: Evolution of the total sensitivity indices of the vertical velocity $v$ obtained
    from \aMR-surrogate with $\Nr=1,\, \No=2$.}
    \label{fig:tot_evo_v}
\end{figure}

In this work, we compute first-order and total-effect Sobol' indices to quantify the influence of individual parameters and their interactions on the surrogate outputs. The distributions of the total Sobol' indices, computed using \aMR{} with  $\Nr=1$, $\No=2$ on $8192$ samples, for the horizontal velocity $u$ and the vertical velocity $v$ in the filtration problem (Case 1) are shown in Fig.~\ref{fig:amr_sobol-u}~and~\ref{fig:amr_sobol-v}, respectively. Here, the total Sobol' indices of both velocity components are  presented for uncertain stress jump parameter~$\beta_\JU$ ($I=1$), permeability in the complex interface $k_\Gamma$ ($I=2$), effective viscosity $\mu_\EF$ ($I=3$), the Beavers--Joseph parameter $\alpha_\BJ$ ($I=4$) and porous-medium permeability $k_\PM$ ($I=5$). 

\begin{figure}
    \centering

    \begin{subfigure}{0.4\linewidth}
    \includegraphics[width=0.95\linewidth]{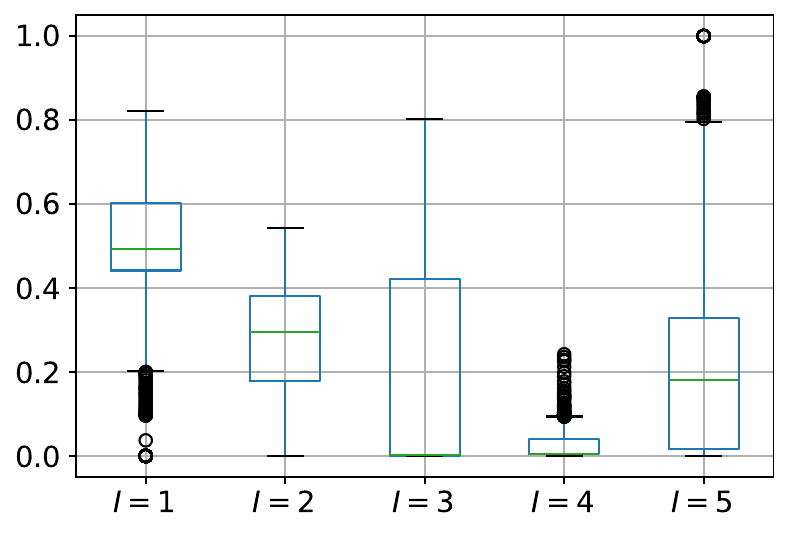}
    \caption{Box plot of total indices of for velocity $u$}
    \end{subfigure}
    \hfil 
    \begin{subfigure}{0.4\linewidth}
    \includegraphics[width=0.95\linewidth]{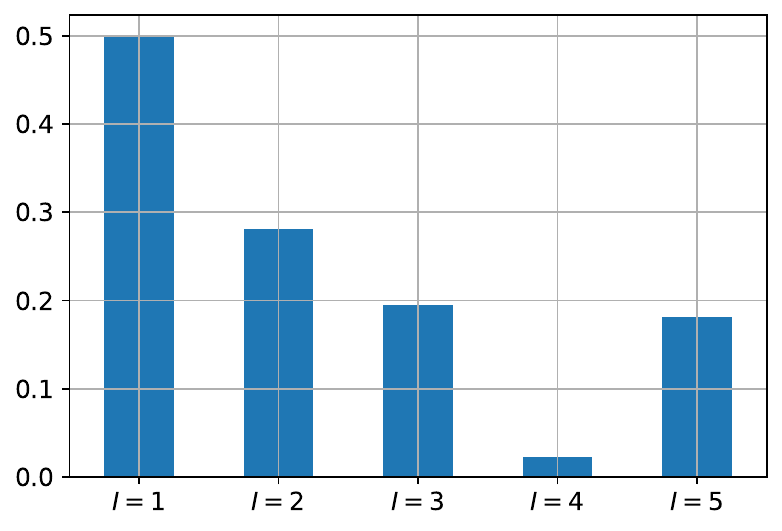}
    \caption{Averaged total indices for velocity $u$}
    \end{subfigure}\\[3mm]

    \begin{subfigure}[t]{0.4\linewidth}
    \includegraphics[width=0.95\linewidth]{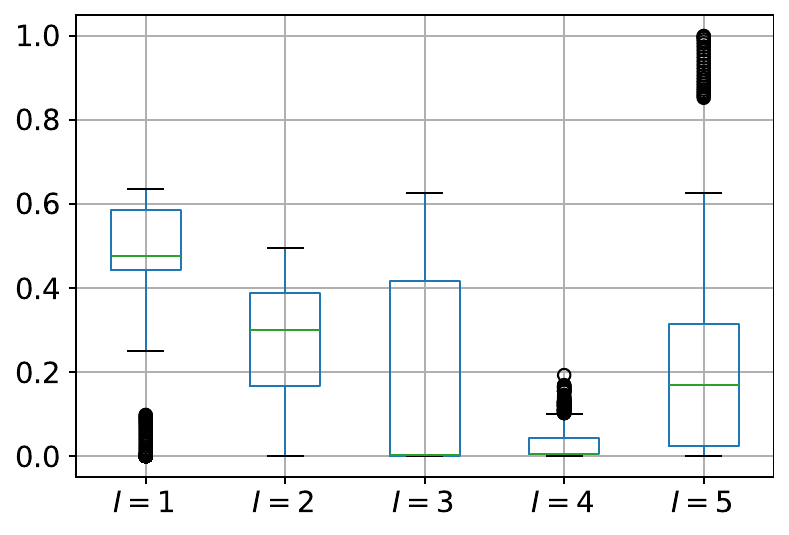}
    \caption{Box plot of total indices for velocity $v$}
    \end{subfigure}
    \hfil 
    \begin{subfigure}[t]{0.4\linewidth}
    \includegraphics[width=0.95\linewidth]{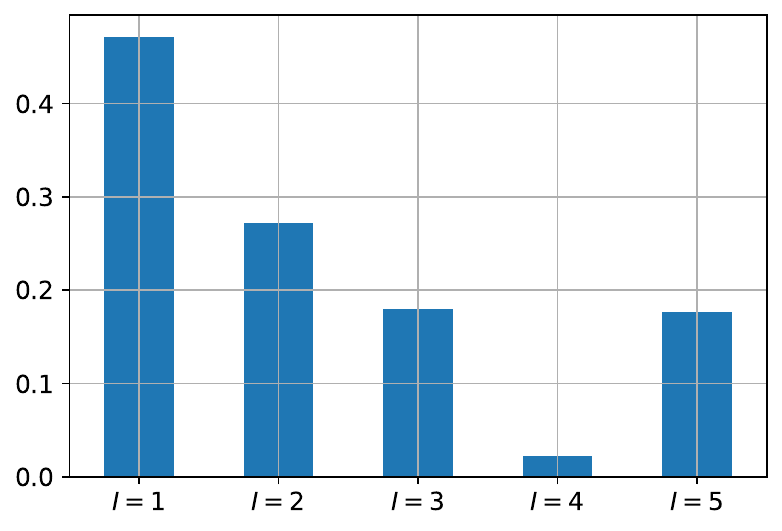}
    \caption{Averaged total indices for velocity $v$}
    \end{subfigure}
    
    \caption{Case 1: Evaluation of the total sensitivity indices 
    over the whole domain,
    obtained from \aMR-surrogate with $\Nr=1, \No=2$,
    trained on 8192 QMC samples.
    }
    \label{fig:sobol-total_av}
\end{figure}
In Fig.~\ref{fig:amr_sobol-u} and \ref{fig:amr_sobol-v}, we observe that both velocity components exhibit highly nonlinear behavior in response to variations in all five uncertain parameters. In particular, sharp fronts can be observed in the middle of the free-flow subdomain. This nonlinearity arises from the boundary conditions of the filtration problem. The stress jump parameter $\beta_\JU$ and the permeability in the complex interface $k_\Gamma$ significantly influence both velocity components in the entire flow domain (Fig.~\ref{fig:amr_sobol-u}(A,B), Fig.~\ref{fig:amr_sobol-v}(A,B)). These two parameters appear in the stress jump terms between the free-flow domain and the complex interface $\Gamma$ (see~\eqref{eqn:AveragedBrinkmanMomentumNormal}--\eqref{eqn:TransmissionGammaFFtangential}). In the filtration problem, the flow is driven from the porous-medium region into the free-flow domain, therefore the velocities in these two flow regions are comparable (Fig.~\ref{fig:testcaseFP}, right). The coupled model is sensitive to variations in the stress jump parameter and the permeability in the complex interface in both subdomains. In contrast, as we will see in Case~2 later, different behavior arises when the velocity in the free-flow region is much higher than in the porous medium.

Figures~\ref{fig:amr_sobol-u}(C) and \ref{fig:amr_sobol-v}(C) indicate that the effective viscosity $\mu_\EF$ mainly affects the flow in the free-flow domain. Although $\mu_\EF$ appears in the higher-order correction terms of the averaged Brinkman equations \eqref{eqn:AveragedBrinkmanMomentumNormal}, \eqref{eqn:AveragedBrinkmanMomentumTangential} and in the transmission conditions~\eqref{eqn:TransmissionGammaFFNormal}--\eqref{eqn:TransmissionGammaPM}, its impact on the porous-medium region remains small. The Beavers--Joseph parameter $\alpha_\BJ$ has the smallest effect among the parameters considered (Fig.~\ref{fig:amr_sobol-u}(D), \ref{fig:amr_sobol-v}(D)), as it mainly reflects the interface roughness at the bottom of the transition region.  Compared to the porous-medium region, the flow in the free-flow domain is slightly more sensitive to variations in the parameter $\alpha_\BJ$. Figures~\ref{fig:amr_sobol-u}(E) and \ref{fig:amr_sobol-v}(E) show that both velocity components in the porous medium are very sensitive to changes in the porous-medium permeability~$k_\PM$, as it appears in Darcy’s law~\eqref{eqn:Darcy2} and directly controls the velocity magnitude in the porous-medium region. These findings are consistent with the results presented in~\cite{kroker2023global}, where the GSA was performed for the classical Stokes--Darcy model.

We note that both Sobol' indices and total sensitivity indices are normalized by the variance (Fig.~\ref{fig:amr_mean_var}(E,F)). Consequently, numerical errors may be amplified in regions where the variance is close to zero, as illustrated by the sharp fronts in the porous-medium region in Fig.~\ref{fig:amr_sobol-u}(E) and \ref{fig:amr_sobol-v}(E). The results in such regions are therefore not necessarily meaningful from a modeling perspective.

To obtain a compact representation of the parameter influence, we employ space-averaged Sobol’ sensitivity indices, which provide a global measure of parameter effects in the flow domain. 
Figures~\ref{fig:tot_evo_u} and \ref{fig:tot_evo_v} display box plots and mean bar plots of the Sobol’ indices for five uncertain parameters across different sample sizes ranging from $512$ to $16,384$.
Convergence of the Sobol’ sensitivity indices with respect to the number of samples is evident for both velocity components (see box plots and histograms in Fig.~\ref{fig:tot_evo_u} and~\ref{fig:tot_evo_v}). From $2048$ samples onward, both the distributions and the averaged values exhibit stable behavior. 
We therefore select $8192$ samples as a representative case and present the corresponding box plots of the Sobol’ indices together with histograms of their averaged values in Fig.~\ref{fig:sobol-total_av}. The results reveal that, within the framework of the newly derived hybrid-dimensional Stokes–Brinkman–Darcy model, the flow model is most sensitive to the stress jump parameter~$\beta_\JU$ (Fig.~\ref{fig:sobol-total_av}(B,D)).  The permeability of the complex interface~$k_\Gamma$ and the porous-medium permeability~$k_\PM$ have also significant influence. The box plot of the effective viscosity~$\mu_\EF$ shows relatively low median values (Fig.~\ref{fig:sobol-total_av}(A,C)). Nevertheless, the averaged Sobol’ index indicates significant influence of $\mu_\EF$, which can be explained by its primary impact in the free-flow domain. The Beavers-–Joseph parameter $\alpha_\BJ$ shows only minor impact in this test case, although its box plot variability indicates occasional strong contributions in specific realizations (outliers).

\subsection{Case~2: Postprocessing and evaluation}\label{sec:postprocessingcase2}
We next proceed with the second test case introduced in Section~\ref{sec:testcase2SF}, which describes the splitting flow problem. In this setting, the flow enters and exits through the boundary of the free-flow domain, resulting in significantly higher velocities in the free-flow region compared to the porous medium. The permeability of the porous medium controls the amount of fluid that can pass through the complex interface. Based on the experience gained from Case~1, we omit the evaluations of the \aMR{} surrogate models for $\Nr=1,\,\No=4$, both with and without hyperbolic truncation-based sparsity.

\subsubsection{Mean and variance}
\begin{figure}
    \centering
    \begin{subfigure}{0.4\linewidth}
            \includegraphics[trim=560pt 185pt 560pt 175pt, clip, width=1.0\linewidth]{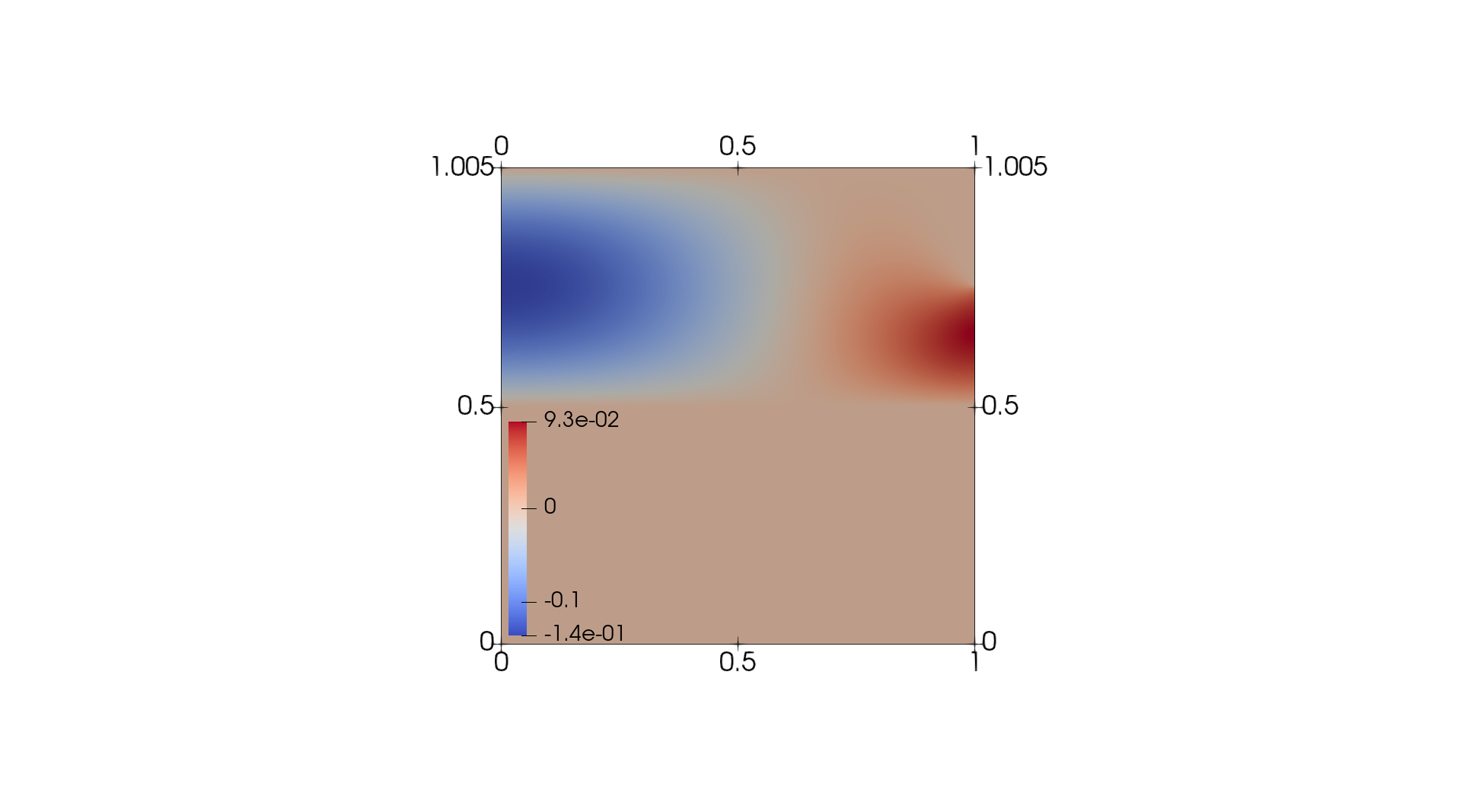}
    \caption{$\E[u]$}
    \end{subfigure}
    \hfil
    \begin{subfigure}{0.4\linewidth}
            \includegraphics[trim=560pt 185pt 560pt 175pt, clip, width=1.0\linewidth]{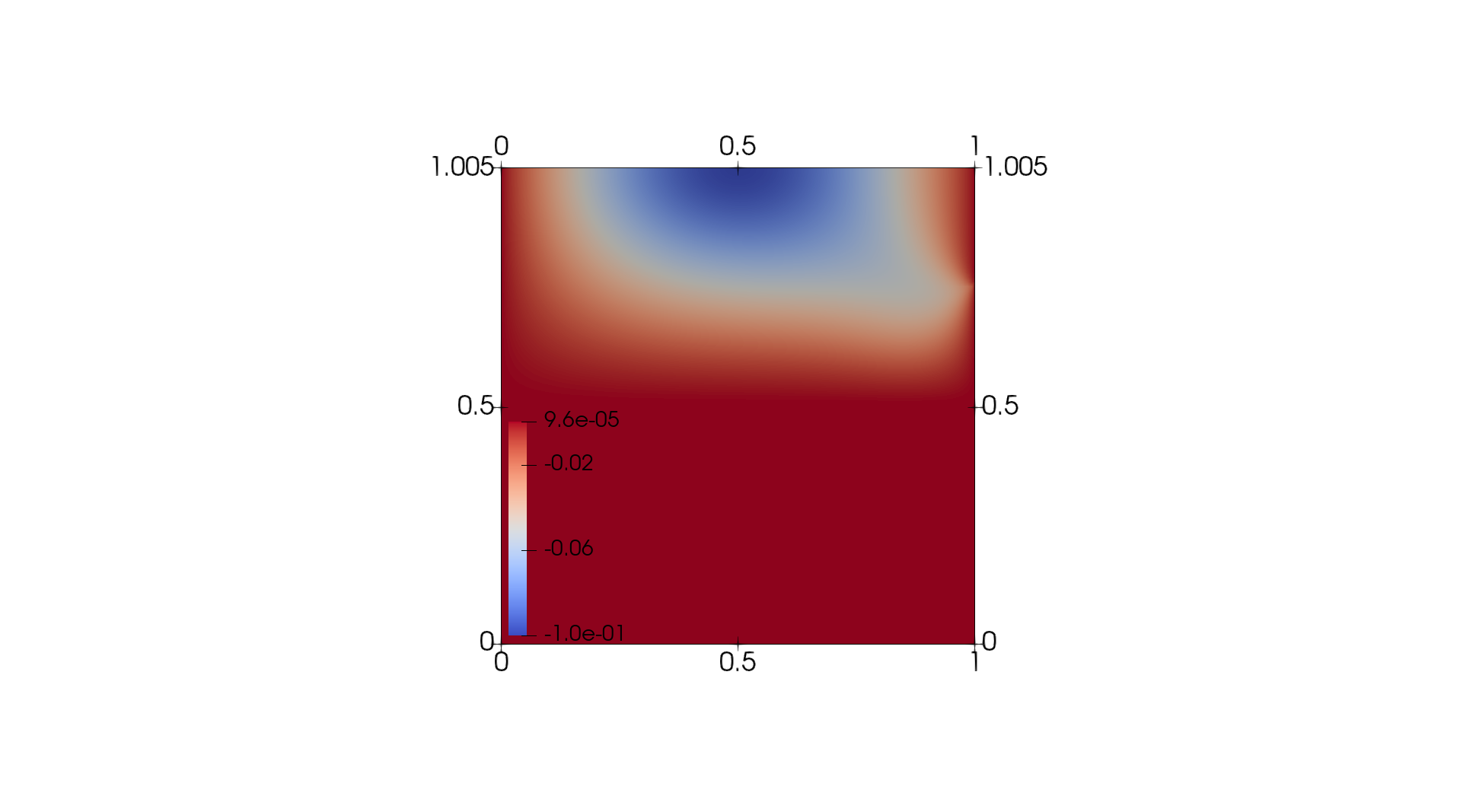}
    \caption{$\E[v]$}
    \end{subfigure}
    
    \begin{subfigure}{0.4\linewidth}
            \includegraphics[trim=560pt 185pt 560pt 175pt, clip, width=1.0\linewidth]{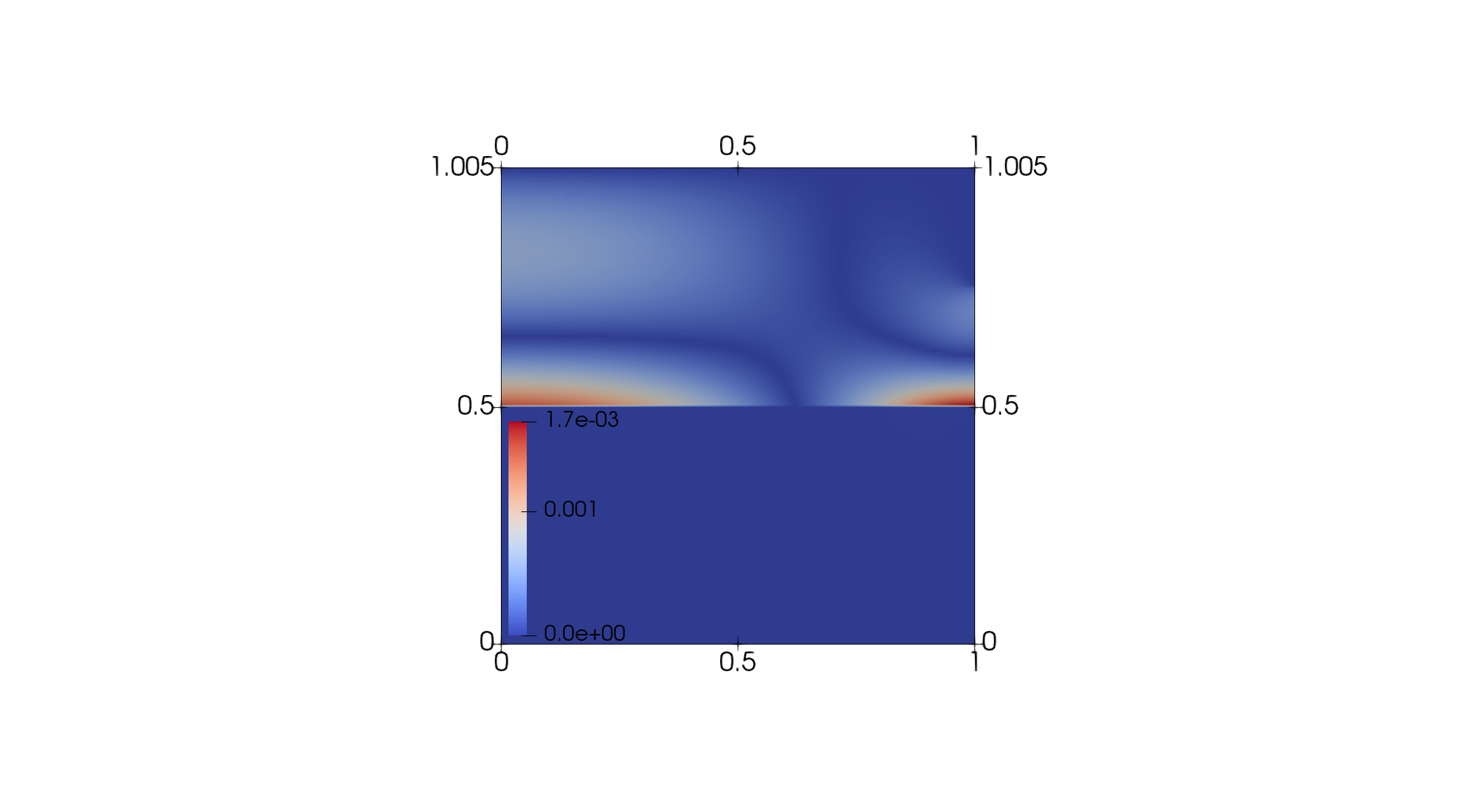}
    \caption{$\std[u]$}
    \end{subfigure}
    \hfil
    \begin{subfigure}{0.4\linewidth}
            \includegraphics[trim=560pt 185pt 560pt 175pt, clip, width=1.0\linewidth]{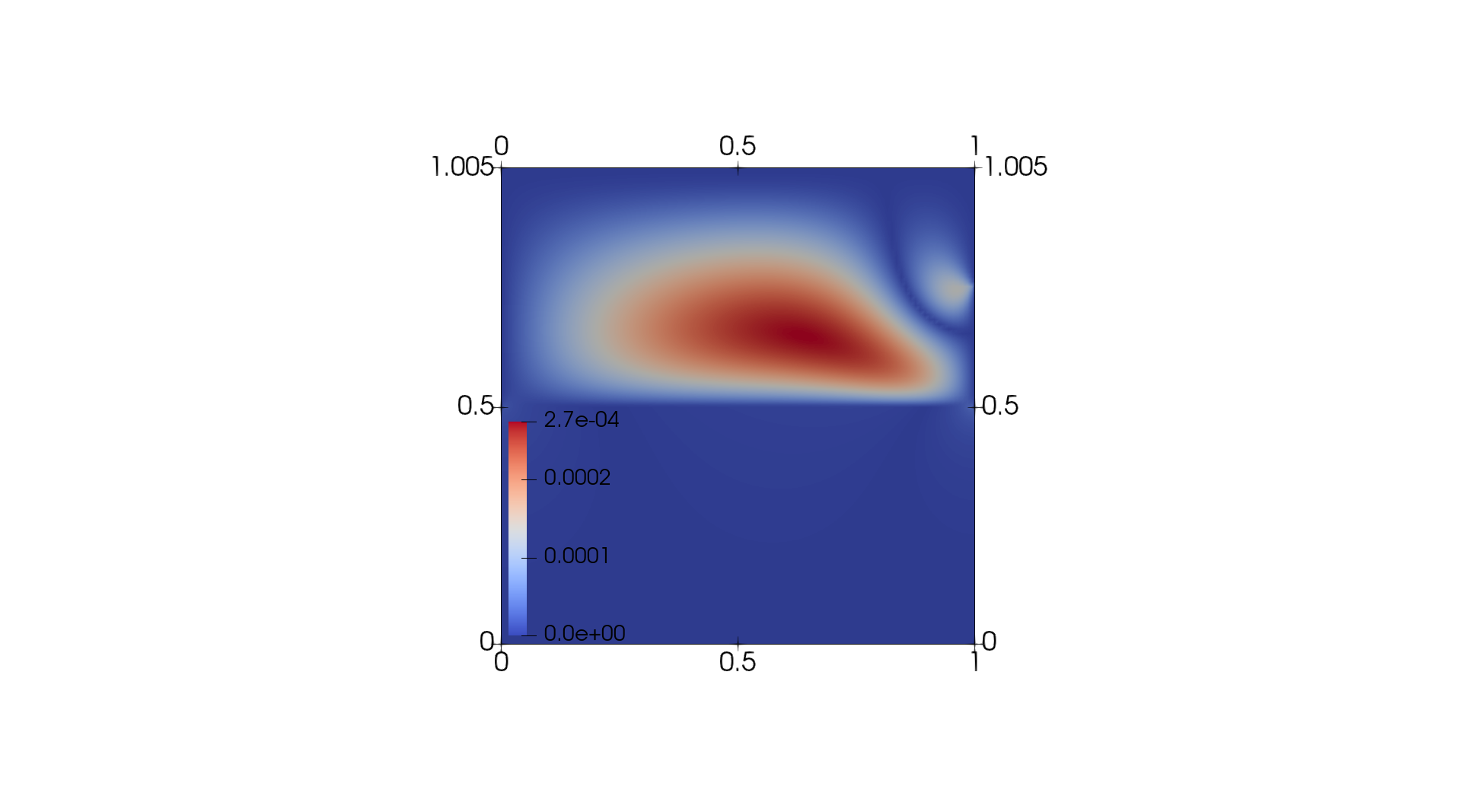}
    \caption{$\std[v]$}
    \end{subfigure}

    \begin{subfigure}{0.4\linewidth}
            \includegraphics[trim=560pt 185pt 560pt 175pt, clip, width=1.0\linewidth]{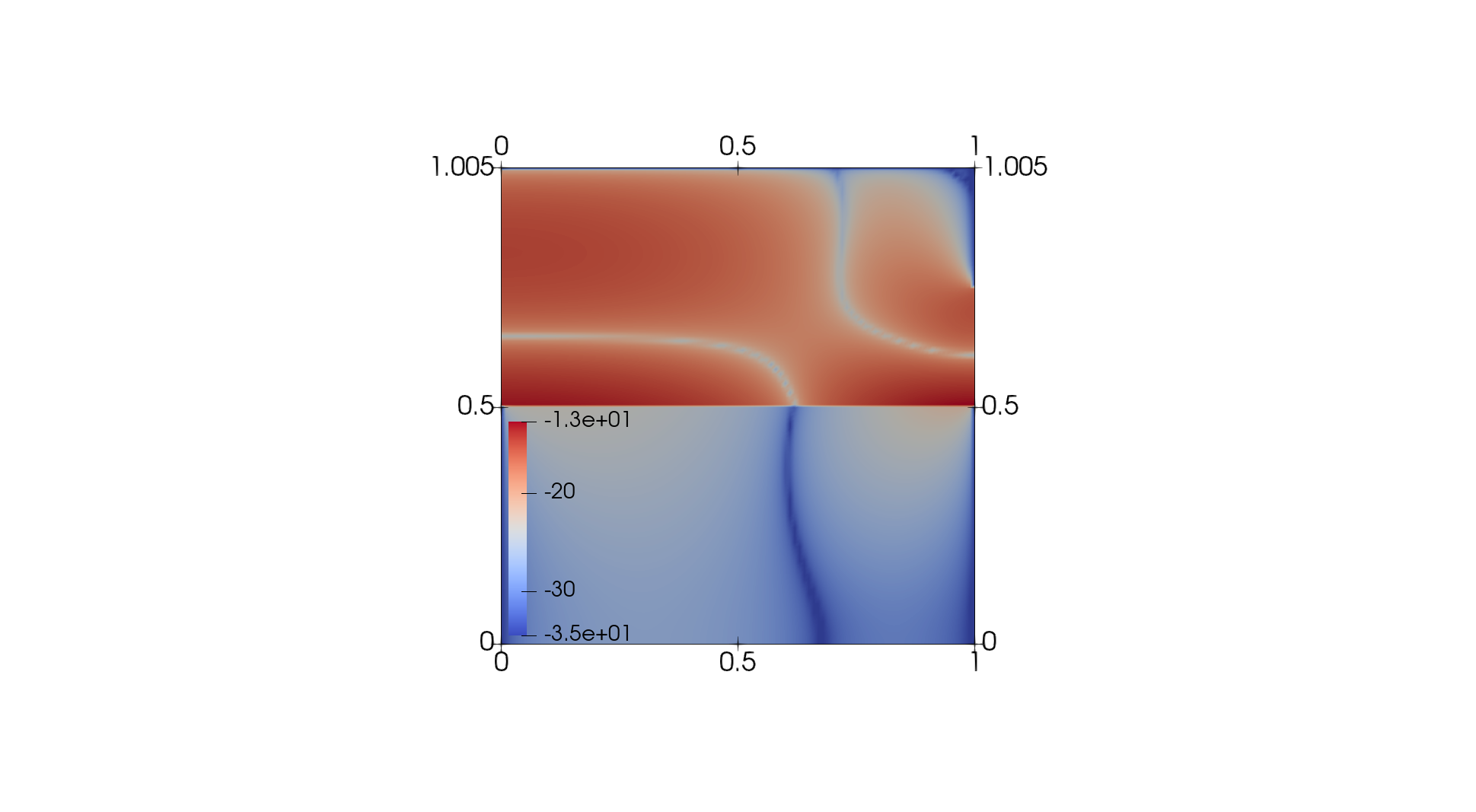}
    \caption{$\log\var[u]$}
    \end{subfigure}
    \hfil
    \begin{subfigure}{0.4\linewidth}
            \includegraphics[trim=560pt 185pt 560pt 175pt, clip, width=1.0\linewidth]{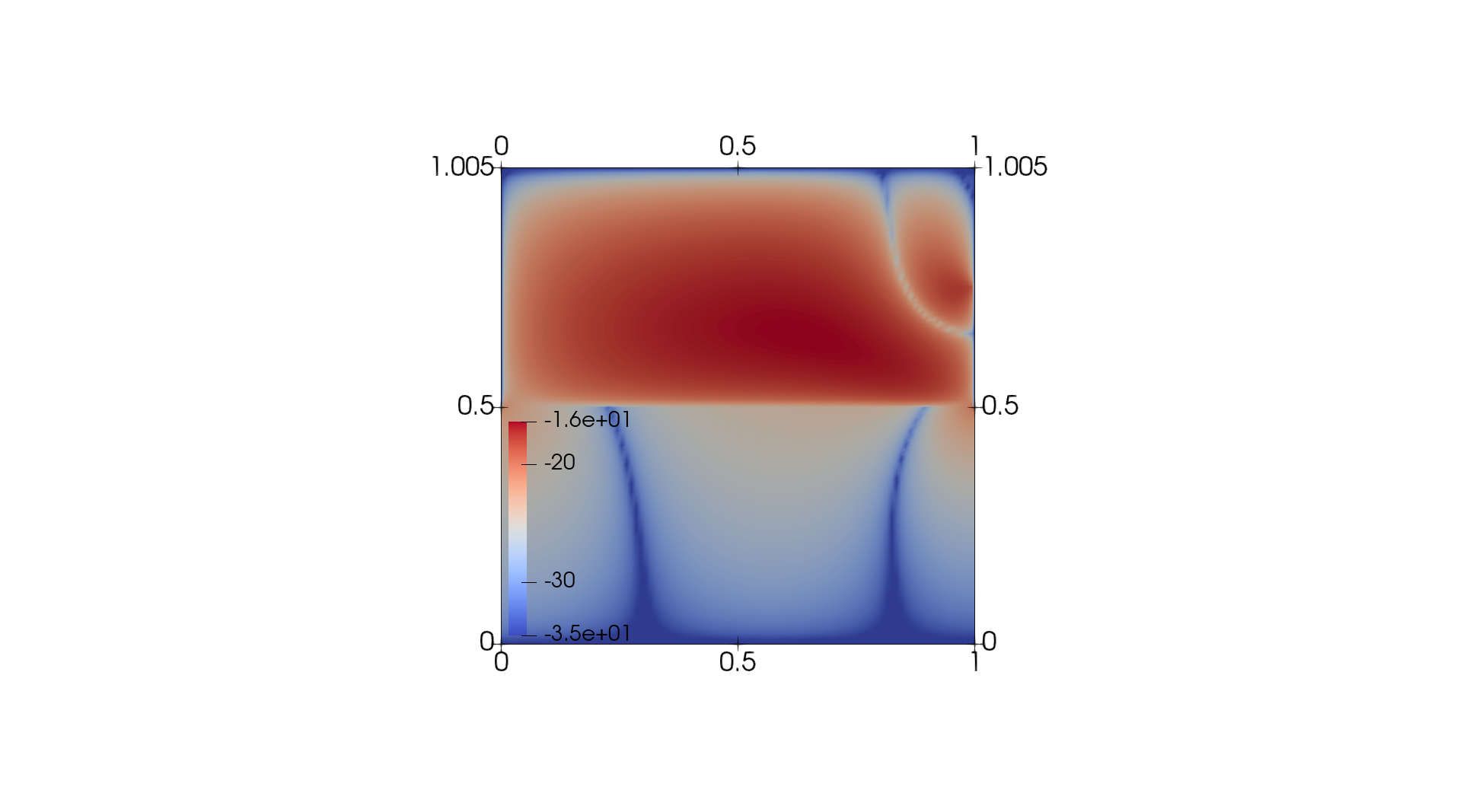}
    \caption{$\log\var[v]$}
    \end{subfigure}
    
    \caption{Case 2: Mean, standard deviation and log-variance of the horizontal
    and vertical velocities $(u, v)$
    obtained from \aMR-surrogate with $\Nr=1, \No=2$,
    trained on 8192 QMC samples.}
    \label{fig:amr_mean_var_ii}
\end{figure}

We visualize the mean, standard deviation, and log-variance obtained from the \aMR{} expansion with polynomial degree $\No=2$ and refinement level $\Nr=1$ in Fig.~\ref{fig:amr_mean_var_ii}. The comparison with MC results using $50{,}000$ samples is omitted here, since the \aMR{} expansion reproduces the MC statistics with high fidelity. The mean and variance in Fig.~\ref{fig:amr_mean_var_ii}(A–D) further indicate that the velocity in the porous-medium region is relatively small, making the variance difficult to detect. The log-variance highlights the regions with comparatively larger relative fluctuations, making the small but non-negligible variability in the porous-medium region more visible (Fig.~\ref{fig:amr_mean_var_ii}(E, F)).

\subsubsection{Accuracy of the approximation}
\begin{figure}
    \centering
    \begin{subfigure}{0.45\linewidth}
    \includegraphics[width=0.95\linewidth]{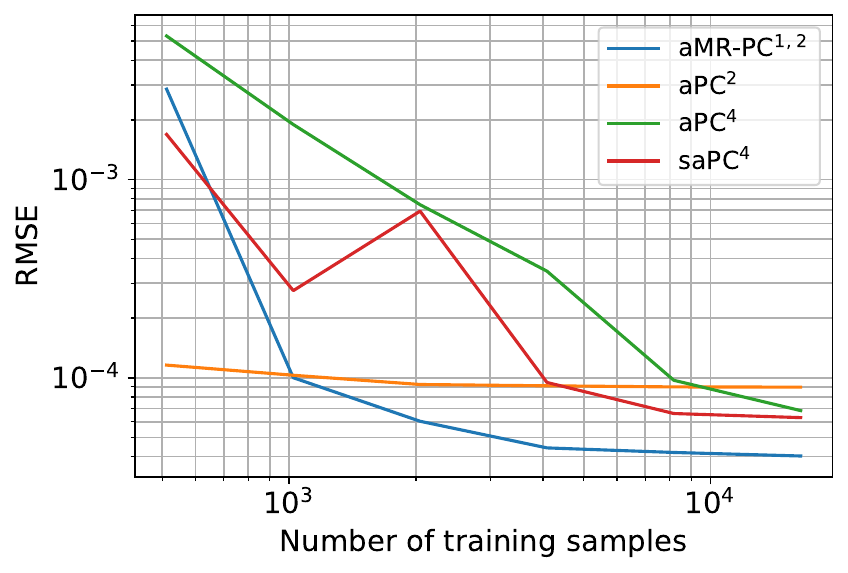}  
    \caption{RMSE of $u$}
    \end{subfigure}
    \hfil 
    \begin{subfigure}{0.45\linewidth}
    \includegraphics[width=0.95\linewidth]{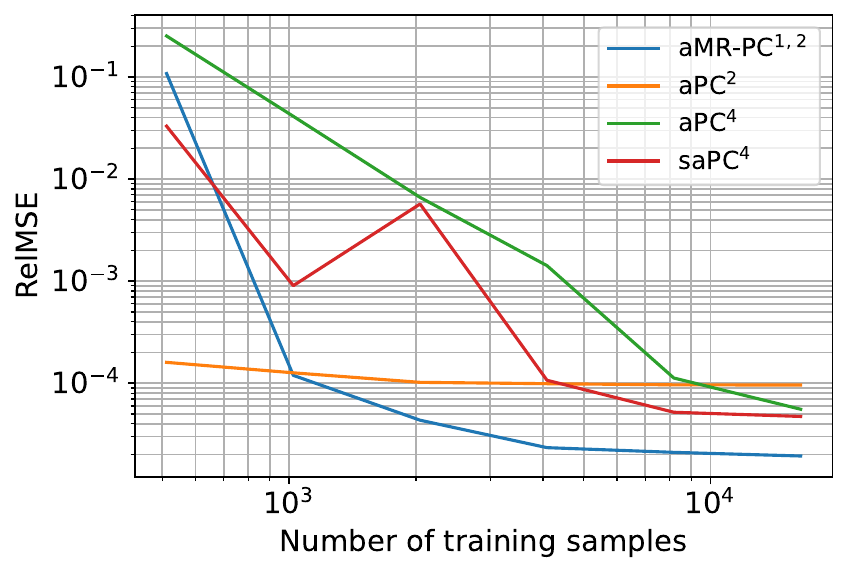}  
    \caption{Relative MSE of $u$}
    \end{subfigure}

    \begin{subfigure}{0.45\linewidth}
    \includegraphics[width=0.95\linewidth]{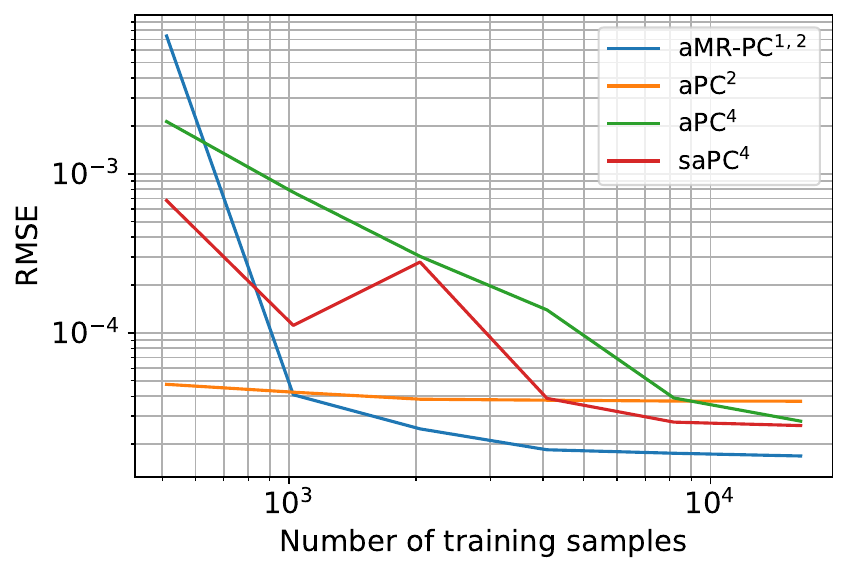}  
    \caption{RMSE of $v$}
    \end{subfigure}
    \hfil 
    \begin{subfigure}{0.45\linewidth}
    \includegraphics[width=0.95\linewidth]{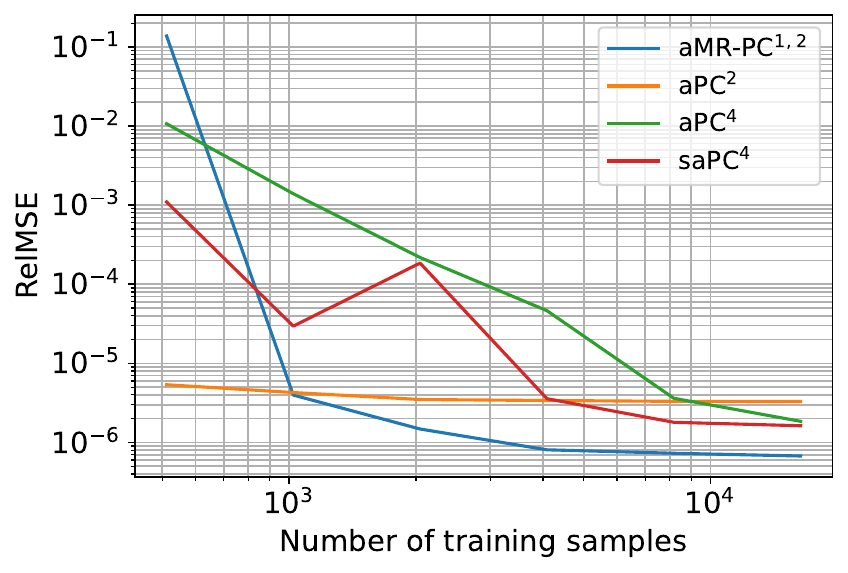}  
    \caption{Relative MSE of $v$}
    \end{subfigure}
    
    \caption{Case 2: RMSE and relative MSE error plots for $512,\ldots,16384$
    training samples.}
    \label{fig:error_ii}
\end{figure}

As in the previous test case, we start with the assessment of the accuracy in terms of
RMSE and relative MSE shown in Fig.~\ref{fig:error_ii}.
Also in this case, errors are computed over the spatial domain using 10,000 randomly selected test parameter sets not used for training.
Here, aPC with $\No=2$ provides the best accuracy using 512 training samples.
However, the improvement for increasing number of training samples is only minimal.
We observe that \aMR{} with $\Nr=1,\, \No=2$ outperforms all other surrogate models
 in terms of accuracy if more than 1024 samples are used for training.

Evaluation of the $\Lp2$-error of the approximation of the mean
using Eq.~\eqref{eq:amr-mean} in comparison with the MC result computed
with 50,000 samples is provided in Tab.~\ref{tab:error-means_ii}.
It demonstrates a similar final precision for all surrogates.
Again, aPC with $\No=2$ shows a better accuracy for lower numbers of
training samples.
Evaluating the $\Lp2$-error in approximation of the standard deviation using Eq.~\eqref{eq:amr-var} against the MC reference solution computed using 50,000 samples, as presented in Tab.~\ref{tab:error-std_ii}, reveals that aPC with $\No=2$ offers optimal accuracy with fewer training samples, while \aMR{} with $\Nr=1,\,\No=2$ provides the best precision when using a larger number of training samples.

\begin{table}
\caption{Case~2: $\Lp2$-error of the approximation of the mean $(\E)$.}
\label{tab:error-means_ii}
\centering

\begin{subtable}{0.95\linewidth}
\caption{Approximation of the mean of the horizontal velocity component $\E[u]$.}
\centering
\begin{tabular}{|r|c|c|c|c|}
\hline
Samples & $N_o=2$ & $N_o=4$ & $N_o=4 / q=0.75$ & $N_r=1 / N_o=2$ \\
\hline
512 & 5.10e-06 & 1.53e-04 & 5.21e-05 & 4.01e-03 \\
1024 & 2.12e-06 & 5.35e-05 & 8.60e-06 & 3.70e-06 \\
2048 & 4.95e-07 & 1.94e-05 & 2.07e-05 & 1.84e-06 \\
4096 & 2.69e-07 & 9.00e-06 & 1.34e-06 & 9.98e-07 \\
8192 & 4.56e-07 & 2.05e-06 & 7.32e-07 & 3.56e-07 \\
16384 & 4.41e-07 & 7.86e-07 & 4.56e-08 & 1.72e-07 \\
\hline
\end{tabular}

\end{subtable}\\[3mm]

\begin{subtable}{0.95\linewidth}
\caption{Approximation of the mean of the vertical velocity component $\E[v]$.}
\centering
\begin{tabular}{|r|c|c|c|c|}
\hline
Samples & $N_o=2$ & $N_o=4$ & $N_o=4 / q=0.75$ & $N_r=1 / N_o=2$ \\
\hline
512 & 1.68e-06 & 5.04e-05 & 1.72e-05 & 3.11e-03 \\
1024 & 6.92e-07 & 1.77e-05 & 2.82e-06 & 1.21e-06 \\
2048 & 1.59e-07 & 6.42e-06 & 6.83e-06 & 6.04e-07 \\
4096 & 9.52e-08 & 2.98e-06 & 4.69e-07 & 3.26e-07 \\
8192 & 1.56e-07 & 6.67e-07 & 2.44e-07 & 1.15e-07 \\
16384 & 1.49e-07 & 2.55e-07 & 1.74e-08 & 6.03e-08 \\
\hline
\end{tabular}

\end{subtable}

\end{table}
\begin{table}
\caption{Case~2: $\Lp2$-error of the approximation of 
the standard deviation $\std$.}
\label{tab:error-std_ii}
\centering

\begin{subtable}{0.95\linewidth}
\caption{Approximation of the standard deviation of the horizontal velocity component $\std[u]$.}
\centering
\begin{tabular}{|r|c|c|c|c|}
\hline
Samples & $N_o=2$ & $N_o=4$ & $N_o=4 / q=0.75$ & $N_r=1 / N_o=2$ \\
\hline
512 & 8.58e-06 & 1.50e-02 & 5.55e-03 & 1.21e-02 \\
1024 & 1.32e-05 & 5.59e-03 & 9.61e-04 & 6.73e-05 \\
2048 & 2.45e-05 & 2.35e-03 & 2.35e-03 & 4.93e-06 \\
4096 & 2.34e-05 & 1.12e-03 & 1.70e-05 & 3.06e-06 \\
8192 & 2.39e-05 & 3.09e-04 & 8.77e-05 & 1.66e-06 \\
16384 & 2.39e-05 & 1.17e-04 & 2.55e-05 & 2.21e-06 \\
\hline
\end{tabular}
\end{subtable}\\[3mm]

\begin{subtable}{0.95\linewidth}
\caption{Approximation of the standard deviation of the 
vertical velocity component $\std[v]$.}
\centering
\begin{tabular}{|r|c|c|c|c|}
\hline
Samples & $N_o=2$ & $N_o=4$ & $N_o=4 / q=0.75$ & $N_r=1 / N_o=2$ \\
\hline
512 & 2.90e-06 & 4.95e-03 & 1.83e-03 & 9.45e-03 \\
1024 & 4.38e-06 & 1.85e-03 & 3.16e-04 & 2.26e-05 \\
2048 & 8.11e-06 & 7.78e-04 & 7.77e-04 & 1.62e-06 \\
4096 & 7.73e-06 & 3.72e-04 & 7.24e-06 & 1.04e-06 \\
8192 & 7.92e-06 & 1.01e-04 & 2.96e-05 & 5.59e-07 \\
16384 & 7.93e-06 & 3.84e-05 & 8.43e-06 & 7.41e-07 \\
\hline
\end{tabular}

\end{subtable}
\end{table}

\subsubsection{GSA using Sobol' indices}\label{sec:totalSobolindicies_ii}
\begin{figure}
    \centering
    \begin{subfigure}{0.4\linewidth}
        \includegraphics[trim=560pt 185pt 560pt 175pt, clip, width=1.0\linewidth]{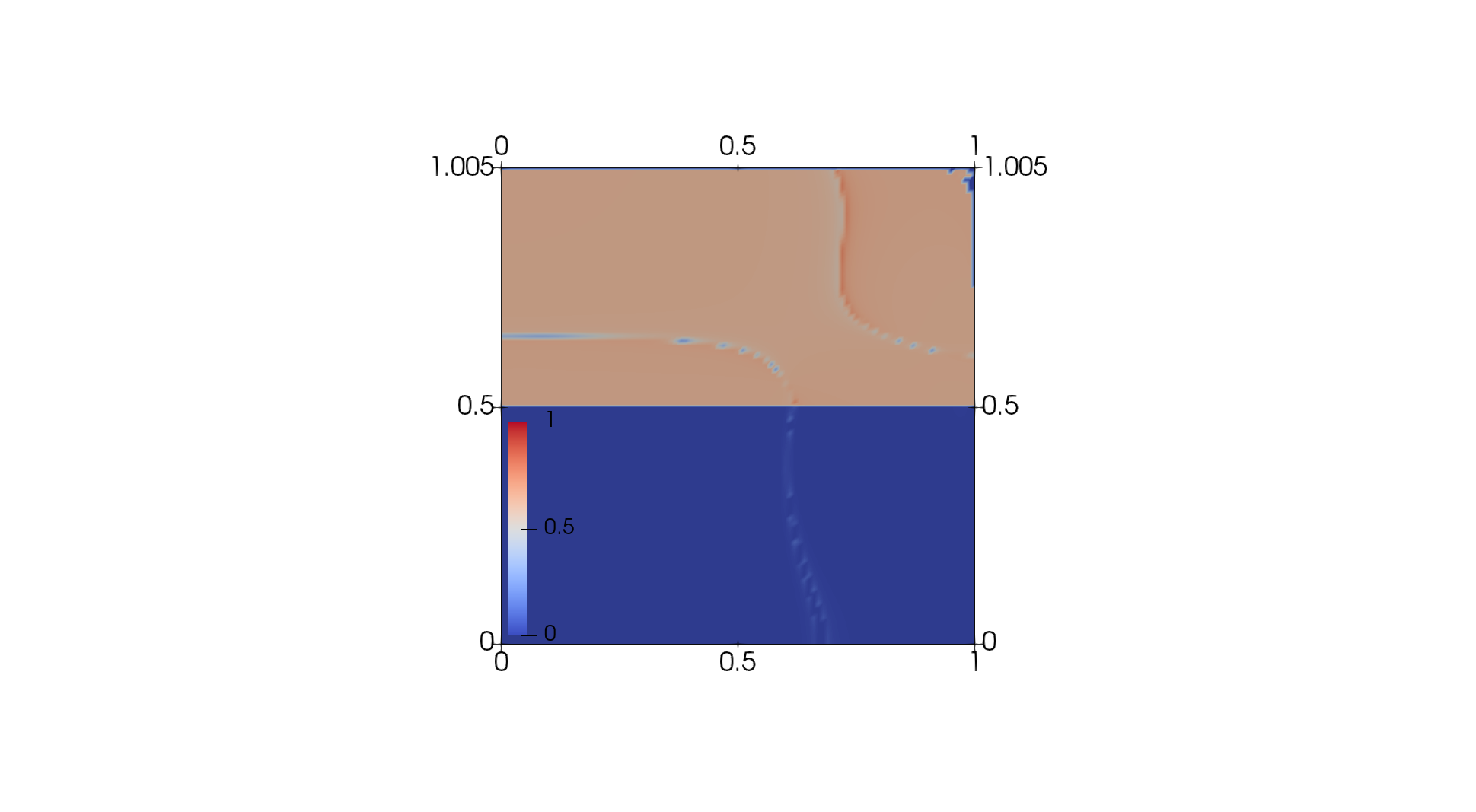}
    \caption{Velocity $u$, input parameter $I=1$}
    \end{subfigure}
    \hfil
    \begin{subfigure}{0.4\linewidth}
        \includegraphics[trim=560pt 185pt 560pt 175pt, clip, width=1.0\linewidth]{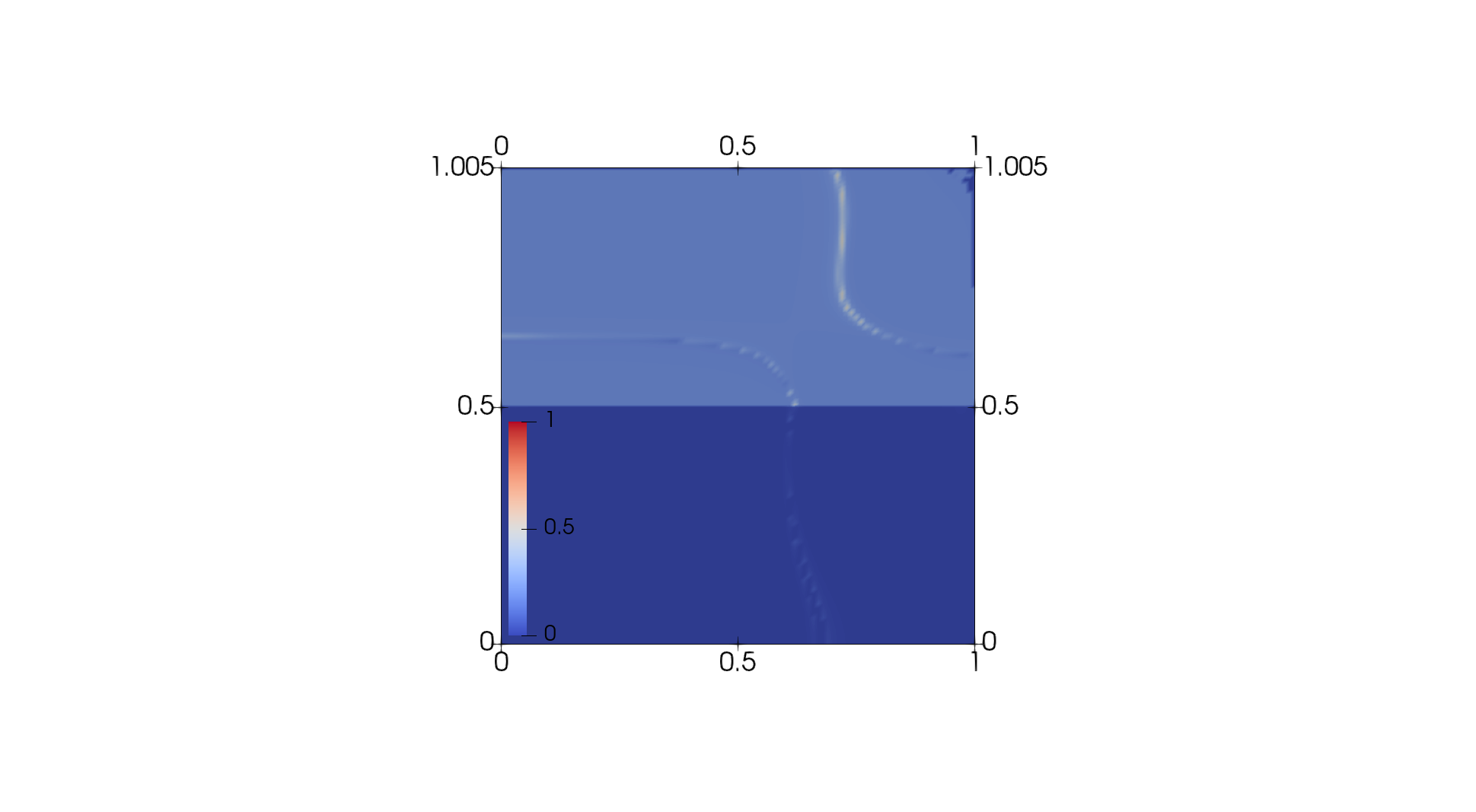}
    \caption{Velocity $u$, input parameter $I=2$}
    \end{subfigure}\\[4mm]
    
    \begin{subfigure}{0.4\linewidth}
        \includegraphics[trim=560pt 185pt 560pt 175pt, clip, width=1.0\linewidth]{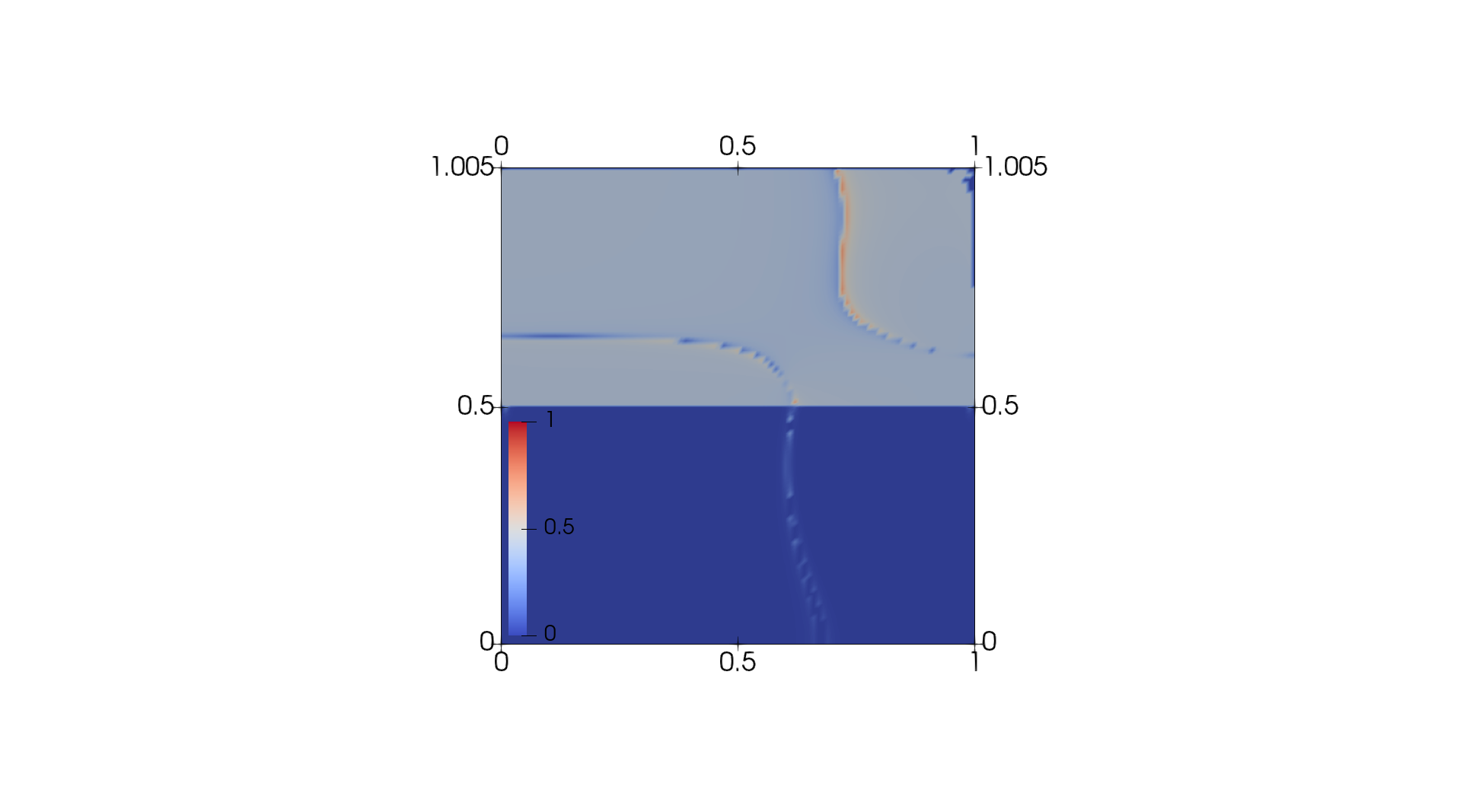}
    \caption{Velocity $u$, input parameter $I=3$}
    \end{subfigure}
    \hfil
    \begin{subfigure}{0.4\linewidth}
        \includegraphics[trim=560pt 185pt 560pt 175pt, clip, width=1.0\linewidth]{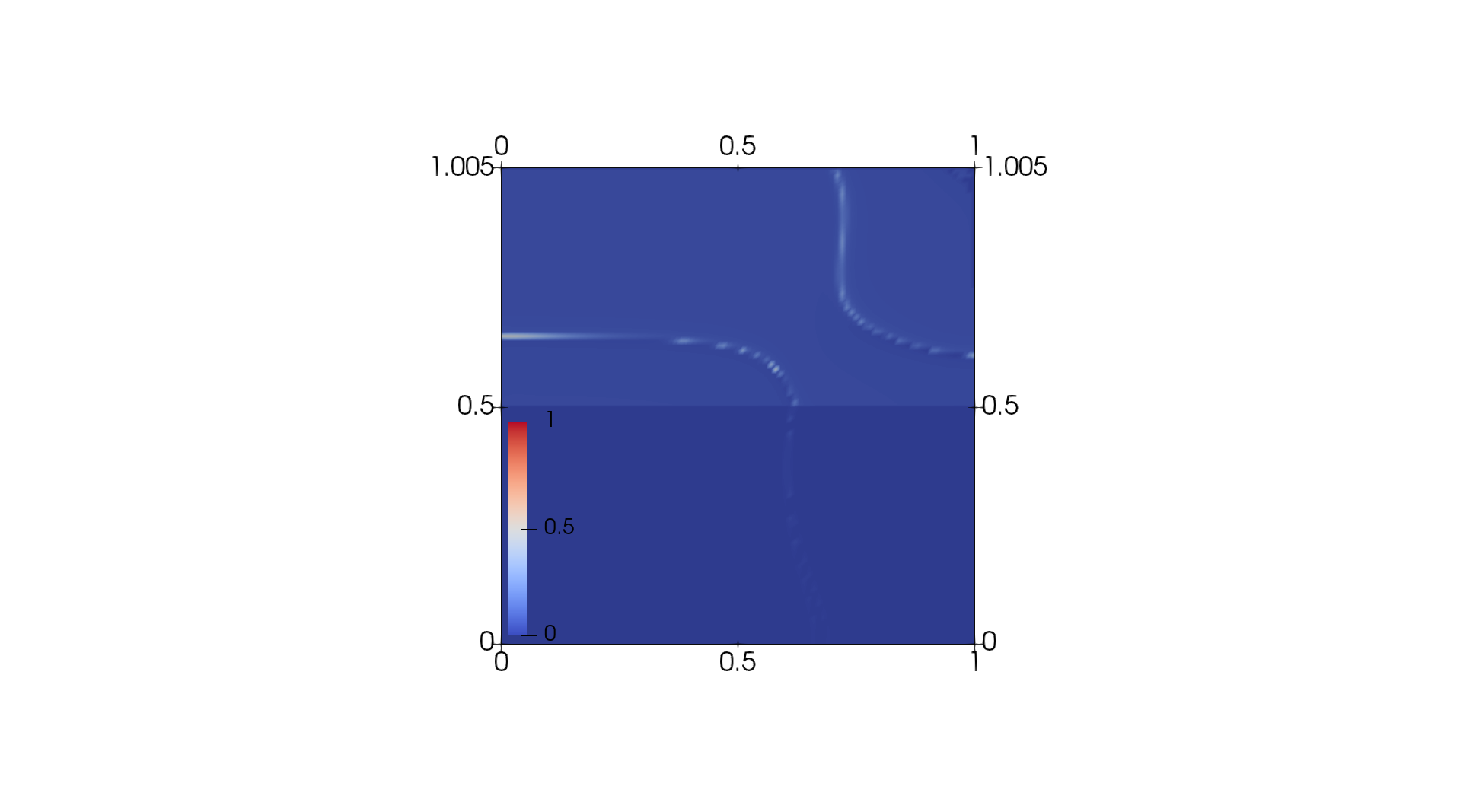}
    \caption{Velocity $u$, input parameter $I=4$}
    \end{subfigure}\\[4mm]

    \begin{subfigure}{0.4\linewidth}
        \includegraphics[trim=560pt 185pt 560pt 175pt, clip, width=1.0\linewidth]{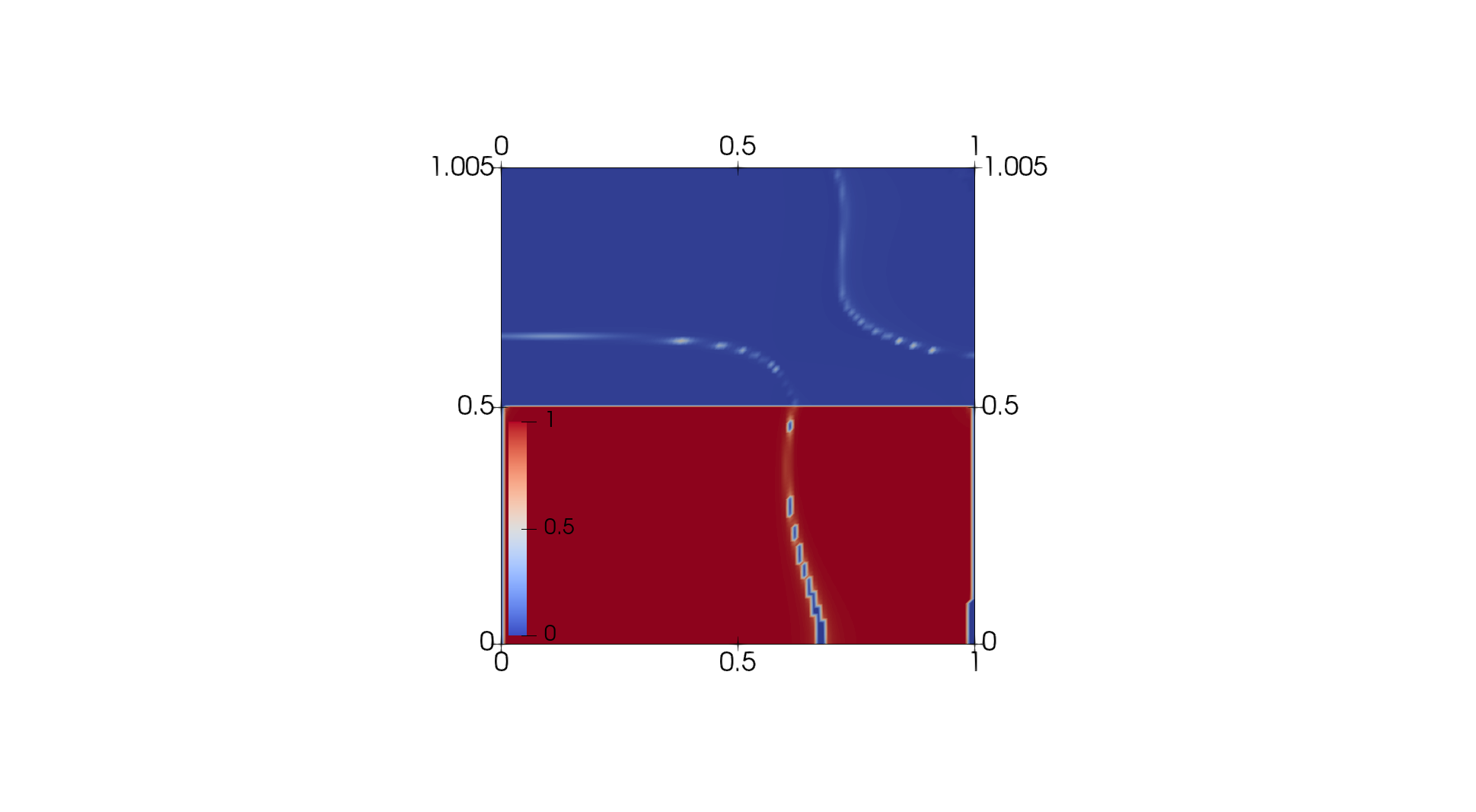}
    \caption{Velocity $u$, input parameter $I=5$}
    \end{subfigure}
   
    \caption{Case 2: Total sensitivity indices of the horizontal velocity $u$ obtained from \aMR-surrogate with $\Nr=1, \No=2$, trained on 8192 QMC samples.}
    \label{fig:amr_sobol-u_ii}
\end{figure}

\begin{figure}
    \centering
    \begin{subfigure}{0.4\linewidth}
        \includegraphics[trim=560pt 185pt 560pt 175pt, clip, width=1.0\linewidth]{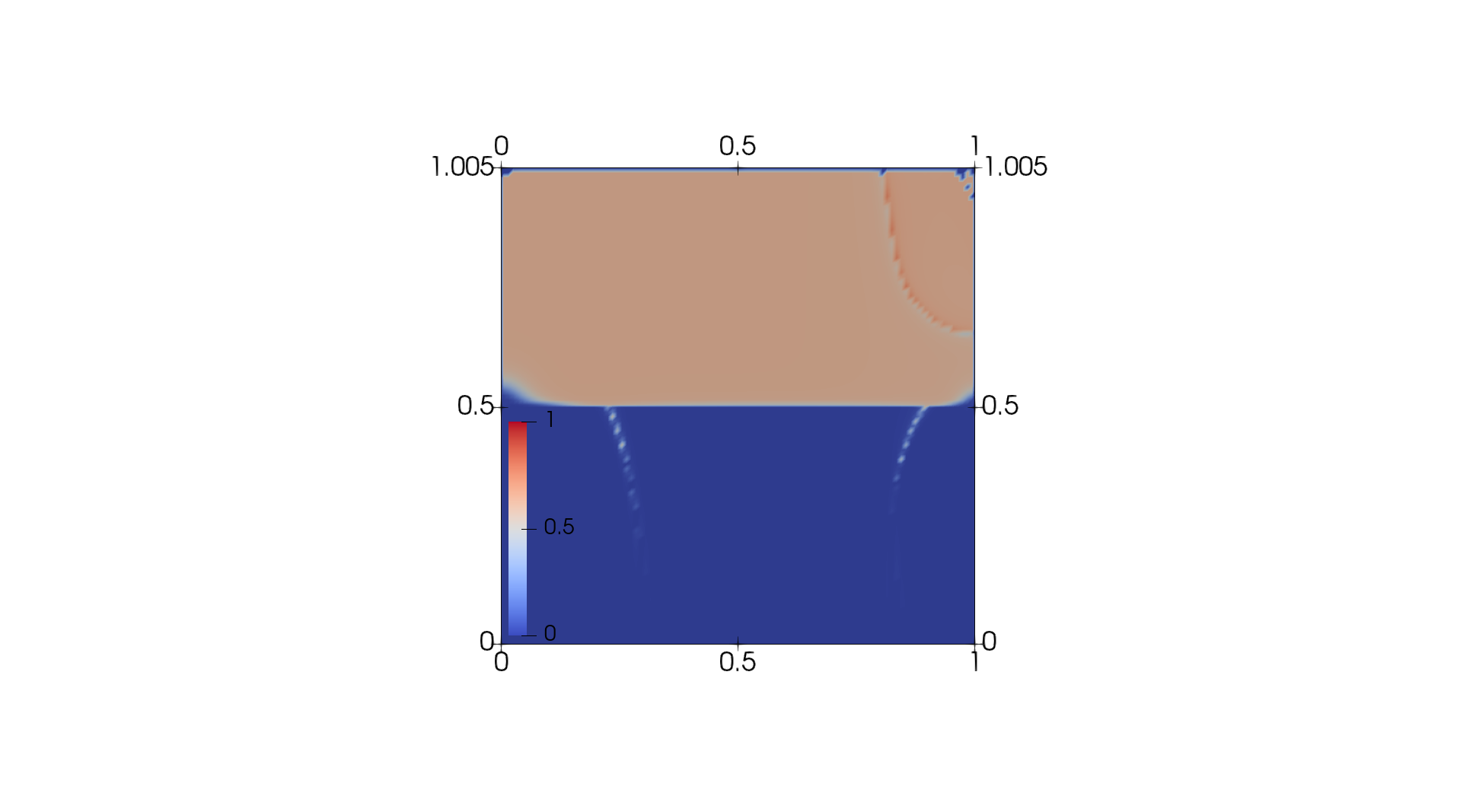}
    \caption{Velocity $v$, input parameter $I=1$}
    \end{subfigure}
    \hfil
    \begin{subfigure}{0.4\linewidth}
        \includegraphics[trim=560pt 185pt 560pt 175pt, clip, width=1.0\linewidth]{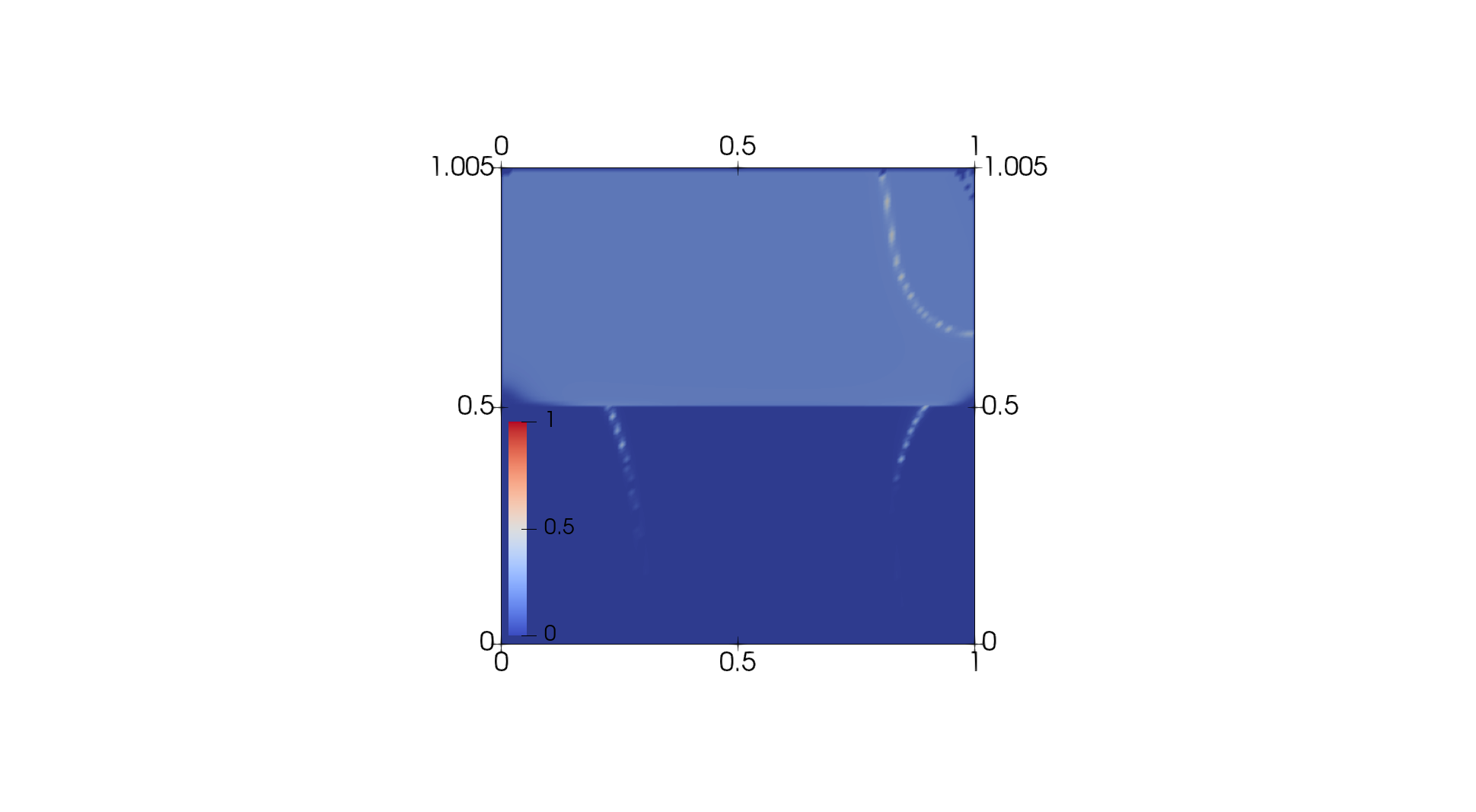}
    \caption{Velocity $v$, input parameter $I=2$}
    \end{subfigure}
    
    \begin{subfigure}{0.4\linewidth}
        \includegraphics[trim=560pt 185pt 560pt 175pt, clip, width=1.0\linewidth]{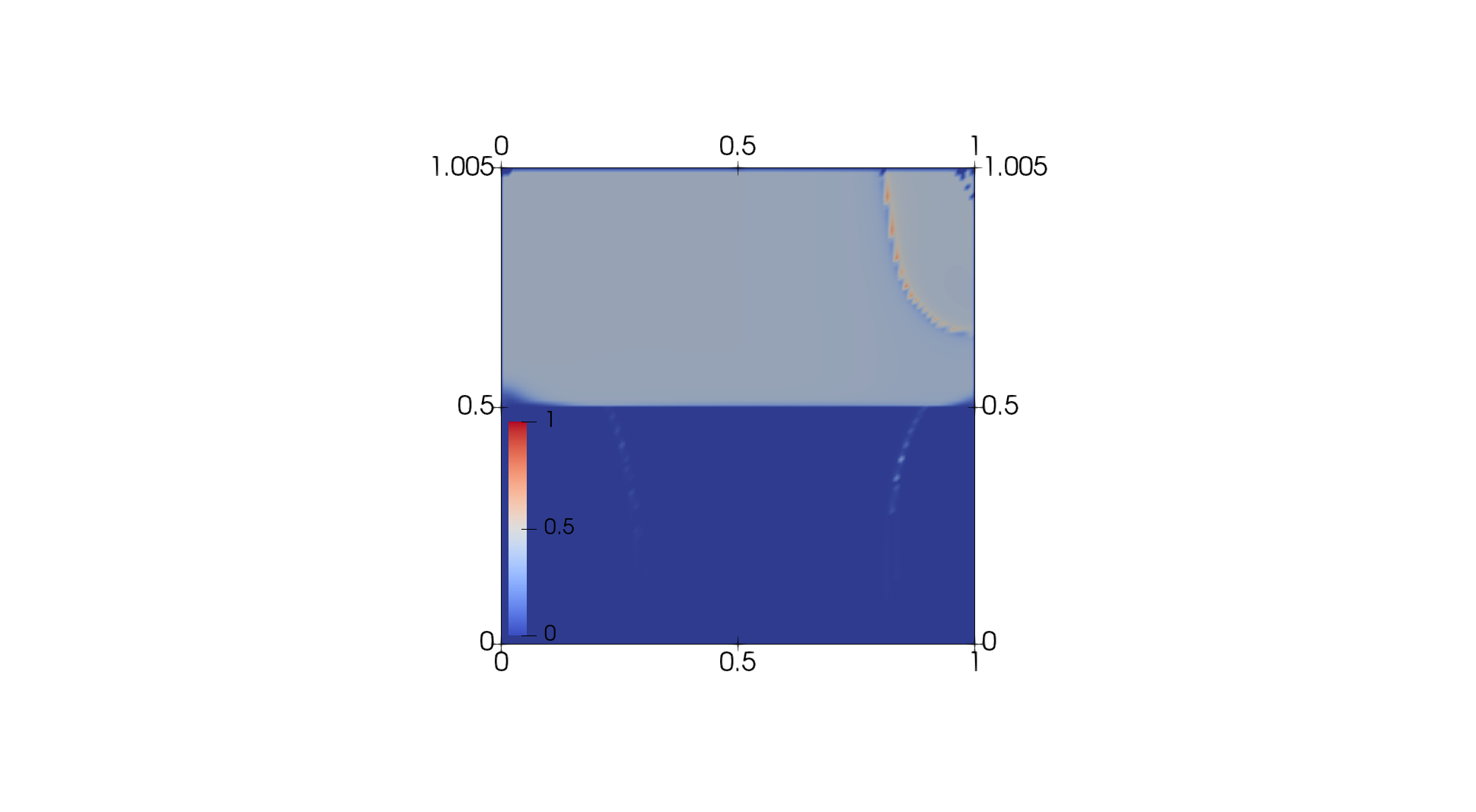}
    \caption{Velocity $v$, input parameter $I=3$}
    \end{subfigure}
    \hfil
    \begin{subfigure}{0.4\linewidth}
        \includegraphics[trim=560pt 185pt 560pt 175pt, clip, width=1.0\linewidth]{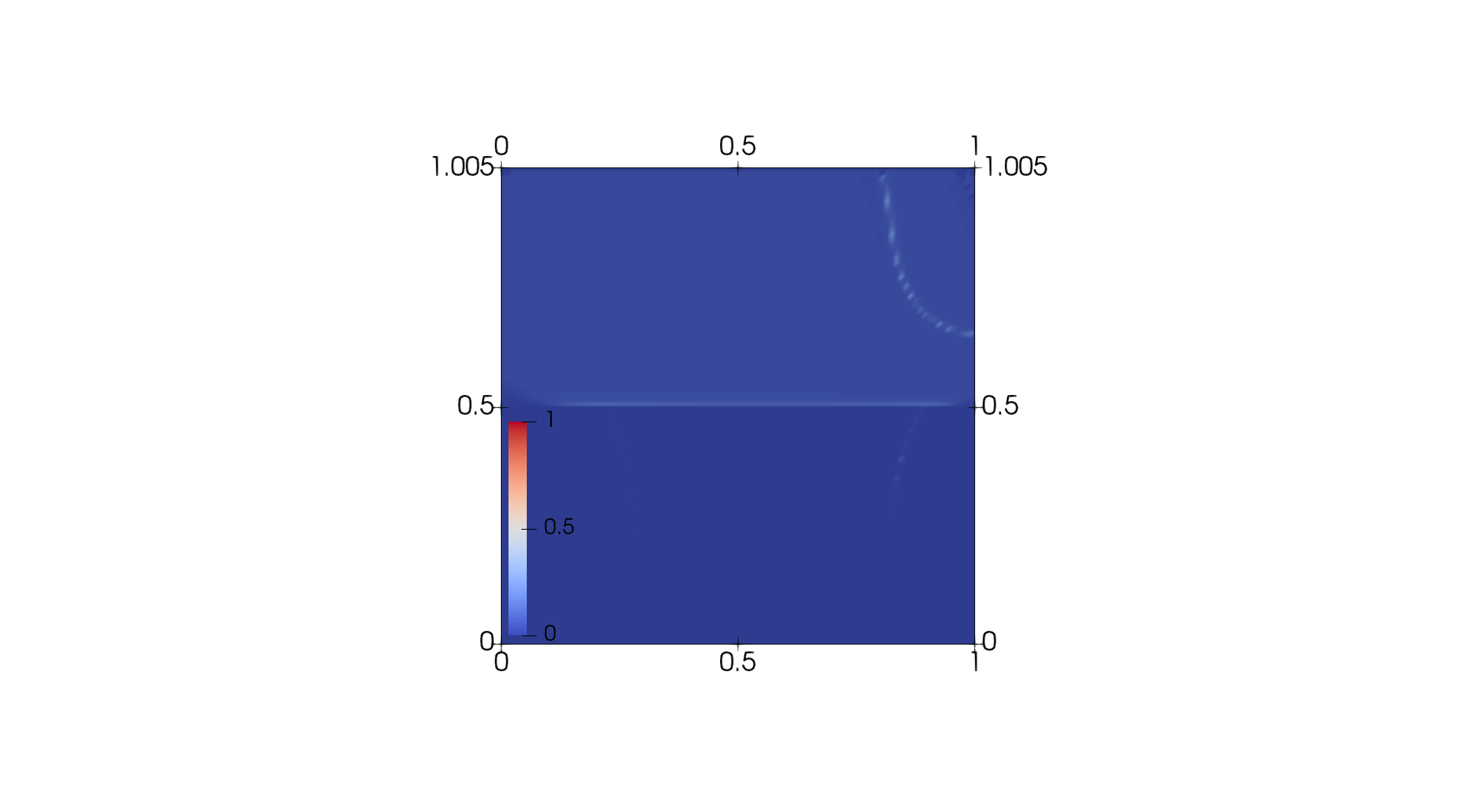}
    \caption{Velocity $v$, input parameter $I=4$}
    \end{subfigure}

    \begin{subfigure}{0.4\linewidth}
        \includegraphics[trim=560pt 185pt 560pt 175pt, clip, width=1.0\linewidth]{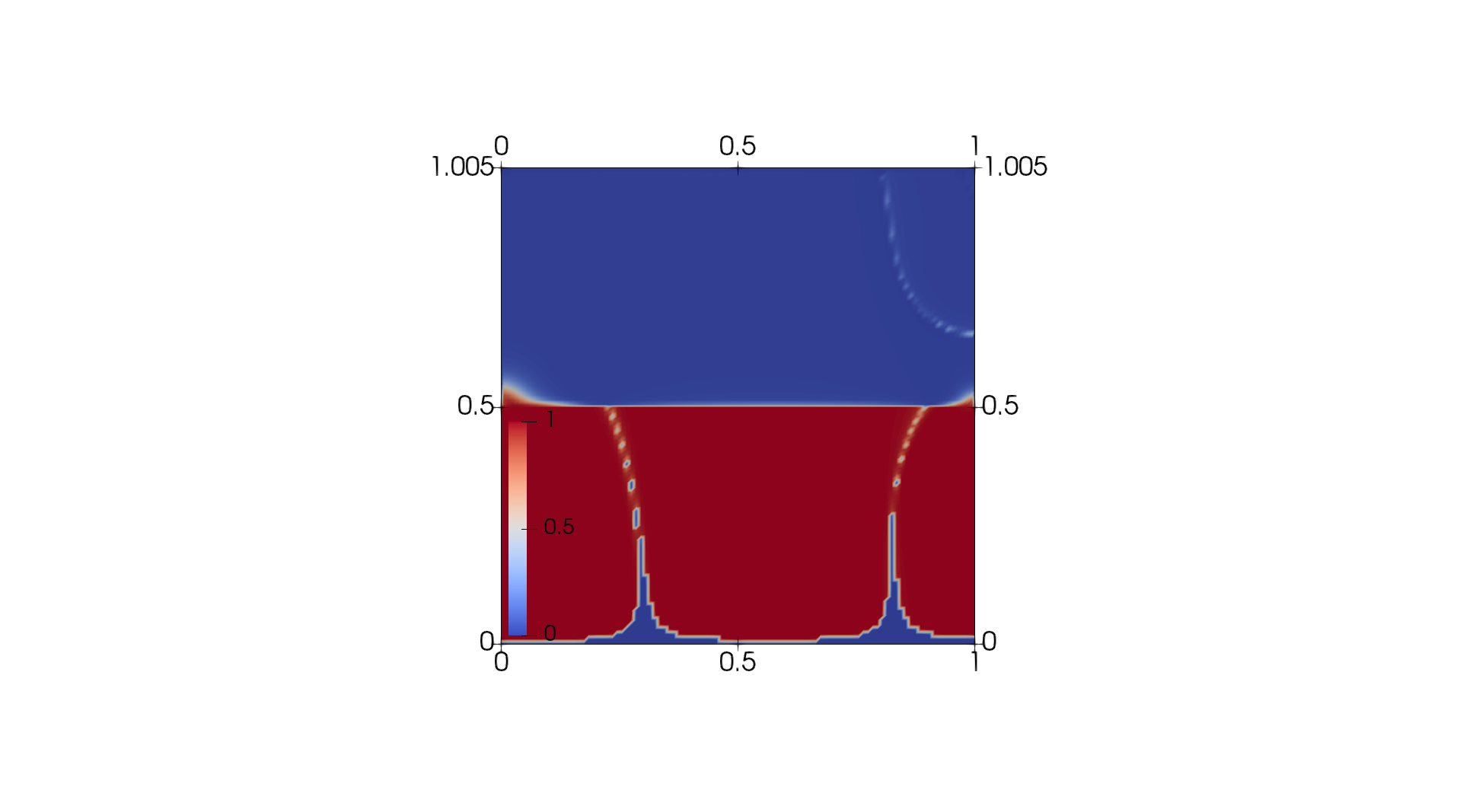}
    \caption{Velocity $v$, input parameter $I=5$}
    \end{subfigure}
   
    \caption{Case 2: Total sensitivity indices 
    of the vertical velocity $v$
    obtained from \aMR-surrogate with $\Nr=1, \No=2$,
    trained on 8192 QMC samples.}
    \label{fig:amr_sobol-v_ii}
\end{figure}

Also in this test case, the \aMR{}-based surrogate with $\Nr=1,\,\No=2$ 
demonstrates a superior accuracy performance.
Therefore, this surrogate is used for the GSA in this section.
Analyzing the convergence behavior of the averaged total
sensitivity indices presented in Fig.~\ref{fig:amr_sobol-u_ii}
and Fig.~\ref{fig:amr_sobol-v_ii} for the velocity components $u$ and $v$, respectively,
we observe only minimal changes after 4096 training samples. 
By examining the total sensitivity indices presented in Tab.~\ref{tab:Sobol-total_ii} along with the Sobol' sensitivity indices shown in Tab.~\ref{tab:Sobol-u_ii} and Tab.~\ref{tab:Sobol-v_ii} for $u$ and $v$, respectively, it is evident that the outcomes derived from the expansion coefficients of the aPC-based surrogate model with $\No=2$ align with the results computed using \aMR{} with $\Nr=1,\, \No=2$. This alignment serves as an indirect confirmation of these findings.
As in Case 1, hyperbolic truncation yields unreliable indices when computed from truncated expansion coefficients.

We again compute first-order and total-effect Sobol' indices (parameters as in Tab.~\ref{table1:UQtable}) to quantify individual and interaction effects for the splitting flow problem. 
The space distribution of the total Sobol' indices of both velocity components $u$ and $v$, computed using \aMR{} with $\Nr=1$, $\No=2$ with $8192$ samples, are visualized in Fig.~\ref{fig:amr_sobol-u_ii} and Fig.\ref{fig:amr_sobol-v_ii}. These results reveal that both velocity components exhibit highly nonlinear behavior with respect to variations in all five uncertain parameters in the splitting flow problem. In particular, sharp fronts emerge in both the free-flow and porous-medium regions. This nonlinearity is attributed to the specific boundary conditions imposed in the splitting flow problem.

\begin{figure}
    \centering
    \begin{subfigure}{.45\linewidth}
    \centering
    \includegraphics[width=0.65\linewidth]{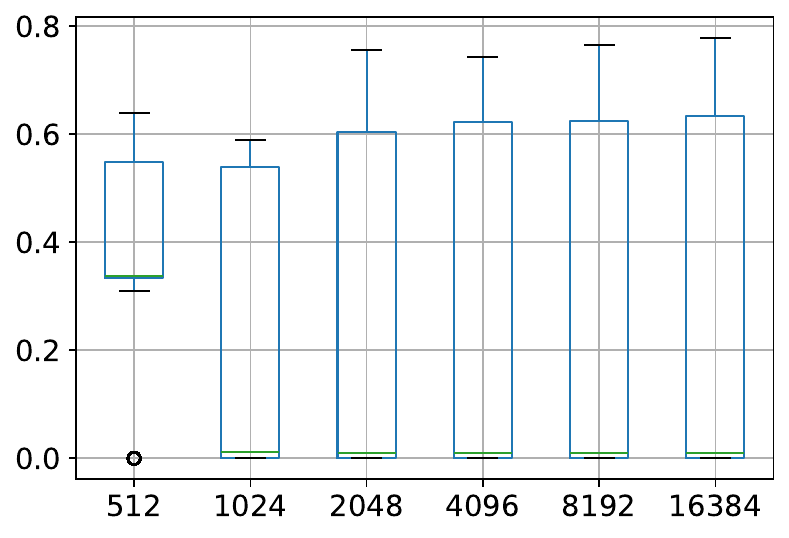}
    \caption{Input parameter $I=1$}
    \end{subfigure}
    \hfil
    \begin{subfigure}{.45\linewidth}
    \centering
    \includegraphics[width=0.65\linewidth]{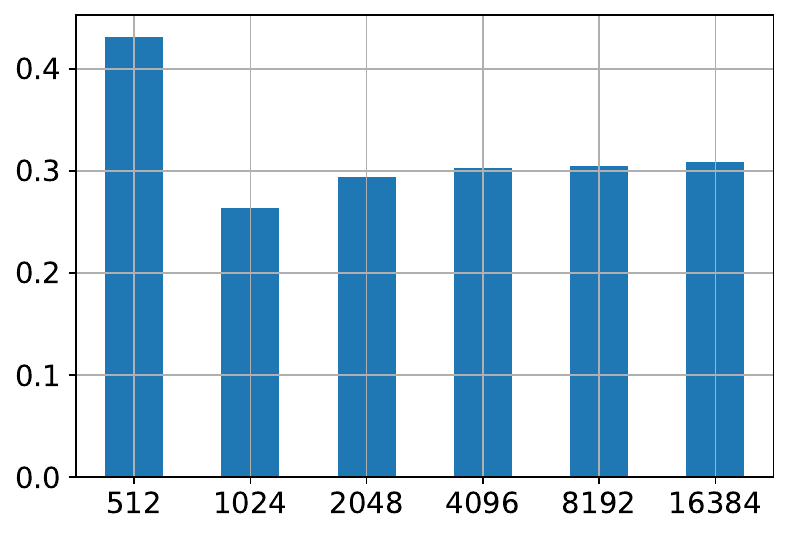}
    \caption{Input parameter $I=1$}
    \end{subfigure}

    \begin{subfigure}{.45\linewidth}
    \centering
    \includegraphics[width=0.65\linewidth]{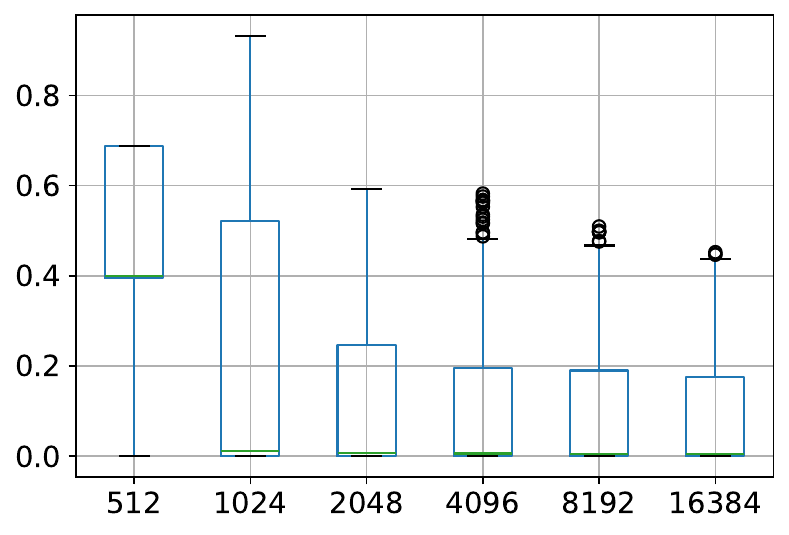}
    \caption{Input parameter $I=2$}
    \end{subfigure}
    \hfil
    \begin{subfigure}{.45\linewidth}
    \centering
    \includegraphics[width=0.65\linewidth]{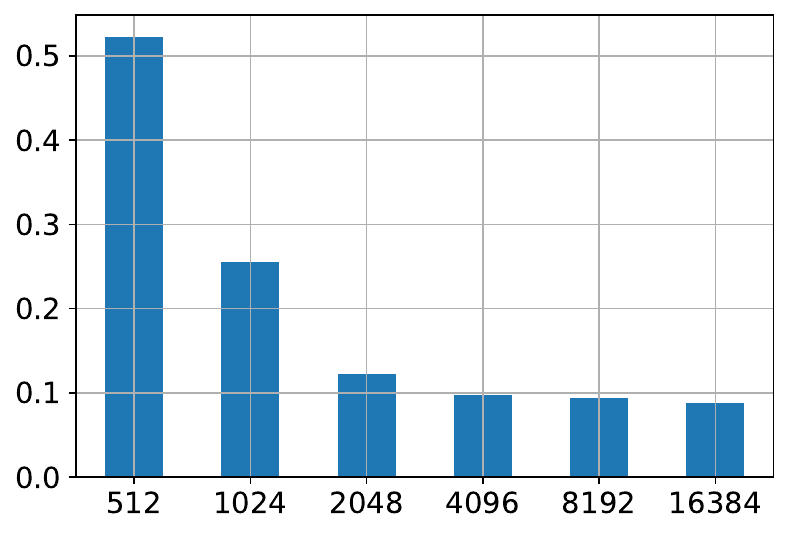}
    \caption{Input parameter $I=2$}
    \end{subfigure}

    \begin{subfigure}{.45\linewidth}
    \centering
    \includegraphics[width=0.65\linewidth]{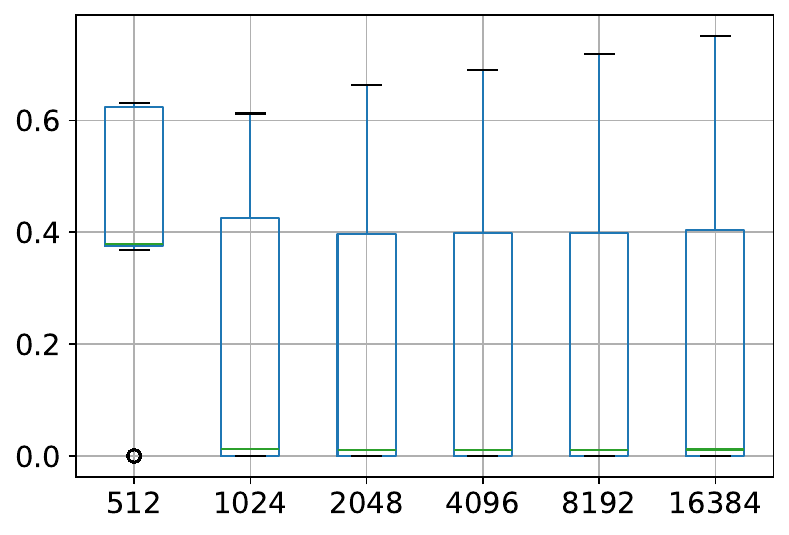}
    \caption{Input parameter $I=3$}
    \end{subfigure}
    \hfil
    \begin{subfigure}{.45\linewidth}
    \centering
    \includegraphics[width=0.65\linewidth]{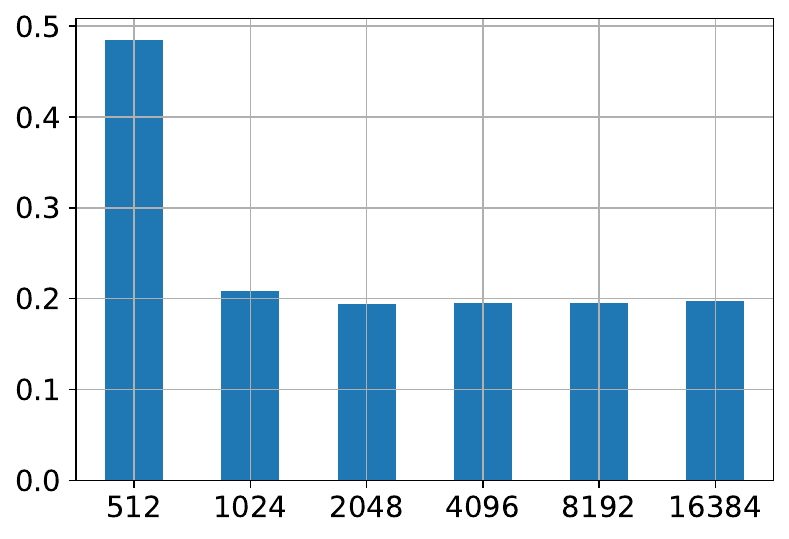}
    \caption{Input parameter $I=3$}
    \end{subfigure}

    \begin{subfigure}{.45\linewidth}
    \centering
    \includegraphics[width=0.65\linewidth]{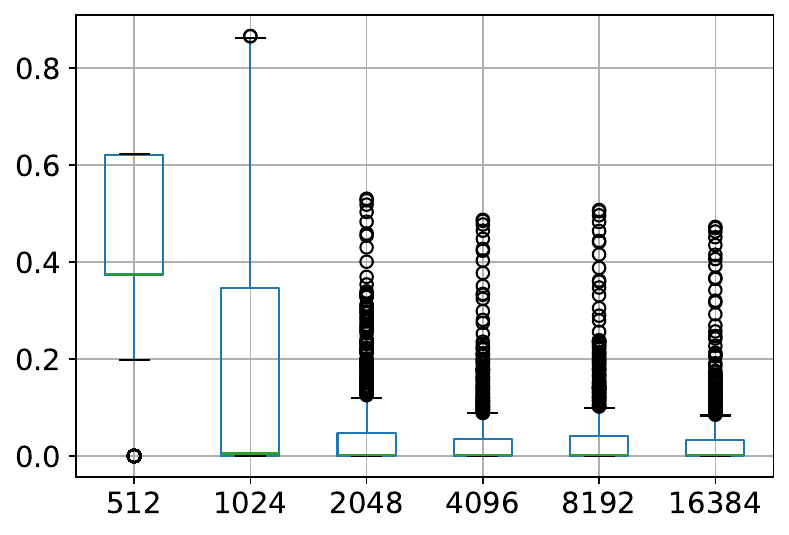}
    \caption{Input parameter $I=4$}
    \end{subfigure}
    \hfil
    \begin{subfigure}{.45\linewidth}
    \centering
    \includegraphics[width=0.65\linewidth]{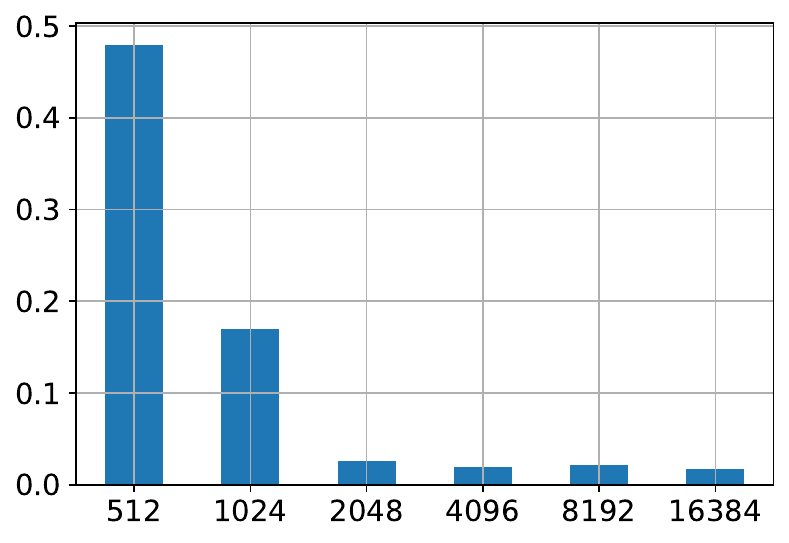}
    \caption{Input parameter $I=4$}
    \end{subfigure}

    \begin{subfigure}{.45\linewidth}
    \centering
    \includegraphics[width=0.65\linewidth]{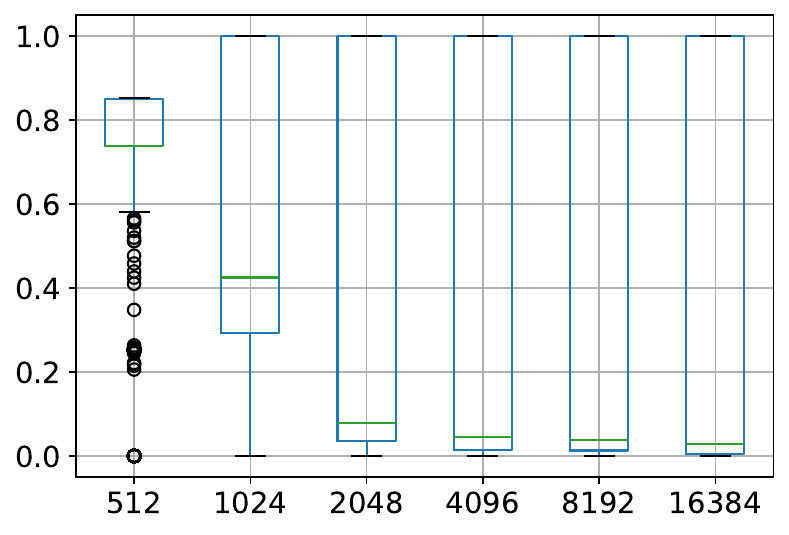}
    \caption{Input parameter $I=5$}
    \end{subfigure}
    \hfil
    \begin{subfigure}{.45\linewidth}
    \centering
    \includegraphics[width=0.65\linewidth]{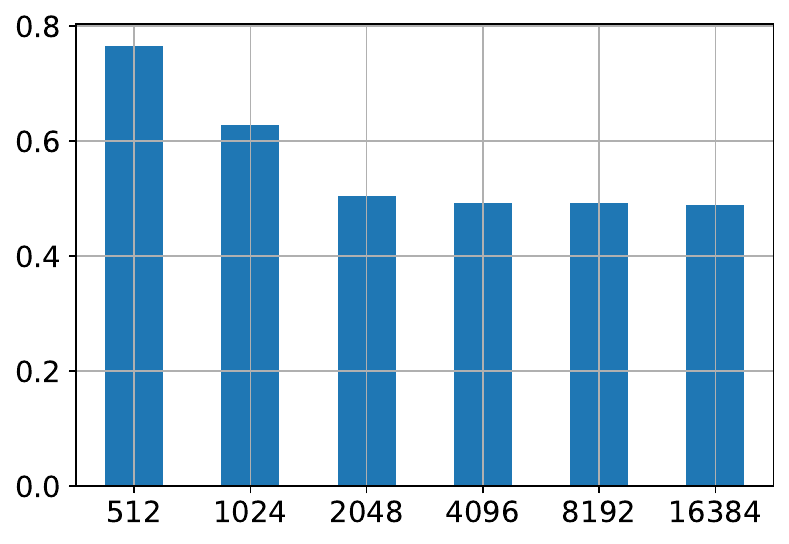}
    \caption{Input parameter $I=5$}
    \end{subfigure}
    
    \caption{Case 2: Evolution of the total sensitivity indices of the horizontal velocity $u$ obtained
    from \aMR-surrogate with $\Nr=1,\, \No=2$.}
    \label{fig:tot_evo_u_ii}
\end{figure}

\begin{figure}
    \centering
    \begin{subfigure}{.45\linewidth}
    \centering
    \includegraphics[width=0.65\linewidth]{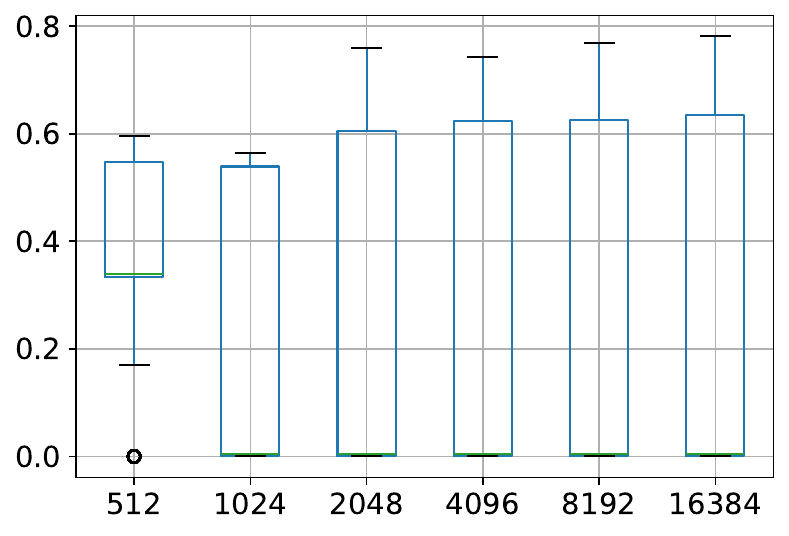}
    \caption{Input parameter $I=1$}
    \end{subfigure}
    \hfil
    \begin{subfigure}{.45\linewidth}
    \centering
    \includegraphics[width=0.65\linewidth]{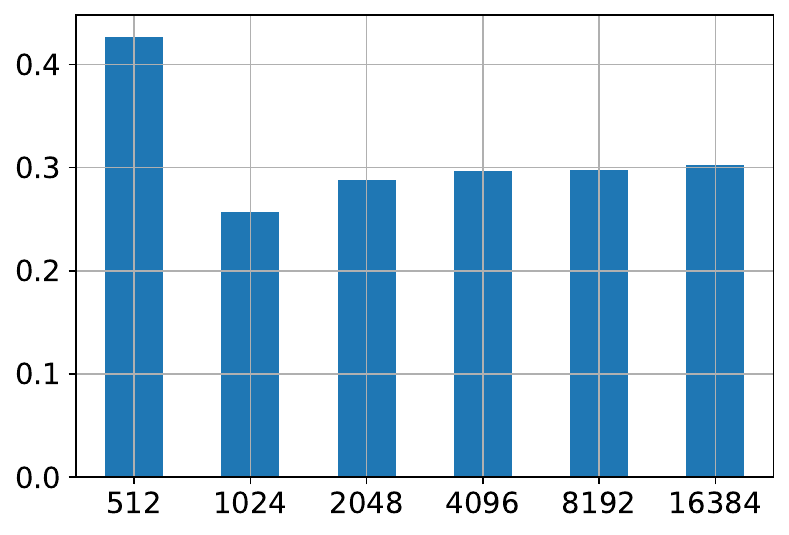}
    \caption{Input parameter $I=1$}
    \end{subfigure}

    \begin{subfigure}{.45\linewidth}
    \centering
    \includegraphics[width=0.65\linewidth]{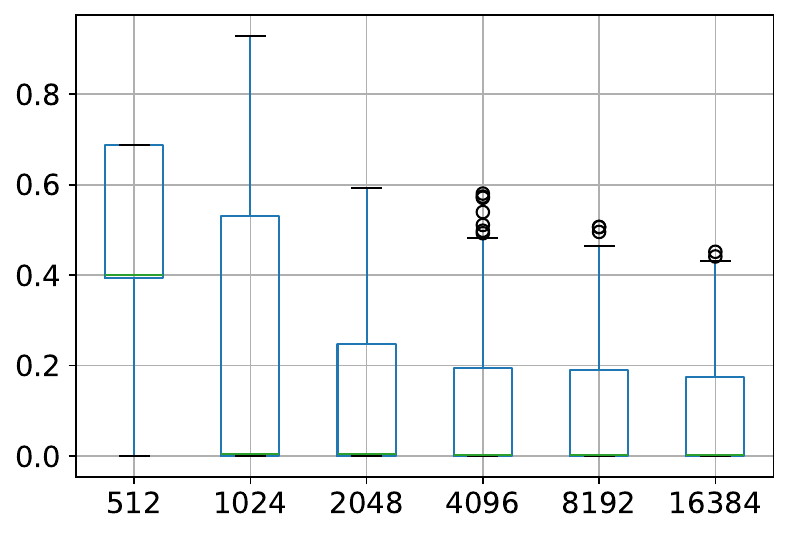}
    \caption{Input parameter $I=2$}
    \end{subfigure}
    \hfil
    \begin{subfigure}{.45\linewidth}
    \centering
    \includegraphics[width=0.65\linewidth]{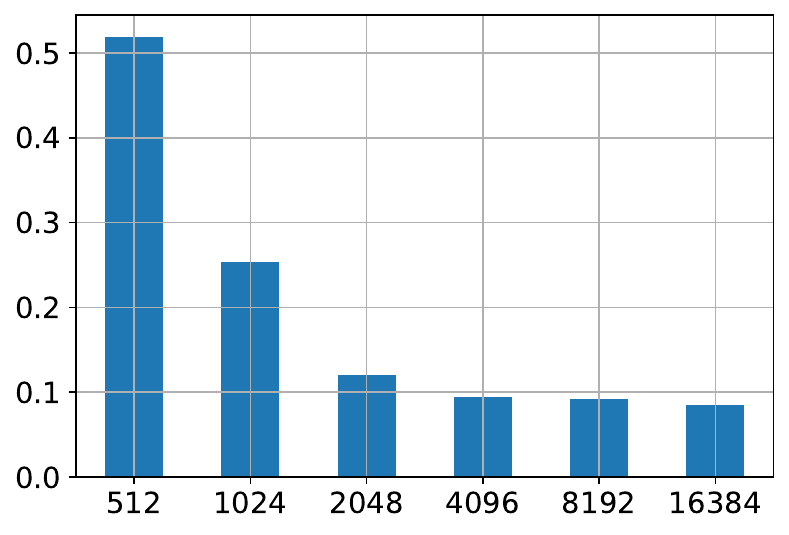}
    \caption{Input parameter $I=2$}
    \end{subfigure}

    \begin{subfigure}{.45\linewidth}
    \centering
    \includegraphics[width=0.65\linewidth]{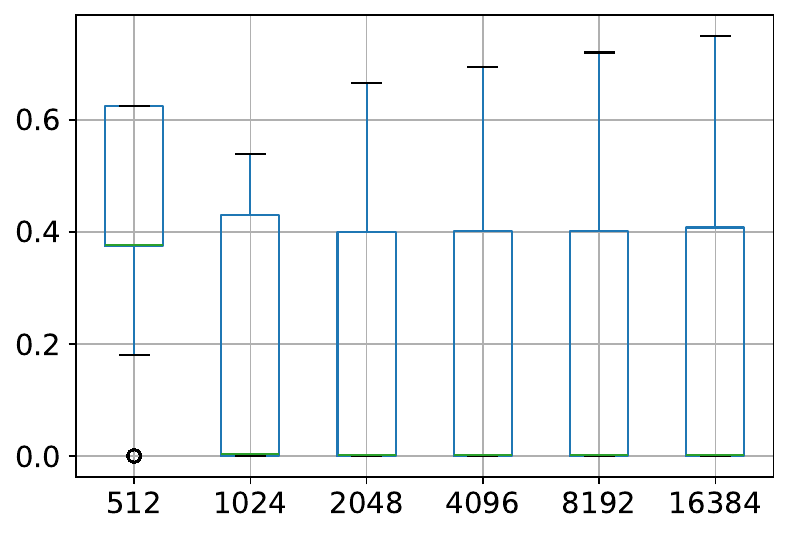}
    \caption{Input parameter $I=3$}
    \end{subfigure}
    \hfil
    \begin{subfigure}{.45\linewidth}
    \centering
    \includegraphics[width=0.65\linewidth]{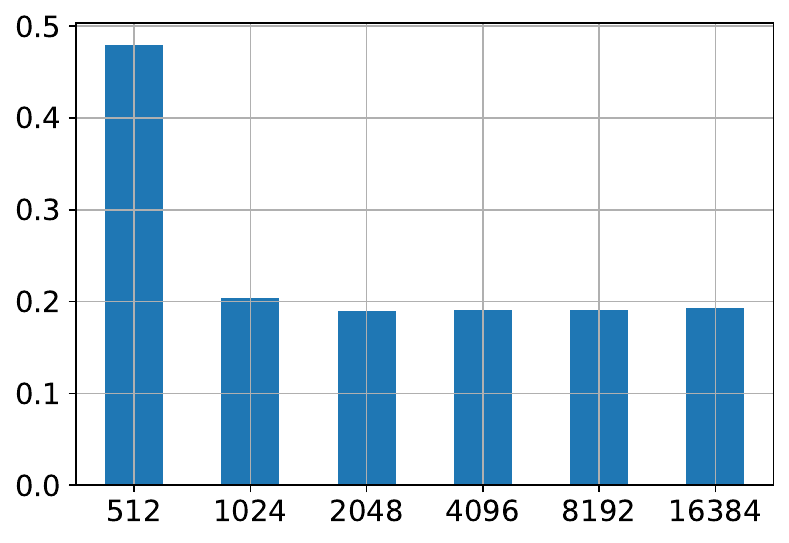}
    \caption{Input parameter $I=3$}
    \end{subfigure}

    \begin{subfigure}{.45\linewidth}
    \centering
    \includegraphics[width=0.65\linewidth]{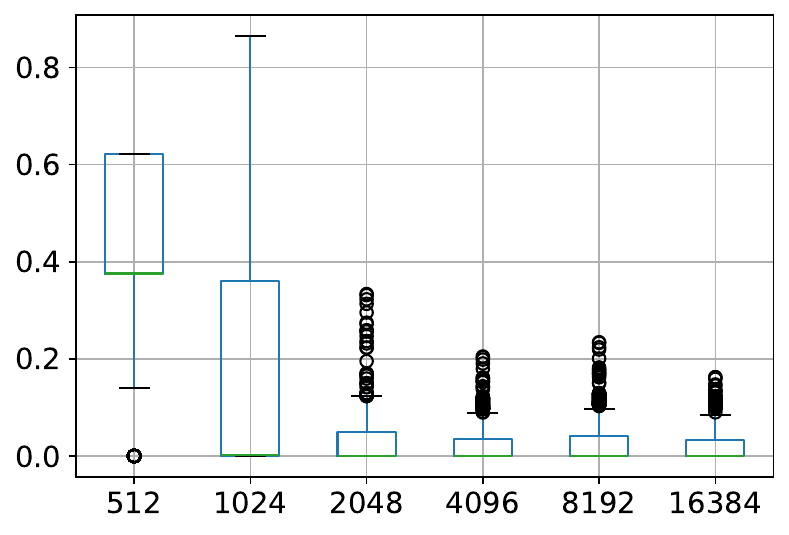}
    \caption{Input parameter $I=4$}
    \end{subfigure}
    \hfil
    \begin{subfigure}{.45\linewidth}
    \centering
    \includegraphics[width=0.65\linewidth]{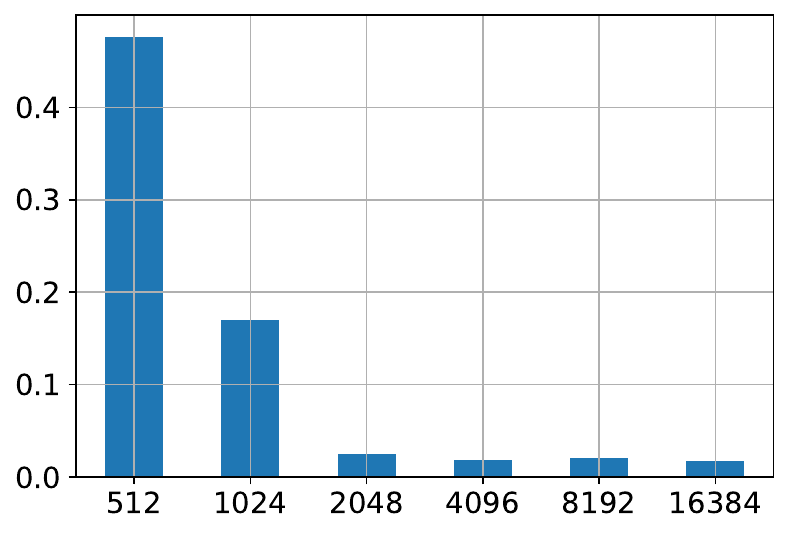}
    \caption{Input parameter $I=4$}
    \end{subfigure}

    \begin{subfigure}{.45\linewidth}
    \centering
    \includegraphics[width=0.65\linewidth]{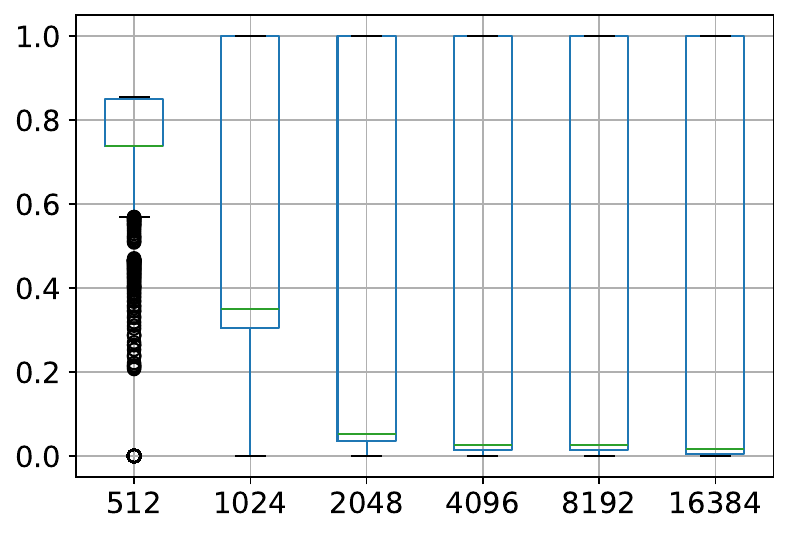}
    \caption{Input parameter $I=5$}
    \end{subfigure}
    \hfil
    \begin{subfigure}{.45\linewidth}
    \centering
    \includegraphics[width=0.65\linewidth]{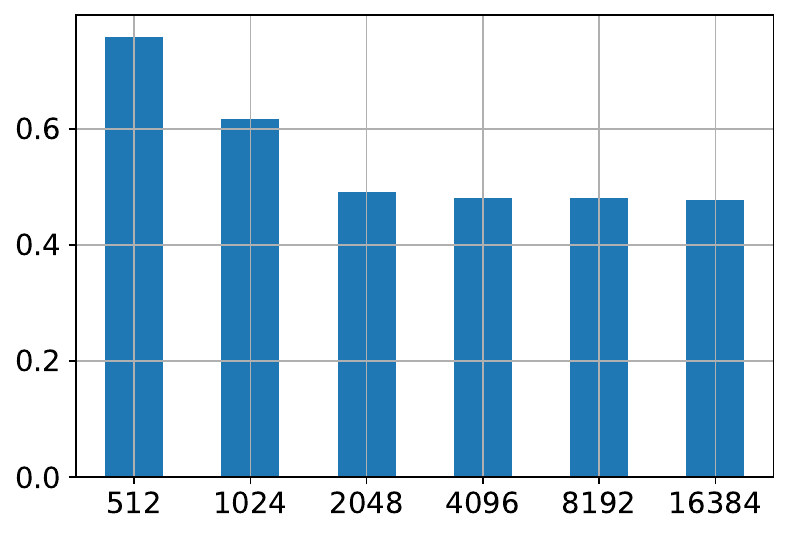}
    \caption{Input parameter $I=5$}
    \end{subfigure}
    
    \caption{Case 2: Evolution of the total sensitivity indices of the vertical velocity $v$ obtained
    from \aMR-surrogate with $\Nr=1,\, \No=2$.}
    \label{fig:tot_evo_v_ii}
\end{figure}

\begin{figure}
    \centering

    \begin{subfigure}{0.4\linewidth}
    \includegraphics[width=0.95\linewidth]{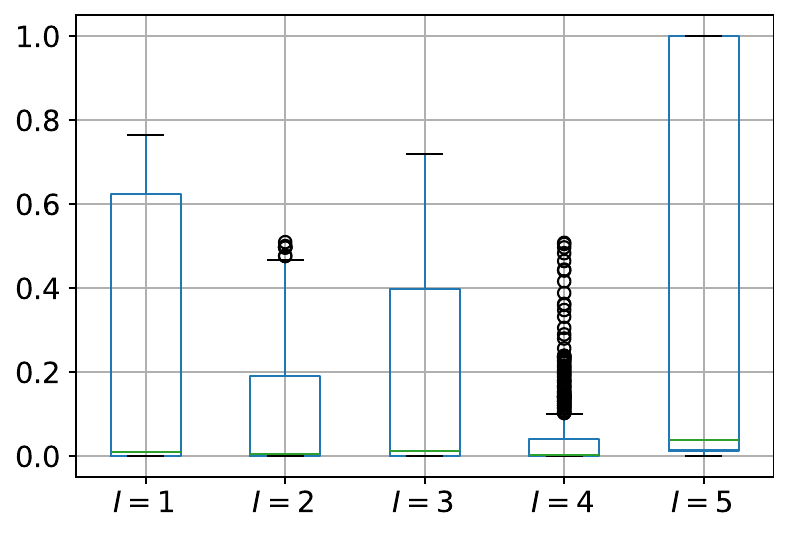}
    \caption{Box plot of total indices for velocity $u$}
    \end{subfigure}
    \hfil 
    \begin{subfigure}{0.4\linewidth}
    \includegraphics[width=0.95\linewidth]{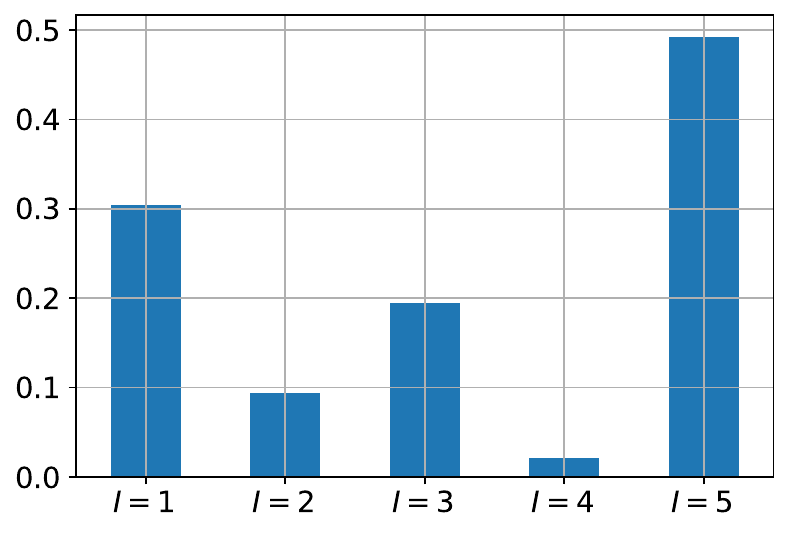}
    \caption{Averaged total indices for velocity $u$}
    \end{subfigure}

    \begin{subfigure}[t]{0.4\linewidth}
    \includegraphics[width=0.95\linewidth]{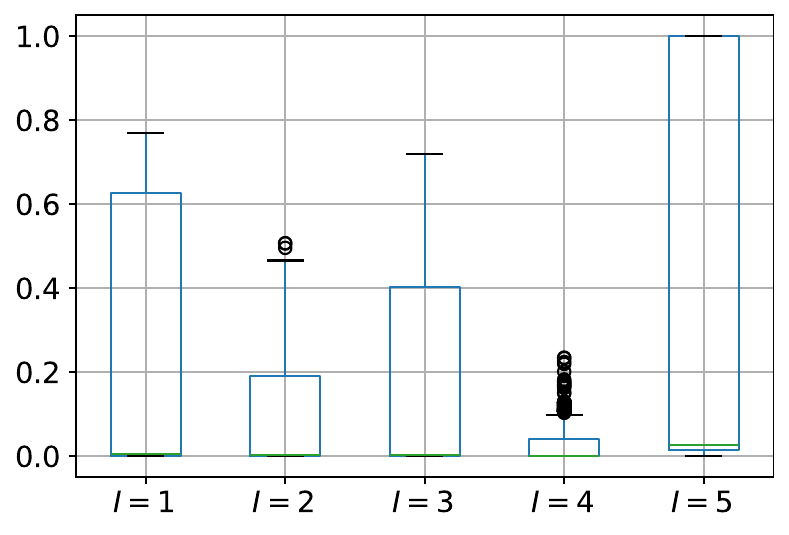}
    \caption{Box plot of total indices for velocity $v$}
    \end{subfigure}
    \hfil 
    \begin{subfigure}[t]{0.4\linewidth}
    \includegraphics[width=0.95\linewidth]{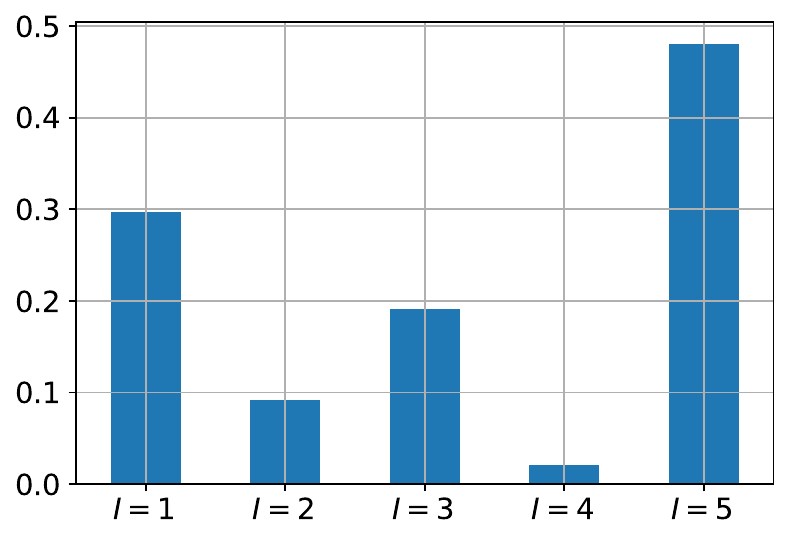}
    \caption{Averaged total indices of velocity $v$}
    \end{subfigure}
    
    \caption{Case 2: Evaluation of the total sensitivity indices 
    over the whole domain,
    obtained from \aMR-surrogate with $\Nr=1, \No=2$,
    trained on 8192 QMC samples.
    }
    \label{fig:sobol-total_av_ii}
\end{figure}

The stress jump parameter $\beta_\JU$ and the permeability in the complex interface $k_\Gamma$ predominantly influence both velocity components in the free-flow domain (Fig.~\ref{fig:amr_sobol-u_ii}(A,B), Fig.~\ref{fig:amr_sobol-v_ii}(A,B)). In the splitting flow problem, the flow is driven and exits through the free-flow region, where the velocities are considerably larger than in the porous-medium regions (Fig.~\ref{fig:testcase2SFgeo}, right). As a result, the flow in the free-flow subdomain is particularly sensitive to changes in the stress jump parameter and the interface permeability.  Figures~\ref{fig:amr_sobol-u_ii}(C) and \ref{fig:amr_sobol-v_ii}(C) indicate that the effective viscosity $\mu_\EF$ mainly affects the free-flow domain as in the filtration problem in Section~\ref{sec:totalSobolindicies_i}.

The Beavers–Joseph parameter $\alpha_\BJ$ has a comparatively smaller effect than the other parameters, with most variations concentrated along the complex interface (Fig.~\ref{fig:amr_sobol-u_ii}(D), Fig.~\ref{fig:amr_sobol-v_ii}(D)). This behavior reflects the role of $\alpha_\BJ$ in characterizing the interface roughness at the bottom of the transition region. As can be seen from Fig.~\ref{fig:amr_sobol-u_ii}(E) and Fig.~\ref{fig:amr_sobol-v_ii}(E), both velocity components in the porous medium are particularly sensitive to variations in the porous-medium permeability $k_\PM$. This is expected, since the relatively small porous-medium permeability directly controls how much flow can enter the porous region through the complex interface, consistent with \cite{kroker2023global}.

We note again that both the Sobol' indices and total Sobol' indices are normalized by variance (Fig.~\ref{fig:amr_mean_var_ii}(E,~F)). As a consequence, numerical errors can be amplified in regions where the variance is close to zero. Hence, the results in these regions are not necessarily meaningful from a modeling perspective.

Figures~\ref{fig:tot_evo_u_ii}, \ref{fig:tot_evo_v_ii} display the box plots (left) and the mean bar plots (right) of the Sobol’ indices for five uncertain parameters across different sample sizes ranging from $512$ to $16384$ in the splitting flow problem. 
Convergence of the Sobol’ sensitivity indices with respect to the number of samples is evident in both the box plots  and the averaged indices. From $2048$ samples onward, both the distributions and the averaged values exhibit stable behavior (Fig.~\ref{fig:tot_evo_u_ii}, \ref{fig:tot_evo_v_ii}).

We therefore select $8192$ samples as a representative case and present the corresponding box plots of the Sobol’ indices together with histograms of their space-averaged values in Fig.~\ref{fig:sobol-total_av_ii}. Compared to the filtration problem (Fig.~\ref{fig:sobol-total_av}(A,C)), the medians of the total indices are smaller across all uncertain parameters (Fig.~\ref{fig:sobol-total_av_ii}(A,C)). This can be explained by the fact that the influence of the parameters is primarily concentrated in one of the subdomains (Fig.~\ref{fig:amr_sobol-u_ii},~\ref{fig:amr_sobol-v_ii}). These results reveal that the newly derived hybrid-dimensional Stokes--Brinkman--Darcy model for the splitting flow problem is most sensitive to the porous-medium permeability~$k_\PM$ (Fig.~\ref{fig:sobol-total_av_ii}(B,D)). The stress jump parameter~$\beta_\JU$ and the effective viscosity~$\mu_\EF$ also have significant influence. In contrast, the Beavers--Joseph parameter~$\alpha_\BJ$ shows only minor impact and is not essential in this test case, although its box plot variability indicates occasional strong contributions in specific realizations (outliers), consistent with the previous findings for the filtration problem (Case~1).

\section{Conclusions}\label{sec:conclusions}
In the current study, we evaluated the appropriateness of various non-regularized PCE-based surrogate models for GSA when training samples are provided through a space-filling Sobol' sequence.
This scenario is particularly relevant if the training data
cannot be generated on demand and are instead provided by observation.
However, we find that surrogates with maximum polynomial degree $\No=2$ yield reliable GSA results. In particular, \aMR{} with $\Nr=1,\,\No=2$ achieves the highest accuracy in the comparison, since a sufficient quantity of training samples was available, making it the preferred choice for conducting GSA of the Stokes--Brinkman--Darcy flow model.

In this work, two test cases have been investigated for GSA: a filtration problem and a splitting flow problem. Analyzing the total sensitivity indices, we observed that the \aMR{} surrogate of the Stokes--Brinkman--Darcy flow model is highly sensitive to the stress jump parameter, the permeability of the complex interface, the effective viscosity, and the porous-medium permeability in both test cases. For the filtration problem (Case~1), the stress jump parameter has the most significant influence, whereas for the splitting flow problem (Case~2), the porous-medium permeability is the dominant factor. 

In both cases, the Beavers--Joseph coefficient does not significantly affect the overall flow behavior, but it influences the fluid behavior locally near the complex interface. The porous-medium permeability, within the considered parameter range, has no strong influence on the free-flow velocity. The influence of the Beavers--Joseph coefficient and the porous-medium permeability is consistent with the results reported in~\cite{kroker2023global}. The effective viscosity primarily impacts the free-flow region rather than the porous-medium region in both test cases. The sensitivity to the stress jump parameter and the interface permeability depends on the specific test case. In the filtration problem (Case~1), where the fluid velocities in the free-flow and porous-medium regions are comparable, the entire flow domain is highly sensitive to both parameters. In contrast, in the splitting flow problem (Case~2), where the free-flow velocity is relatively higher, both parameters are influential only in the free-flow region. 
 
Several directions remain open for future research. 
One possible extension of this work is to consider finer mesh resolutions. 
Although finer meshes would improve the accuracy of the reference solutions, they would also increase the computational effort.
Since the present study serves as a benchmark for surrogate models, the insights obtained here can guide future reductions in the number of required training samples and the set of surrogate configurations to compare. 
Another natural direction is to incorporate a larger set of input parameters in GSA, for instance by accounting for anisotropic porous media.
In this case, the stress jump tensor $\ten{\beta}$ as well as the permeability tensors $\ten{K}_\Gamma$ and $\ten{K}_\PM$ would be treated as fully anisotropic. 
Such an extension would allow for a more comprehensive characterization of uncertainties in the hybrid-dimensional Stokes--Brinkman--Darcy flow model and, at the same time, open opportunities for developing and testing advanced surrogate strategies capable of efficiently handling the expanded parameter space.

\subsection*{Data Availability Statement}
The research data associated with this article are available in DaRUS \cite{DARUS-5251_2025}.

\bibliographystyle{abbrv} 
\bibliography{manuscipt} 

@ARTICLE{Angot_etal_17,
  author = {Angot, P. and Goyeau, B. and Ochoa-Tapia, J. A.},
  title = {Asymptotic modeling of transport phenomena at the interface 
between a fluid and a porous layer: jump conditions},
  journal = {Phys. Rev. E},
  year = {2017},
  volume = {95},
  number ={6},
  pages = {063302},
  doi = {10.1103/PhysRevE.95.063302}
}

@article{Angot_etal_20,
  author = {Angot, P. and Goyeau, B. and Ochoa-Tapia, J. A.},
  title = {A nonlinear asymptotic model for the inertial flow at a fluid-porous interface},
  journal = {Adv. Water Res.},
  year = {2021},
  volume = {149},
  pages = {103798},
  doi = {10.1016/j.advwatres.2020.103798}
}

@ARTICLE{Beavers_Joseph_67,
  author = {Beavers, G.~S. and Joseph, D.~D.},
  title = {Boundary conditions at a naturally permeable wall},
  journal = {J. Fluid Mech.},
  year = {1967},
  volume = {30},
  doi = {10.1017/S0022112067001375},
  pages = {197--207}
}

@ARTICLE{Discacciati_Miglio_Quarteroni_02,
  author = {M. Discacciati and E. Miglio and A. Quarteroni},
  title = {Mathematical and numerical models for coupling surface and
  groundwater flows},
  journal = {Appl. Num. Math.},
  year = {2002},
  volume = {43},
  pages = {57--74},
  doi = {10.1016/S0168-9274(02)00125-3}
}

@ARTICLE{Discacciati_Quarteroni_09,
  author = {Discacciati, M. and Quarteroni, A.},
  title = {Navier--{S}tokes/{D}arcy coupling: modeling, analysis, and
  numerical approximation},
  journal = {Rev. Mat. Complut.},
  year = {2009},
  volume = {22},
  pages = {315-426},
  issue ={2},
  doi = {10.5209/rev_REMA.2009.v22.n2.16263}
}

@ARTICLE{Eggenweiler_Rybak_20,
   author = {Eggenweiler, E. and Rybak, I.},
   title = {Unsuitability of the {B}eavers--{J}oseph interface condition for filtration problems},
   journal = {J. Fluid Mech.},
   volume = {892},
   pages = {A10},
   year = {2020},
   doi = {10.1017/jfm.2020.194}
}

@article{Eggenweiler_Rybak_MMS20,
  author = {Eggenweiler, E. and Rybak, I.},
  journal = {Multiscale Model. Simul.},
  title = {Effective coupling conditions for arbitrary flows in {S}tokes--{D}arcy systems},
  doi = {10.1137/20M1346638},
  year = 2021,
  volume = {19},
  numberN = {2},
  pages ={731--757}
}

@ARTICLE{Goyeau_Lhuillier_etal_03,
  author = {Goyeau, B. and Lhuillier, D. and Gobin, D. and Velarde, M.},
  title = {Momentum transport at a fluid--porous interface},
  journal = {Int. J. Heat Mass Transfer},
  year = {2003},
  volume = {46},
  pages = {4071--4081},
  doi = {10.1016/S0017-9310(03)00241-2}
}

@ARTICLE{Jaeger_Mikelic_09,
  author = {J\"{a}ger, W. and Mikeli\'{c}, A.},
  title = {Modeling effective interface laws for transport phenomena
  between an unconfined fluid and a porous medium using
  homogenization},
  journal = {Transp. Porous Media},
  year = {2009},
  volume = {78},
  doi = {10.1007/s11242-009-9354-9},
  pages = {489--508}
}

@ARTICLE{Lacis_etal_20,
   author = {L\={a}cis, U. and Sudhakar, Y. and Pasche, S. and Bagheri, S.},
   title = {Transfer of mass and momentum at rough and porous surfaces},
   journal = {J. Fluid Mech.},
   volume = {884},
   pages = {A21}, 
   year = {2020},
   doi = {10.1017/jfm.2019.897}
}

@ARTICLE{Saffman,
  author = {Saffman, P.~G.},
  title = {On the boundary condition at the surface of a porous
  medium},
  journal = {Stud. Appl. Math.},
  year = {1971},
  volume = {50},
  numberN ={2},
  pages = {93--101},
  doi = {10.1002/sapm197150293}
}

@article{Sudhakar_21,
	author = {Y. Sudhakar and U. Lacis and S. Pasche and S. Bagheri},
	title = {Higher-Order Homogenized Boundary Conditions for Flows Over Rough and Porous Surfaces},
	journal = {Transp. Porous Media},
	year = {2021},
  volume = {136},
  doi = {10.1007/s11242-020-01495-w},
  pages = {1--42}
}

@ARTICLE{Lacis_Bagheri_17,
  author = {L\={a}cis, U. and Bagheri, S.},
  title = {A framework for computing effective boundary conditions at the
  interface between free fluid and a porous medium},
  journal = {J. Fluid Mech.},
  volume = {812},
  pages = {866--889},
  year = {2017},
  doi = {10.1017/jfm.2016.838}
}

@article{zampogna2020effective,
  title={Effective stress jump across membranes},
  author={Zampogna, G.~A. and Gallaire, F.},
  journal={J. Fluid Mech.},
  volume={892},
  pages={A9},
  year={2020},
  doi = {10.1017/jfm.2020.144}
}

@inproceedings{Ruan-Rybak-23,
  author = {Ruan, L. and Rybak, I.},
  booktitle = {Finite Volumes for Complex Applications X~–~Volume~1, Elliptic and Parabolic Problems},
  title = {{S}tokes--{B}rinkman--{D}arcy models for coupled free--flow and porous--medium systems},
  editor = {Franck, E. and et~al.},
  year = {2023},
  doi = {10.1007/978-3-031-40864-9\_31},
  series = {Springer Proceedings in Mathematics \& Statistics},
  volume = {432},
  pages= {365--373}
}

@article{Strohbeck-Eggenweiler-Rybak-23,
  author = {Strohbeck, P. and Eggenweiler, E. and Rybak, I.},
  journal = {Transp. Porous Media},
  title = {A modification of the {B}eavers--{J}oseph condition for arbitrary flows to the fluid--porous interface},
  volume = {147},
  pages = {605--628},
  year = {2023},
doi = {10.1007/s11242-023-01919-3}
}

@article{Lasseux_Valdes-Parada_Bottaro_2021, 
title={Upscaled model for unsteady slip flow in porous media}, 
volume={923}, 
pages = {A37},
doi={10.1017/jfm.2021.606}, 
journal={J. Fluid Mech.}, 
author={Lasseux, D. and Valdés-Parada, F. J. and Bottaro, A.}, 
year={2021}
}

@article{Ruan_Rybak_AMC,
  author = {Ruan, L. and Rybak, I.},
  journal = {Appl. Math. Comp. },
  title = {Stokes--{B}rinkman--{D}arcy models for fluid--porous systems: derivation, analysis and validation},
volume = {510},
pages = {129687},
year = {2026},
issn = {0096-3003},
doi = {10.1016/j.amc.2025.129687},
}

@article{Ruan_Rybak_TIPM, 
title={A hybrid-dimensional {S}tokes--{B}rinkman--{D}arcy model for arbitrary flows to the fluid--porous interface}, 
author={Ruan, L. and Rybak, I.}, 
journal={Transp. Porous Media}, 
volume={152}, 
year={2025},
issn = {1573-1634},
pages = {75},
doi={10.1007/s11242-025-02220-1}
}

@BOOK{vafai2005handbook,
  title={Handbook of Porous Media},
  author={Vafai, K.},
  url={https://books.google.de/books?id=5JNjIM2pRzUC},
  year={2005},
  publisher={CRC Press},
}

@article{OCHOATAPIA19952647,
title = {Momentum transfer at the boundary between a porous medium and a homogeneous fluid--{II}.~{C}omparison with experiment},
journal = {Int. J. Heat Mass Transfer},
volume = {38},
pages = {2647-2655},
year = {1995},
issn = {0017-9310},
doi = {10.1016/0017-9310(94)00347-X},
author = {J.Alberto Ochoa-Tapia and Stephen Whitaker},
}

@article{OCHOATAPIA19952635,
title = {Momentum transfer at the boundary between a porous medium and a homogeneous fluid--{I}.~{T}heoretical development},
journal = {Int. J. Heat Mass Transfer},
volume = {38},
pages = {2635-2646},
year = {1995},
issn = {0017-9310},
doi = {10.1016/0017-9310(94)00346-W},
author = {J.Alberto Ochoa-Tapia and Stephen Whitaker},
}

@article{RICCIARDI2005248,
title = {Comparison of the lognormal and beta distribution functions to describe the uncertainty in permeability},
journal = {J. Hydrology},
volume = {313},
number = {3},
pages = {248-256},
year = {2005},
issn = {0022-1694},
doi = {10.1016/j.jhydrol.2005.03.007},
author = {Karen L. Ricciardi and George F. Pinder and Kenneth Belitz},
}

@article{kroker2023global,
  author = {Kr{\"o}ker, Ilja and Oladyshkin, Sergey and Rybak, Iryna},
  doi = {10.1007/s10596-023-10236-z},
  journal = {Comput. Geosci.},
  volume = 27, 
  pages = {805--827},
  title = {Global sensitivity analysis using multi-resolution polynomial chaos expansion for coupled {S}tokes--{D}arcy flow problems},
  year = 2023
}

@article{kroker2022arbitrary,
  author = {Kr{\"o}ker, Ilja and Oladyshkin, Sergey},
  doi = {10.1016/j.ress.2022.108376},
  journal = {Reliab. Eng. Syst. Saf.},
  fjournal = {Reliability Engineering \& System Safety},
  pages = 108376,
  title = {Arbitrary Multi-Resolution Multi-Wavelet-based Polynomial Chaos Expansion for Data-Driven Uncertainty Quantification},
  volume = 222,
  year = 2022
}

@article{kroeker2025bayesian3,
  author = {Kröker, Ilja and Brünnette, Tim and Wildt, Nils and Oreamuno, Maria Fernanda Morales and Kohlhaas, Rebecca and Oladyshkin, Sergey and Nowak, Wolfgang},
  doi = {10.1615/Int.J.UncertaintyQuantification.2024052675},
  issn = {2152-5080},
  journal = {Int. J. Uncertain. Quantif.},
  number = 3,
  pages = {21--54},
  title = {Bayesian$^3$ active learning for regularized multi-resolution arbitrary polynomial chaos using information theory},
  volume = 15,
  year = 2025
}

@article{kohlhaas2023gaussian,
  author = {Kohlhaas, Rebecca and Kr{\"o}ker, Ilja and Oladyshkin, Sergey and Nowak, Wolfgang},
  doi = {doi:10.1007/s10596-023-10199-1},
  journal = {Comput. Geosci.},
  number = 3,
  pages = {1--21},
  title = {Gaussian active learning on multi-resolution arbitrary polynomial chaos emulator: concept for bias correction, assessment of surrogate reliability and its application to the carbon dioxide benchmark},
  volume = 27,
  year = 2023
}

@article{Buerger2012,
  author = {B{\"u}rger, Raimund and Kr{\"o}ker, Ilja and Rohde, Christian},
  doi = {10.1002/zamm.201200174},
  journal = {ZAMM Z. Angew. Math. Mech.},
  number = 10,
  pages = {793-817},
  title = {A hybrid stochastic {G}alerkin method for uncertainty quantification
 applied to a conservation law modelling a clarifier-thickener unit},
  url = {http://dx.doi.org/10.1002/zamm.201200174},
  volume = 94,
  year = 2014
}

@article{Koeppel2017,
  author = {K{\"o}ppel, M. and Kr{\"o}ker, I. and Rohde, C.},
  issn = {1573-1499},
  journal = {Comput. Geosci.},
  number = 4,
  pages = {807--832},
  title = {{Intrusive uncertainty quantification for hyperbolic-elliptic systems governing two-phase flow in heterogeneous porous media}},
  volume = 21,
  year = 2017
}

@article{buerkner2022,
title = {A fully {B}ayesian sparse polynomial chaos expansion approach with joint priors on the coefficients and global selection of terms},
journal = {J. Comput. Phys.},
pages = {112210},
year = {2023},
issn = {0021-9991},
doi = {https://doi.org/10.1016/j.jcp.2023.112210},
url = {https://www.sciencedirect.com/science/article/pii/S0021999123003054},
author = {Paul-Christian Bürkner and Ilja Kröker and Sergey Oladyshkin and Wolfgang Nowak}
}

@article{MR2855645,
	title        = {On the convergence of generalized polynomial chaos expansions},
	author       = {Ernst, Oliver G. and Mugler, Antje and Starkloff, Hans-J\"org and Ullmann, Elisabeth},
	year         = 2012,
	journal      = {ESAIM Math. Model. Numer. Anal.},
	volume       = 46,
	number       = 2,
	pages        = {317--339},
	doi          = {10.1051/m2an/2011045},
	issn         = {0764-583X},
	url          = {https://doi.org/10.1051/m2an/2011045},
	fjournal     = {ESAIM. Mathematical Modelling and Numerical Analysis},
	mrclass      = {60H35 (33C45 35R60 60G07 65N75)},
	mrnumber     = 2855645,
	mrreviewer   = {Renato G. C. Spigler}
}

@article{favard1935polynomes,
	title        = {Sur les polynomes de {T}chebicheff},
	author       = {Favard, Jean},
	year         = 1935,
	journal      = {CR Acad. Sci. Paris},
	volume       = 200,
	number       = {2052-2055},
	pages        = 11
}

@article{Luethen2021,
	title        = {Sparse Polynomial Chaos Expansions: Literature Survey and Benchmark},
	author       = {L\"{u}then, Nora and Marelli, Stefano and Sudret, Bruno},
	year         = 2021,
	journal      = {SIAM/ASA J. Uncertain. Quantif.},
	volume       = 9,
	number       = 2,
	pages        = {593--649},
	doi          = {10.1137/20M1315774},
	eprint       = {https://doi.org/10.1137/20M1315774},
}

@book{Ghanem91StochFE,
	title        = {Stochastic Finite Elements: a Spectral Approach},
	author       = {R. G. Ghanem and P. D. Spanos},
	year         = 1991,
	publisher    = {Springer-Verlag},
	address      = {New York},
	pages        = {x+214},
	isbn         = {0-387-97456-3},
	mrclass      = {73V05 (65U05 73K35 73V30)},
	mrnumber     = {1083354 (91k:73102)}
}

@article{wiener1938,
	title        = {The homogeneous chaos},
	author       = {Wiener, Norbert},
	year         = 1938,
	journal      = {Amer. J. Math.},
	publisher    = {JSTOR},
	volume       = 60,
	number       = 4,
	pages        = {897--936}
}

@article{SUDRET2008964,
title = "Global sensitivity analysis using polynomial chaos expansions",
journal = "Reliab. Eng. Syst. Saf.",
volume = "93",
number = "7",
pages = "964--979",
year = "2008",
noteN = "Bayesian Networks in Dependability",
issn = "0951-8320",
doi = "10.1016/j.ress.2007.04.002",
author = "Bruno Sudret"
}

@article{HOMMA19961,
    author = "Toshimitsu Homma and Andrea Saltelli",
	doi = "10.1016/0951-8320(96)00002-6",
	journal = "Reliab. Eng. Syst. Saf.",
	number = "1",
	pages = "1--17",
	title = "Importance measures in global sensitivity analysis of nonlinear models",
	volume = "52",
	year = "1996"
}

@article {MR1052836,
    AUTHOR = {Sobol, I. M.},
     TITLE = {Estimation of the sensitivity of nonlinear mathematical models},
   JOURNAL = {Mat. Model.},
  FJOURNAL = {Matematicheskoe Modelirovanie},
    VOLUME = {2},
      YEAR = {1990},
    NUMBER = {1},
     PAGES = {112--118},
      ISSN = {0234-0879},
   MRCLASS = {65C05 (65D30)},
  MRNUMBER = {1052836},
MRREVIEWER = {I. V\u{a}duva}
}

@ARTICLE{MR1823119,
    AUTHOR = {Sobol, I. M.},
     TITLE = {Global sensitivity indices for nonlinear mathematical models
              and their {M}onte {C}arlo estimates},
      NOTEN = {The Second IMACS Seminar on Monte Carlo Methods (Varna, 1999)},
   JOURNAL = {Math. Comput. Simulation},
  FJOURNAL = {Mathematics and Computers in Simulation},
    VOLUME = {55},
      YEAR = {2001},
    NUMBERN = {1-3},
     PAGES = {271--280},
      ISSN = {0378-4754},
   MRCLASS = {65C05 (65D30)},
  MRNUMBER = {1823119},
       DOI = {10.1016/S0378-4754(00)00270-6}
}

@article{OLADYSHKIN2012179,
title = "Data-driven uncertainty quantification using the arbitrary polynomial chaos expansion",
journal = "Reliab. Eng. Syst. Saf.",
volume = "106",
pages = "179 - 190",
year = "2012",
doi = "10.1016/j.ress.2012.05.002",
author = "S. Oladyshkin and W. Nowak"
}

@ARTICLE{MR2063905,
  author = {O. P. {Le Ma{\^{\i}}tre} and H. N. Najm and R. G. Ghanem and O. M.
	Knio},
  title = {Multi-resolution analysis of {W}iener-type uncertainty propagation
	schemes},
  journal = {J. Comput. Phys.},
  year = {2004},
  volume = {197},
  pages = {502--531},
  number = {2},
  coden = {JCTPAH},
  doi = {10.1016/j.jcp.2003.12.020},
  fjournal = {Journal of Computational Physics},
  issn = {0021-9991},
  mrclass = {65P99 (42C40)},
  mrnumber = {2063905 (2005b:65142)},
  urlN = {http://dx.doi.org/10.1016/j.jcp.2003.12.020}
}

@book {MR2061539,
    AUTHOR = {Gautschi, Walter},
     TITLE = {Orthogonal Polynomials: Computation and Approximation},
    SERIESN = {Numerical Mathematics and Scientific Computation},
      NOTEN = {Oxford Science Publications},
 PUBLISHER = {Oxford University Press, New York},
      YEAR = {2004},
     PAGES = {x+301},
      ISBN = {0-19-850672-4},
   MRCLASS = {42-02 (33C45 33F05 42C05)},
  MRNUMBER = {2061539},
MRREVIEWER = {Leonid B. Golinski\u\i },
}

@article{Marelli_2021,
	author  = {Stefano Marelli and Paul-Remo Wagner and Christos Lataniotis and Bruno Sudret},
	title   = {Stochastic spectral embedding},
	journal = {Int. J. Uncertain. Quantif.},
	issn    = {2152-5080},
	year    = {2021},
	volume  = {11},
	number  = {2},
	pages   = {25--47},
	doi = {10.1615/Int.J.UncertaintyQuantification.2020034395}
}

@article {MR2240796,
    AUTHOR = {Wan, Xiaoliang and Karniadakis, George Em},
     TITLE = {Multi-element generalized polynomial chaos for arbitrary
              probability measures},
   JOURNAL = {SIAM J. Sci. Comput.},
  FJOURNAL = {SIAM Journal on Scientific Computing},
    VOLUME = {28},
      YEAR = {2006},
    NUMBER = {3},
     PAGES = {901--928},
      ISSN = {1064-8275},
   MRCLASS = {65C30},
  MRNUMBER = {2240796},
MRREVIEWER = {Vassil S. Grozdanov},
       DOI = {10.1137/050627630},
       URLN = {https://doi.org/10.1137/050627630},
}

@BOOK{sullivan_book,
    AUTHOR = {Sullivan, T. J.},
     TITLE = {Introduction to Uncertainty Quantification},
    SERIESN = {Texts in Applied Mathematics},
    VOLUMEN = {63},
 PUBLISHER = {Springer},
      YEAR = {2015},
     PAGES = {xii+342},
      ISBN = {978-3-319-23394-9; 978-3-319-23395-6},
   MRCLASS = {60H30 (41-01 62F15 62Jxx 65-02 81P45 94A15)},
  MRNUMBER = {3364576},
MRREVIEWER = {Elisabeth Ullmann},
       DOI = {10.1007/978-3-319-23395-6},
       URLN = {http://dx.doi.org/10.1007/978-3-319-23395-6},
}

@article{sobol1967,
  title={On the distribution of points in a cube and the approximate evaluation of integrals},
  author={Sobol', Il'ya Meerovich},
  journal={Zhurnal Vychislitel'noi Matematiki i Matematicheskoi Fiziki},
  volume={7},
  number={4},
  pages={784--802},
  year={1967},
  publisher={Russian Academy of Sciences, Branch of Mathematical Sciences}
}

@article{WAGNER2022102179,
title = {Rare event estimation using stochastic spectral embedding},
journal = {Structural Safety},
volume = {96},
pages = {102179},
year = {2022},
issn = {0167-4730},
doi = {https://doi.org/10.1016/j.strusafe.2021.102179},
url = {https://www.sciencedirect.com/science/article/pii/S0167473021000990},
author = {P.-R. Wagner and S. Marelli and I. Papaioannou and D. Straub and B. Sudret}
}

@article{Kang_Wang_2024, 
title={Brinkman double-layer model for flow at a free-porous interface}, 
author={Kang, J. and Wang, M.}, 
journal={Int. J. Mech. Sci.}, 
volume={263}, 
year={2024},
pages= {108770},
doi={10.1016/j.ijmecsci.2023.108770}
}

@data{DARUS-5251_2025,
author = {Ruan, Linheng and Kröker, Ilja and  Rybak, Iryna and Oladyshkin, Sergey},
publisher = {DaRUS},
title = {Replication Data for: Surrogate-assisted Global sensitivity analysis for hybrid-dimensional {S}tokes-{B}rinkman-{D}arcy model},
year = {2025},
version = {DRAFT VERSION},
doi = {10.18419/DARUS-5251},
url = {https://doi.org/10.18419/DARUS-5251},
note = {doi: 10.18419/DARUS-5251},
}

@book {daveiga_book_2021,
    AUTHOR = {Da Veiga, S\'ebastien and Gamboa, Fabrice and Iooss, Bertrand
              and Prieur, Cl\'ementine},
     TITLE = {Basics and trends in sensitivity analysis--theory and
              practice in {R}},
    SERIES = {Computational Science \& Engineering},
    VOLUME = {23},
 PUBLISHER = {Society for Industrial and Applied Mathematics (SIAM),
              Philadelphia, PA},
      YEAR = {2021},
     PAGES = {xvi+291},
      ISBN = {978-1-611976-68-7},
   MRCLASS = {62-01 (62-04 62-08)},
  MRNUMBER = {4359984},
       DOI = {10.1137/1.9781611976694},
       URL = {https://doi.org/10.1137/1.9781611976694},
}

@incollection {shapley1953,
    AUTHOR = {Shapley, L. S.},
     TITLE = {A value for {$n$}-person games},
 BOOKTITLE = {Contributions to the theory of games, vol. 2},
    SERIES = {Ann. of Math. Stud.},
    VOLUME = {no. 28},
     PAGES = {307--317},
 PUBLISHER = {Princeton Univ. Press, Princeton, NJ},
      YEAR = {1953},
   MRCLASS = {90.0X},
  MRNUMBER = {53477},
MRREVIEWER = {D.\ Gale},
}

@article {owen2014,
    AUTHOR = {Owen, Art B.},
     TITLE = {Sobol' indices and {S}hapley value},
   JOURNAL = {SIAM/ASA J. Uncertain. Quantif.},
  FJOURNAL = {SIAM/ASA Journal on Uncertainty Quantification},
    VOLUME = {2},
      YEAR = {2014},
    NUMBER = {1},
     PAGES = {245--251},
      ISSN = {2166-2525},
   MRCLASS = {91B82 (60G25 91A12)},
  MRNUMBER = {3283908},
MRREVIEWER = {R\u azvan\ R\u aducanu},
       DOI = {10.1137/130936233},
       URL = {https://doi.org/10.1137/130936233},
}

@article{fisher1923manurial,
  title={The manurial response of different potato varieties},
  author={Fisher, RA and Mackenzie, WA},
  journal={J. Agric. Sci.},
  volume={13},
  pages={311--320},
  year={1923}
}

@article{Fisher_1921,
title={Studies in crop variation. {I.}~{A}n examination of the yield of dressed grain from {B}roadbalk}, volume={11},
DOI={10.1017/S0021859600003750},
number={2},
journal={J. Agric. Sci.},
author={Fisher, R. A.},
year={1921},
pages={107–135}
}

@article {MR783216,
    AUTHOR = {Askey, Richard and Wilson, James},
     TITLE = {Some basic hypergeometric orthogonal polynomials that
              generalize {J}acobi polynomials},
   JOURNAL = {Mem. Amer. Math. Soc.},
  FJOURNAL = {Memoirs of the American Mathematical Society},
    VOLUME = {54},
      YEAR = {1985},
    NUMBER = {319},
     PAGES = {iv+55},
      ISSN = {0065-9266},
   MRCLASS = {05A30 (33A65)},
  MRNUMBER = {783216},
MRREVIEWER = {Mourad E. H. Ismail},
       DOI = {10.1090/memo/0319},
       URL = {https://doi.org/10.1090/memo/0319},
}

@article {MR1951058,
    AUTHOR = {Xiu, Dongbin and Karniadakis, George Em},
     TITLE = {The {W}iener-{A}skey polynomial chaos for stochastic
              differential equations},
   JOURNAL = {SIAM J. Sci. Comput.},
  FJOURNAL = {SIAM Journal on Scientific Computing},
    VOLUME = {24},
      YEAR = {2002},
    NUMBER = {2},
     PAGES = {619--644},
      ISSN = {1064-8275},
   MRCLASS = {60H15 (33C45 33C90 65C30)},
  MRNUMBER = {1951058},
       DOI = {10.1137/S1064827501387826},
       URL = {https://doi.org/10.1137/S1064827501387826},
}

@article {MR2121216,
    AUTHOR = {Matthies, Hermann G. and Keese, Andreas},
     TITLE = {Galerkin methods for linear and nonlinear elliptic stochastic
              partial differential equations},
   JOURNAL = {Comput. Methods Appl. Mech. Engrg.},
  FJOURNAL = {Computer Methods in Applied Mechanics and Engineering},
    VOLUME = {194},
      YEAR = {2005},
    NUMBER = {12-16},
     PAGES = {1295--1331},
      ISSN = {0045-7825},
   MRCLASS = {65N30 (60H15 60H35)},
  MRNUMBER = {2121216},
       DOI = {10.1016/j.cma.2004.05.027},
       URL = {https://doi.org/10.1016/j.cma.2004.05.027},
}

@ARTICLE{Xiu2003137,
  author = {D. Xiu and G. E. Karniadakis},
  title = {Modeling uncertainty in flow simulations via generalized polynomial
	chaos},
  journal = {J. Comput. Phys.},
  year = {2003},
  volume = {187},
  pages = {137--167},
  number = {1},
  doi = {10.1016/S0021-9991(03)00092-5},
  issn = {0021-9991},
  keywords1 = {Uncertainty},
  keywords2 = {Fluids},
  keywords3 = {Stochastic modeling}
}

@article{Pettersson2016367,
title = {Stochastic {G}alerkin framework with locally reduced bases for nonlinear two-phase transport in heterogeneous formations},
journal = {Comput. Methods Appl. Mech. Engrg.},
volume = 310,
pages = {367--387},
year = 2016,
doi = {10.1016/j.cma.2016.07.013},
author = {Per Pettersson and Hamdi A. Tchelepi}
}

@ARTICLE{MR2501693,
  author = {G. Po{\"e}tte and B. Despr{\'e}s and D. Lucor},
  title = {Uncertainty quantification for systems of conservation laws},
  journal = {J. Comput. Phys.},
  year = {2009},
  volume = {228},
  pages = {2443--2467},
  number = {7},
  coden = {JCTPAH},
  doi = {10.1016/j.jcp.2008.12.018},
  fjournal = {Journal of Computational Physics},
  issn = {0021-9991},
  mrclass = {65C30 (35L65)},
  mrnumber = {2501693 (2010c:65014)},
  mrreviewer = {Carlos M. Mora Gonz{\'a}lez},
  url = {http://dx.doi.org/10.1016/j.jcp.2008.12.018}
}

@article{KUSCH2020109073,
title = {Filtered stochastic {G}alerkin methods for hyperbolic equations},
journal = {J. Comput. Phys.},
fjournal = {Journal of Computational Physics},
volume = {403},
pages = {109073},
year = {2020},
issn = {0021-9991},
doi = {10.1016/j.jcp.2019.109073},
author = {Jonas Kusch and Ryan G. McClarren and Martin Frank}
}

@article{Schobi_2015,
	author  = {Roland Schobi and Bruno Sudret and Joe  Wiart},
	title   = {POLYNOMIAL-CHAOS-BASED KRIGING},
	journal = {Int. J. Uncertain. Quantif.},
	issn    = {2152-5080},
	year    = {2015},
	volume  = {5},
	number  = {2},
	pages   = {171--193}
}

@article{BLATMAN20101216,
author = {Géraud Blatman and Bruno Sudret},
title = {Efficient computation of global sensitivity indices using sparse polynomial chaos expansions},
journal = {Reliab. Eng. Syst. Saf.	},
volume = {95},
number = {11},
pages = {1216-1229},
year = {2010},
issn = {0951-8320},
doi = {https://doi.org/10.1016/j.ress.2010.06.015},
url = {https://www.sciencedirect.com/science/article/pii/S0951832010001493},
}

@article {MR2764550,
    AUTHOR = {Blatman, G\'eraud and Sudret, Bruno},
     TITLE = {Adaptive sparse polynomial chaos expansion based on least
              angle regression},
   JOURNAL = {J. Comput. Phys.},
  FJOURNAL = {Journal of Computational Physics},
    VOLUME = {230},
      YEAR = {2011},
    NUMBER = {6},
     PAGES = {2345--2367},
      ISSN = {0021-9991,1090-2716},
   MRCLASS = {65N30 (60G35 62G08 65N75)},
  MRNUMBER = {2764550},
       DOI = {10.1016/j.jcp.2010.12.021},
       URL = {https://doi.org/10.1016/j.jcp.2010.12.021},
}

@inproceedings{penrose1956best,
  title={On best approximate solutions of linear matrix equations},
  author={Penrose, Roger},
  booktitle={Mathematical Proceedings of the Cambridge Philosophical Society},
  volume={52},
  cnumber={1},
  pages={17--19},
  year={1956},
  organization={Cambridge University Press}
}

@article {MR3023712,
    AUTHOR = {Tryoen, J. and Le Ma\^itre, O. and Ern, A.},
     TITLE = {Adaptive anisotropic spectral stochastic methods for uncertain
              scalar conservation laws},
   JOURNAL = {SIAM J. Sci. Comput.},
  FJOURNAL = {SIAM Journal on Scientific Computing},
    VOLUME = {34},
      YEAR = {2012},
    NUMBER = {5},
     PAGES = {A2459--A2481},
      ISSN = {1064-8275,1095-7197},
   MRCLASS = {65M75 (35L65 35R60 60H35 65C20 68U20)},
  MRNUMBER = {3023712},
MRREVIEWER = {Guangming\ Yao},
       DOI = {10.1137/120863927},
       URL = {https://doi.org/10.1137/120863927},
}

@software{ikroeker_aMR-PC_2023,
  author       = {Ilja Kröker and Sergey Oladyshkin},
  title        = {{aMR-PC}: Arbitrary Multi-Resolution Polynomial Chaos Python Toolbox},
  url          = {https://github.com/ikroeker/aMR-PC},
  version      = {v1.0.2},
  date         = {2023-07-07},
  license      = {MIT},
  year         = {2023}, 
  note         = {https://github.com/ikroeker/aMR-PC},
  note         = {Python package for multi-resolution polynomial chaos expansion.}
}

@article{DREAU2023117920,
title = {Multi-element polynomial chaos expansion based on automatic discontinuity detection for nonlinear systems},
journal = {J. Sound Vib.},
volume = {567},
pages = {117920},
year = {2023},
issn = {0022-460X},
doi = {https://doi.org/10.1016/j.jsv.2023.117920},
url = {https://www.sciencedirect.com/science/article/pii/S0022460X23003693},
author = {Juliette Dréau and Benoit Magnain and Alain Batailly}
}

\appendix
\section{Supplementary Tables and Figures}\label{appendixA}
In this section, we provide additional tabulated data supporting the results of GSA presented in Sections~\ref{sec:totalSobolindicies_i} and \ref{sec:totalSobolindicies_ii}. In particular, Tab.~\ref{tab:Sobol-total} reports the total sensitivity indices computed with \aMR{} for different choices of polynomial degree~$\No$, refinement level~$\Nr$, and sparsity by hyperbolic truncation parameter~$q$ for the filtration problem (Case~1). The corresponding results for the splitting flow problem (Case~2) are summarized in Tab.~\ref{tab:Sobol-total_ii}. 
Furthermore, Tabs.~\ref{tab:Sobol-u} and \ref{tab:Sobol-v} present the space-averaged Sobol’ indices of the tangential velocity~$u$ and normal velocity~$v$, respectively, for Case~1, while Tabs.~\ref{tab:Sobol-u_ii} and \ref{tab:Sobol-v_ii} report the corresponding results for Case~2.

\begin{table}[h!]
    \caption{Case 1: Averaged total sensitivity indices.}
    \label{tab:Sobol-total}
    \centering

\begin{subtable}{0.99\linewidth}
    \caption{Horizontal velocity component $u$.}
    \begin{tabular}{|c|c|c|c|c|c|c|}
\hline
 & $N_o=2$ & $N_o=4$ & $N_o=4 / q=0.75$ & $N_r=1 / N_o=2$ & $N_r=1 / N_o=4$ & $N_r=1 / N_o=4 / q=0.75$ \\
\hline
$I=1$ & 5.18e-01 & 1.93e-01 & 3.00e-02 & 5.05e-01 & 2.32e-01 & 6.91e-02 \\
$I=2$ & 1.85e-01 & 7.46e-01 & 1.41e-02 & 2.70e-01 & 9.57e-01 & 3.86e-02 \\
$I=3$ & 1.74e-01 & 1.13e-01 & 6.96e-02 & 1.97e-01 & 1.76e-01 & 3.20e-02 \\
$I=4$ & 1.18e-02 & 1.29e-02 & 8.77e-02 & 1.78e-02 & 1.04e-01 & 2.29e-02 \\
$I=5$ & 1.91e-01 & 3.14e-02 & 6.25e-02 & 1.78e-01 & 1.59e-01 & 2.55e-02 \\
\hline
\end{tabular}

    \end{subtable}\\[4mm]

    \begin{subtable}{0.99\linewidth}
    \caption{Vertical velocity component $v$.}
    \begin{tabular}{|c|c|c|c|c|c|c|}
\hline
 & $N_o=2$ & $N_o=4$ & $N_o=4 / q=0.75$ & $N_r=1 / N_o=2$ & $N_r=1 / N_o=4$ & $N_r=1 / N_o=4 / q=0.75$ \\
\hline
$I=1$ & 4.87e-01 & 1.94e-01 & 2.80e-02 & 4.77e-01 & 2.42e-01 & 6.75e-02 \\
$I=2$ & 1.77e-01 & 7.29e-01 & 1.37e-02 & 2.62e-01 & 9.58e-01 & 3.86e-02 \\
$I=3$ & 1.63e-01 & 1.18e-01 & 6.65e-02 & 1.82e-01 & 1.72e-01 & 3.07e-02 \\
$I=4$ & 1.19e-02 & 1.35e-02 & 8.45e-02 & 1.77e-02 & 1.03e-01 & 2.26e-02 \\
$I=5$ & 1.87e-01 & 3.10e-02 & 6.20e-02 & 1.74e-01 & 1.61e-01 & 2.44e-02 \\
\hline
\end{tabular}

    \end{subtable}
\end{table}

\begin{table}
    \caption{Case 1: Averaged Sobol' sensitivity indices.}
    \label{tab:Sobol-u}
\begin{subtable}{0.99\linewidth}
    \caption{Horizontal velocity component $u$.}
    \centering
    \begin{tabular}{|l|c|c|c|c|c|c|}
\hline
 & $N_o=2$ & $N_o=4$ & $N_o=4 / q=0.75$ & $N_r=1 / N_o=2$ & $N_r=1 / N_o=4$ & $N_r=1 / N_o=4 / q=0.75$ \\
\hline
$I=1$ & 4.17e-01 & 1.06e-01 & 1.93e-02 & 3.53e-01 & 6.47e-03 & 2.54e-02 \\
$I=12$ & 1.91e-02 & 3.99e-02 & 4.45e-04 & 3.49e-02 & 9.37e-02 & 1.62e-02 \\
$I=123$ & 0.00e+00 & 1.29e-03 & 0.00e+00 & 4.40e-03 & 2.35e-02 & 4.57e-04 \\
$I=1234$ & 0.00e+00 & 9.11e-05 & 0.00e+00 & 3.97e-04 & 4.03e-03 & 1.55e-03 \\
$I=12345$ & 0.00e+00 & 0.00e+00 & 0.00e+00 & 7.28e-04 & 4.85e-02 & 1.28e-03 \\
$I=1235$ & 0.00e+00 & 1.94e-05 & 0.00e+00 & 2.35e-04 & 1.47e-02 & 1.55e-03 \\
$I=124$ & 0.00e+00 & 7.79e-04 & 0.00e+00 & 7.13e-04 & 6.29e-03 & 2.34e-03 \\
$I=1245$ & 0.00e+00 & 8.92e-07 & 0.00e+00 & 1.69e-04 & 2.11e-02 & 9.65e-04 \\
$I=125$ & 0.00e+00 & 7.66e-04 & 0.00e+00 & 3.92e-03 & 1.07e-02 & 2.85e-03 \\
$I=13$ & 5.02e-02 & 3.92e-02 & 0.00e+00 & 7.77e-02 & 2.41e-03 & 3.74e-03 \\
$I=134$ & 0.00e+00 & 1.18e-04 & 0.00e+00 & 2.99e-04 & 1.48e-05 & 3.59e-03 \\
$I=1345$ & 0.00e+00 & 5.32e-06 & 0.00e+00 & 7.69e-05 & 1.01e-05 & 1.29e-03 \\
$I=135$ & 0.00e+00 & 2.29e-05 & 0.00e+00 & 6.05e-05 & 4.05e-06 & 1.59e-03 \\
$I=14$ & 2.18e-03 & 2.11e-03 & 8.69e-04 & 4.34e-03 & 1.59e-04 & 2.81e-03 \\
$I=145$ & 0.00e+00 & 5.05e-05 & 0.00e+00 & 1.27e-04 & 8.00e-06 & 1.02e-03 \\
$I=15$ & 2.90e-02 & 3.34e-03 & 9.36e-03 & 2.42e-02 & 2.63e-04 & 2.48e-03 \\
$I=2$ & 1.47e-01 & 6.70e-01 & 7.58e-04 & 1.87e-01 & 6.03e-01 & 8.70e-03 \\
$I=23$ & 1.22e-02 & 2.29e-02 & 1.35e-05 & 1.59e-02 & 6.01e-02 & 8.82e-04 \\
$I=234$ & 0.00e+00 & 2.12e-04 & 0.00e+00 & 2.16e-04 & 5.70e-03 & 6.39e-04 \\
$I=2345$ & 0.00e+00 & 1.17e-06 & 0.00e+00 & 2.33e-04 & 2.07e-03 & 1.13e-05 \\
$I=235$ & 0.00e+00 & 2.33e-04 & 0.00e+00 & 1.89e-04 & 1.23e-02 & 7.66e-06 \\
$I=24$ & 1.21e-03 & 5.10e-03 & 1.29e-02 & 1.60e-03 & 6.37e-03 & 4.87e-04 \\
$I=245$ & 0.00e+00 & 1.13e-04 & 0.00e+00 & 1.48e-04 & 9.38e-03 & 4.38e-05 \\
$I=25$ & 5.66e-03 & 4.23e-03 & 0.00e+00 & 1.90e-02 & 3.61e-02 & 6.19e-04 \\
$I=3$ & 1.12e-01 & 4.84e-02 & 8.73e-04 & 9.60e-02 & 2.58e-03 & 9.96e-03 \\
$I=34$ & 2.59e-05 & 1.99e-05 & 6.82e-02 & 8.20e-05 & 3.19e-06 & 5.47e-03 \\
$I=345$ & 0.00e+00 & 3.75e-06 & 0.00e+00 & 4.73e-05 & 2.20e-06 & 1.45e-05 \\
$I=35$ & 3.89e-05 & 2.07e-05 & 5.51e-04 & 4.19e-05 & 2.05e-06 & 1.46e-05 \\
$I=4$ & 8.25e-03 & 4.21e-03 & 3.95e-03 & 8.49e-03 & 2.69e-04 & 1.21e-03 \\
$I=45$ & 1.24e-04 & 7.41e-05 & 1.88e-03 & 1.55e-04 & 7.36e-06 & 2.07e-04 \\
$I=5$ & 1.57e-01 & 2.25e-02 & 5.07e-02 & 1.29e-01 & 4.27e-03 & 1.16e-02 \\
\hline
\end{tabular}

    \end{subtable}
\end{table}

\begin{table}
    \caption{Case 1: Averaged Sobol' sensitivity indices.}
    \label{tab:Sobol-v}
    \begin{subtable}{0.99\linewidth}
    \caption{Vertical velocity component $v$.}
    \centering
    \begin{tabular}{|l|c|c|c|c|c|c|}
\hline
 & $N_o=2$ & $N_o=4$ & $N_o=4 / q=0.75$ & $N_r=1 / N_o=2$ & $N_r=1 / N_o=4$ & $N_r=1 / N_o=4 / q=0.75$ \\
\hline
$I=1$ & 3.89e-01 & 9.73e-02 & 1.74e-02 & 3.30e-01 & 6.09e-03 & 2.31e-02 \\
$I=12$ & 2.00e-02 & 4.73e-02 & 4.59e-04 & 3.62e-02 & 1.07e-01 & 1.66e-02 \\
$I=123$ & 0.00e+00 & 1.33e-03 & 0.00e+00 & 4.13e-03 & 2.30e-02 & 4.68e-04 \\
$I=1234$ & 0.00e+00 & 1.01e-04 & 0.00e+00 & 3.71e-04 & 3.63e-03 & 1.57e-03 \\
$I=12345$ & 0.00e+00 & 0.00e+00 & 0.00e+00 & 6.48e-04 & 4.67e-02 & 1.22e-03 \\
$I=1235$ & 0.00e+00 & 2.36e-05 & 0.00e+00 & 2.12e-04 & 1.38e-02 & 1.56e-03 \\
$I=124$ & 0.00e+00 & 8.52e-04 & 0.00e+00 & 6.86e-04 & 5.33e-03 & 2.44e-03 \\
$I=1245$ & 0.00e+00 & 1.04e-06 & 0.00e+00 & 1.45e-04 & 2.16e-02 & 9.20e-04 \\
$I=125$ & 0.00e+00 & 8.45e-04 & 0.00e+00 & 3.95e-03 & 1.20e-02 & 3.09e-03 \\
$I=13$ & 4.71e-02 & 4.04e-02 & 0.00e+00 & 7.23e-02 & 2.45e-03 & 3.60e-03 \\
$I=134$ & 0.00e+00 & 1.31e-04 & 0.00e+00 & 3.04e-04 & 1.58e-05 & 3.52e-03 \\
$I=1345$ & 0.00e+00 & 4.83e-06 & 0.00e+00 & 6.89e-05 & 9.97e-06 & 1.22e-03 \\
$I=135$ & 0.00e+00 & 2.18e-05 & 0.00e+00 & 5.06e-05 & 3.94e-06 & 1.60e-03 \\
$I=14$ & 2.18e-03 & 2.24e-03 & 8.74e-04 & 4.31e-03 & 1.69e-04 & 2.93e-03 \\
$I=145$ & 0.00e+00 & 5.39e-05 & 0.00e+00 & 1.26e-04 & 8.49e-06 & 9.75e-04 \\
$I=15$ & 2.87e-02 & 3.54e-03 & 9.26e-03 & 2.41e-02 & 2.86e-04 & 2.63e-03 \\
$I=2$ & 1.39e-01 & 6.43e-01 & 7.49e-04 & 1.80e-01 & 5.89e-01 & 8.05e-03 \\
$I=23$ & 1.15e-02 & 2.58e-02 & 1.03e-05 & 1.48e-02 & 5.83e-02 & 8.38e-04 \\
$I=234$ & 0.00e+00 & 2.37e-04 & 0.00e+00 & 2.05e-04 & 5.31e-03 & 6.21e-04 \\
$I=2345$ & 0.00e+00 & 1.03e-06 & 0.00e+00 & 1.97e-04 & 2.07e-03 & 1.12e-05 \\
$I=235$ & 0.00e+00 & 2.43e-04 & 0.00e+00 & 1.56e-04 & 1.38e-02 & 7.40e-06 \\
$I=24$ & 1.15e-03 & 5.02e-03 & 1.25e-02 & 1.61e-03 & 6.50e-03 & 4.84e-04 \\
$I=245$ & 0.00e+00 & 1.16e-04 & 0.00e+00 & 1.32e-04 & 1.09e-02 & 4.59e-05 \\
$I=25$ & 5.56e-03 & 4.00e-03 & 0.00e+00 & 1.89e-02 & 3.81e-02 & 6.62e-04 \\
$I=3$ & 1.04e-01 & 4.96e-02 & 9.22e-04 & 8.80e-02 & 2.59e-03 & 9.30e-03 \\
$I=34$ & 3.61e-05 & 1.92e-05 & 6.50e-02 & 7.92e-05 & 3.09e-06 & 5.16e-03 \\
$I=345$ & 0.00e+00 & 4.28e-06 & 0.00e+00 & 4.42e-05 & 2.26e-06 & 1.42e-05 \\
$I=35$ & 1.01e-05 & 1.30e-05 & 5.85e-04 & 1.60e-05 & 1.68e-06 & 1.25e-05 \\
$I=4$ & 8.45e-03 & 4.64e-03 & 4.11e-03 & 8.60e-03 & 2.95e-04 & 1.24e-03 \\
$I=45$ & 1.23e-04 & 7.98e-05 & 1.97e-03 & 1.49e-04 & 7.89e-06 & 2.17e-04 \\
$I=5$ & 1.52e-01 & 2.21e-02 & 5.02e-02 & 1.25e-01 & 2.18e-03 & 1.02e-02 \\
\hline
\end{tabular}

    \end{subtable}
\end{table}

\begin{table}
    \caption{Case 2: Averaged total sensitivity indices.}
    \label{tab:Sobol-total_ii}
    \centering

\begin{subtable}{0.99\linewidth}
    \caption{Horizontal velocity component $u$.}
    \centering
    \begin{tabular}{|c|c|c|c|c|}
\hline
 & $N_o=2$ & $N_o=4$ & $N_o=4 / q=0.75$ & $N_r=1 / N_o=2$ \\
\hline
$I=1$ & 3.01e-01 & 1.55e-01 & 2.83e-02 & 3.09e-01 \\
$I=2$ & 7.35e-02 & 3.01e-01 & 1.37e-02 & 8.74e-02 \\
$I=3$ & 1.70e-01 & 9.82e-02 & 6.45e-02 & 1.98e-01 \\
$I=4$ & 1.12e-02 & 8.29e-03 & 8.17e-02 & 1.78e-02 \\
$I=5$ & 4.87e-01 & 4.84e-01 & 4.93e-01 & 4.89e-01 \\
\hline
\end{tabular}

    \end{subtable}

    \begin{subtable}{0.99\linewidth}
    \caption{Vertical velocity component $v$.}
    \centering
    \begin{tabular}{|c|c|c|c|c|}
\hline
 & $N_o=2$ & $N_o=4$ & $N_o=4 / q=0.75$ & $N_r=1 / N_o=2$ \\
\hline
$I=1$ & 2.94e-01 & 1.55e-01 & 2.80e-02 & 3.02e-01 \\
$I=2$ & 7.18e-02 & 2.92e-01 & 1.33e-02 & 8.53e-02 \\
$I=3$ & 1.66e-01 & 9.75e-02 & 6.29e-02 & 1.93e-01 \\
$I=4$ & 1.09e-02 & 8.62e-03 & 7.97e-02 & 1.74e-02 \\
$I=5$ & 4.75e-01 & 4.70e-01 & 4.80e-01 & 4.76e-01 \\
\hline
\end{tabular}

    \end{subtable}
\end{table}

\begin{table}
    \caption{Case 2: Averaged Sobol' sensitivity indices.}
    \label{tab:Sobol-u_ii}
\begin{subtable}{0.99\linewidth}
\centering
    \caption{Horizontal velocity component $u$.}
    \begin{tabular}{|l|c|c|c|c|}
\hline
 & $N_o=2$ & $N_o=4$ & $N_o=4 / q=0.75$ & $N_r=1 / N_o=2$ \\
\hline
$I=1$ & 2.45e-01 & 1.01e-01 & 1.95e-02 & 2.12e-01 \\
$I=12$ & 1.43e-03 & 1.02e-02 & 4.42e-04 & 2.22e-03 \\
$I=123$ & 0.00e+00 & 1.39e-03 & 0.00e+00 & 4.85e-03 \\
$I=1234$ & 0.00e+00 & 8.67e-05 & 0.00e+00 & 4.33e-04 \\
$I=12345$ & 0.00e+00 & 0.00e+00 & 0.00e+00 & 7.08e-04 \\
$I=1235$ & 0.00e+00 & 9.38e-06 & 0.00e+00 & 1.88e-04 \\
$I=124$ & 0.00e+00 & 1.88e-04 & 0.00e+00 & 7.27e-04 \\
$I=1245$ & 0.00e+00 & 2.46e-08 & 0.00e+00 & 1.93e-04 \\
$I=125$ & 0.00e+00 & 8.44e-05 & 0.00e+00 & 1.62e-04 \\
$I=13$ & 5.15e-02 & 3.93e-02 & 0.00e+00 & 8.15e-02 \\
$I=134$ & 0.00e+00 & 1.13e-04 & 0.00e+00 & 3.02e-04 \\
$I=1345$ & 0.00e+00 & 4.70e-06 & 0.00e+00 & 8.46e-05 \\
$I=135$ & 0.00e+00 & 4.21e-05 & 0.00e+00 & 7.10e-05 \\
$I=14$ & 2.42e-03 & 2.20e-03 & 9.16e-04 & 4.83e-03 \\
$I=145$ & 0.00e+00 & 5.26e-05 & 0.00e+00 & 1.36e-04 \\
$I=15$ & 3.33e-04 & 2.69e-04 & 7.46e-03 & 3.89e-04 \\
$I=2$ & 5.84e-02 & 2.73e-01 & 8.04e-04 & 5.89e-02 \\
$I=23$ & 1.23e-02 & 1.24e-02 & 1.37e-05 & 1.61e-02 \\
$I=234$ & 0.00e+00 & 1.86e-04 & 0.00e+00 & 2.11e-04 \\
$I=2345$ & 0.00e+00 & 9.06e-07 & 0.00e+00 & 2.47e-04 \\
$I=235$ & 0.00e+00 & 3.96e-05 & 0.00e+00 & 2.15e-04 \\
$I=24$ & 1.26e-03 & 1.83e-03 & 1.24e-02 & 1.77e-03 \\
$I=245$ & 0.00e+00 & 1.06e-04 & 0.00e+00 & 1.47e-04 \\
$I=25$ & 1.12e-04 & 1.46e-03 & 0.00e+00 & 2.74e-04 \\
$I=3$ & 1.06e-01 & 4.44e-02 & 7.69e-04 & 9.27e-02 \\
$I=34$ & 2.69e-05 & 1.89e-05 & 6.32e-02 & 8.47e-05 \\
$I=345$ & 0.00e+00 & 9.13e-07 & 0.00e+00 & 4.67e-05 \\
$I=35$ & 1.02e-04 & 1.21e-04 & 4.86e-04 & 1.15e-04 \\
$I=4$ & 7.35e-03 & 3.41e-03 & 3.47e-03 & 7.71e-03 \\
$I=45$ & 1.12e-04 & 8.18e-05 & 1.65e-03 & 1.56e-04 \\
$I=5$ & 4.86e-01 & 4.82e-01 & 4.83e-01 & 4.85e-01 \\
\hline
\end{tabular}

    \end{subtable}
\end{table}

\begin{table}
    \caption{Case 2: Averaged Sobol' sensitivity indices.}
    \label{tab:Sobol-v_ii}
    \begin{subtable}{0.99\linewidth}
    \centering
    \caption{Vertical velocity component $v$.}
    \begin{tabular}{|l|c|c|c|c|}
\hline
 & $N_o=2$ & $N_o=4$ & $N_o=4 / q=0.75$ & $N_r=1 / N_o=2$ \\
\hline
$I=1$ & 2.39e-01 & 1.01e-01 & 1.92e-02 & 2.07e-01 \\
$I=12$ & 1.21e-03 & 1.01e-02 & 4.51e-04 & 2.11e-03 \\
$I=123$ & 0.00e+00 & 1.37e-03 & 0.00e+00 & 4.68e-03 \\
$I=1234$ & 0.00e+00 & 8.72e-05 & 0.00e+00 & 4.01e-04 \\
$I=12345$ & 0.00e+00 & 0.00e+00 & 0.00e+00 & 6.84e-04 \\
$I=1235$ & 0.00e+00 & 8.95e-06 & 0.00e+00 & 1.46e-04 \\
$I=124$ & 0.00e+00 & 1.91e-04 & 0.00e+00 & 7.05e-04 \\
$I=1245$ & 0.00e+00 & 1.90e-08 & 0.00e+00 & 1.54e-04 \\
$I=125$ & 0.00e+00 & 8.58e-05 & 0.00e+00 & 1.37e-04 \\
$I=13$ & 5.05e-02 & 3.89e-02 & 0.00e+00 & 7.98e-02 \\
$I=134$ & 0.00e+00 & 1.20e-04 & 0.00e+00 & 2.99e-04 \\
$I=1345$ & 0.00e+00 & 4.84e-06 & 0.00e+00 & 7.32e-05 \\
$I=135$ & 0.00e+00 & 3.21e-05 & 0.00e+00 & 6.49e-05 \\
$I=14$ & 2.43e-03 & 2.31e-03 & 9.35e-04 & 4.89e-03 \\
$I=145$ & 0.00e+00 & 5.44e-05 & 0.00e+00 & 1.36e-04 \\
$I=15$ & 5.30e-04 & 3.26e-04 & 7.40e-03 & 5.61e-04 \\
$I=2$ & 5.69e-02 & 2.64e-01 & 7.50e-04 & 5.76e-02 \\
$I=23$ & 1.22e-02 & 1.27e-02 & 9.88e-06 & 1.60e-02 \\
$I=234$ & 0.00e+00 & 1.95e-04 & 0.00e+00 & 1.90e-04 \\
$I=2345$ & 0.00e+00 & 9.36e-07 & 0.00e+00 & 2.03e-04 \\
$I=235$ & 0.00e+00 & 3.83e-05 & 0.00e+00 & 1.86e-04 \\
$I=24$ & 1.28e-03 & 1.88e-03 & 1.21e-02 & 1.76e-03 \\
$I=245$ & 0.00e+00 & 9.73e-05 & 0.00e+00 & 1.00e-04 \\
$I=25$ & 1.25e-04 & 1.05e-03 & 0.00e+00 & 3.64e-04 \\
$I=3$ & 1.04e-01 & 4.40e-02 & 7.56e-04 & 9.05e-02 \\
$I=34$ & 1.55e-05 & 1.60e-05 & 6.17e-02 & 6.14e-05 \\
$I=345$ & 0.00e+00 & 8.93e-07 & 0.00e+00 & 4.26e-05 \\
$I=35$ & 4.80e-05 & 4.86e-05 & 4.81e-04 & 5.47e-05 \\
$I=4$ & 7.06e-03 & 3.59e-03 & 3.35e-03 & 7.52e-03 \\
$I=45$ & 1.07e-04 & 7.27e-05 & 1.61e-03 & 1.47e-04 \\
$I=5$ & 4.74e-01 & 4.68e-01 & 4.70e-01 & 4.73e-01 \\
\hline
\end{tabular}

    \end{subtable}
\end{table}

\end{document}